\input amstex
\documentstyle{amsppt}
\pagewidth{5.4in}
\pageheight{7.6in}
\magnification=1200
\TagsOnRight
\NoRunningHeads
\topmatter
\title
\bf generalized $\Cal{L}$-geodesic and monotonicity\\ 
of the generalized reduced\\ 
volume in the Ricci flow  
\endtitle
\author
Shu-Yu Hsu
\endauthor
\affil
Department of Mathematics\\
National Chung Cheng University\\
168 University Road, Min-Hsiung\\
Chia-Yi 621, Taiwan, R.O.C.\\
e-mail:syhsu\@math.ccu.edu.tw
\endaffil
\date
March 3, 2008
\enddate
\address
e-mail address:syhsu\@math.ccu.edu.tw
\endaddress
\abstract
Suppose $M$ is a complete n-dimensional manifold, $n\ge 2$, with a metric 
$\overline{g}_{ij}(x,t)$ that evolves by the Ricci flow $\partial_t
\overline{g}_{ij}=-2\overline{R}_{ij}$ in $M\times (0,T)$. For any $0<p<1$, 
$(p_0,t_0)\in M\times (0,T)$, $q\in M$, we define the $\Cal{L}_p$-length 
between $p_0$ and $q$, $\Cal{L}_p$-geodesic, the generalized reduced 
distance $l_p$ and the generalized reduced volume $\widetilde{V}_p(\tau)$, 
$\tau=t_0-t$, 
corresponding to the $\Cal{L}_p$-geodesic at the point $p_0$ at time $t_0$. 
Under the condition $\overline{R}_{ij}\ge -c_1\overline{g}_{ij}$ on 
$M\times (0,t_0)$ for some constant $c_1>0$, we will prove the existence of 
a $\Cal{L}_p$-geodesic which minimize the $\Cal{L}_p(q,\overline{\tau})$-length
between $p_0$ and $q$ for any $\overline{\tau}>0$. This result for the case 
$p=1/2$ is conjectured and used many times but no proof of it was given in 
Perelman's papers on Ricci flow. Let $g(\tau)=\overline{g}(t_0-\tau)$ and 
let $\widetilde{V}_p^{\overline{\tau}}(\tau)$ be the rescaled generalized 
reduced volume. Suppose $M$ also has nonnegative curvature operator with 
respect to the metric $\overline{g}(t)$ for any $t\in (0,T)$ and when 
$1/2<p<1$, $M$ has uniformly bounded scalar curvature on $(0,T)$. 
Let $0<c<1$ and let $\tau_0=\min ((2(1-p))^{-1/(2p-1)},t_0)$. For any 
$1/2\le p<1$ we prove that there exists a constant $A_0\ge 0$ with $A_0=0$ 
for $p=1/2$ such that $e^{-A_0\tau}\widetilde{V}_p(\tau)$ is a monotone 
decreasing function in $(0,\overline{\tau}_1)$ where $\overline{\tau}_1
=(1-c)\tau_0$ if $1/2<p<1$ and $\overline{\tau}_1=t_0$ if $p=1/2$.
When $(M,\overline{g})$ is an ancient $\kappa$-solution of the Ricci 
flow, we will
prove a monotonicity property of the rescaled generalized volume
$\widetilde{V}_p^{\overline{\tau}}(\tau)$ with respect to $\overline{\tau}$
for any $1/2\le p<1$. When $p=1/2$, the 
$\Cal{L}_p$-length, $\Cal{L}_p$-geodesic, the $l_p$ function and 
$\widetilde{V}_p(\tau)$ are equal to the $\Cal{L}$-length, 
$\Cal{L}$-geodesic, the reduced distance $l$ and the reduced volume 
$\widetilde{V}(\tau)$ introduced by Perelman in his papers on Ricci flow. 
We will also prove a conjecture on the reduced distance $l$ and the reduced
volume $\widetilde{V}$ which was used by Perelman without proof in \cite{P1}. 
\endabstract
\keywords
Ricci flow, $\Cal{L}_p$-length, $\Cal{L}_p$-geodesic, existence of minimizer,
$\Cal{L}_p$-cut point, generalized reduced distance, generalized reduced 
volume, monotonicity of generalized reduced 
volume
\endkeywords
\subjclass
Primary 58J35, 53C44 Secondary 58C99
\endsubjclass
\endtopmatter
\NoBlackBoxes
\define \pd#1#2{\frac{\partial #1}{\partial #2}}
\define \1{\partial}
\define \2{\overline}
\define \3{\varepsilon}
\define \4{\widetilde}
\define \5{\underline}
\document

Recently there is a lot of study on the Ricci flow on manifold by R.~Hamilton
\cite{H1--6}, S.Y.~Hsu \cite{Hs1--5}, G.~Perelman \cite{P1}, \cite{P2}, 
W.X.~Shi \cite{S1}, \cite{S2}, L.F.~Wu \cite{W1}, \cite{W2}, and others.
We refer the readers to the lecture notes by B.~Chow \cite{Ch} and the book
\cite{CK} by B.~Chow and D.~Knopf on the basics of Ricci flow and the papers
\cite{P1}, \cite{P2} of G.~Perelman for the most recent results on Ricci flow.
 
In the paper \cite{H1} R.~Hamilton proved that if $M$ is a compact manifold 
with a metric $\2{g}$ that evolves by the Ricci flow
$$ 
\frac{\1 }{\1 t}\2{g}_{ij}=-2\2{R}_{ij}\tag 0.1
$$
where $\2{R}_{ij}$ is the Ricci curvature of $\2{g}$ and $\2{g}_{ij}(x,0)
=\2{g}_{ij}(x)$ is a metric of strictly positive Ricci curvature,
then the evolving metric will converge modulo scaling to a metric 
of constant positive curvature. Similiar result was obtained by R.~Hamilton
\cite{H2} for compact 4-dimensional manifolds with positive curvature 
operator. Harnack inequality for the Ricci flow was proved by R.~Hamilton 
in \cite{H4}.

Short time existence of solutions of the Ricci flow on complete
non-compact Riemannian manifold with bounded curvature was proved by W.X~Shi
\cite{S1}. Global existence and uniqueness of solutions of the Ricci flow 
on non-compact manifold $\Bbb{R}^2$ was obtained by S.Y.~Hsu in \cite{Hs1}. 
Asymptotic behaviour of solutions of the Ricci flow equation on $\Bbb{R}^2$ 
was proved by by S.Y.~Hsu in \cite{Hs2}, \cite{Hs3}, \cite{Hs4}.

In \cite{P1}, \cite{P2}, G.~Perelman introduced the concept of 
$\Cal{L}$-length, $\Cal{L}$-geodesic, the reduced distance $l$ and the 
reduced volume $\4{V}(\tau)$
for Ricci flow on complete manifolds with postive bounded curvature operator. 
G.~Perelman found that these are very useful tools in studying Ricci flow 
on manifolds.  He used these tools to proved many new properties for the 
Ricci flow in \cite{P1}, \cite{P2}. These included the non-local collapsing 
theorem and the asymptotic convergence of a subsequence of the rescaled 
solution of an ancient $\kappa$-solution to a soliton solution of the Ricci 
flow on complete manifold with postive bounded curvature operator. Recently 
R.~Ye (\cite{Ye1}, \cite{Ye2}) 
extended the concept of $\Cal{L}$-geodesic, the reduced distance $l$ and the 
reduced volume $\4{V}(\tau)$ to manifolds with a lower bound on the Ricci 
curvature.

In this paper we will generalize the notion of $\Cal{L}$-length, 
$\Cal{L}$-geodesic, the reduced distance $l$ and the reduced volume 
$\4{V}(\tau)$ of G.~Perelman. For any $0<p<1$, $\tau>0$, we will define the 
$\Cal{L}_p(q,\tau)$-length, $L_p(q,\tau)$-length, $\Cal{L}_p$-exponential map,
$\Cal{L}_p$-geodesic, and reduced volume $\4{V}_p(\tau)$ and prove various 
properties 
of them in this paper. When $p=1/2$, $\Cal{L}_p(q,\tau)$, $L_p(q,\tau)$, 
$\Cal{L}_p$-exponential map, $\4{V}_p(\tau)$ are equal to the $\Cal{L}
(q,\tau)$, $L(q,\tau)$, $\Cal{L}$-exponential map, and $\4{V}(\tau)$ defined 
by G.~Perelman in \cite{P1}. 

For any $q\in M$, $\2{\tau}>0$, we will prove the existence of a 
$\Cal{L}_p$-geodesic which minimize the $\Cal{L}_p(q,\2{\tau})$-length. This 
result for the case $p=1/2$ is conjectured and used many times in 
G.~Perelman's paper on Ricci flow \cite{P1}, \cite{P2}, but no proof of it 
is given in his papers. There is also no detail proof of this important 
conjecture
in the recent book of J.W.~Morgan and G.~Tian \cite{MT} and the paper of
H.D.~Cao and X.P.~Zhu \cite{CZ} on Ricci flow. My result is new and
answers in affirmative the existence of such $\Cal{L}$-geodesic minimizer 
for the $L_p(q,\tau)$-length which is crucial to the proof of many results
in \cite{P1}, \cite{P2}. We also prove that for any 
$\Cal{L}_p(q,\2{\tau})$-length minimizing $\Cal{L}_p$-geodesic there does 
not exists any $\Cal{L}_p$-conjugate points along the curve.

One remarkable property of the reduced volume $\4{V}(\tau)$, $\tau=t_0-t$, 
with respect to any point $(p_0,t_0)\in M\times (0,T)$ proved by 
G.~Perelman is that it is a monotone decreasing function of $\tau\in 
(0,t_0)$. 
Surprisingly in this paper we find that the generalized reduced volume 
$\4{V}_p(\tau)$ with respect to any point $(p_0,t_0)\in M\times (0,T)$
also has similar monotonicity property. Suppose $M$ is complete and has 
nonnegative curvature operator with respect to the metric $\2{g}(t)$ for any 
$t\in (0,T)$ and when $1/2<p<1$, $M$ has uniformly bounded scalar curvature 
on $(0,T)$ . Let $0<c<1$ and let
$$
\tau_0=\min ((2(1-p))^{-\frac{1}{2p-1}},t_0).\tag 0.2
$$
For any $1/2\le p<1$ we prove that there exists a constant $A_0\ge 0$ with 
$A_0=0$ for $p=1/2$ such that $e^{-A_0\tau}\4{V}_p(\tau)$ is a monotone 
decreasing function in $(0,\2{\tau}_1)$ where $\2{\tau}_1=(1-c)\tau_0$ if 
$1/2<p<1$ and $\2{\tau}_1=t_0$ if $p=1/2$.

Suppose $(M,\2{g})$ is an ancient $\kappa$-solution of the Ricci flow
and $g(\tau)=\2{g}(t_0-\tau)$ for some constant $t_0<0$. 
Let $\2{\tau}_0>0$ for $1/2<p<1$ and $\2{\tau}_0=0$ for $p=1/2$.
When $1/2<p<1$, suppose also that 
$(M,g(\tau))$ is compact. Let $\4{V}_p^{\2{\tau}}(\rho)$ be the rescaled 
generalized volume. Then for any $1/2\le p<1$ there exist constants 
$A_0\ge 0$, $A_1\ge 0$, $A_2\ge 0$, $c_2>0$, such that 
$e^{-W(\2{\tau},\rho)}\4{V}_p^{\2{\tau}}(\rho)$ is a monotone decreasing 
function of $\2{\tau}>\2{\tau}_0$ for any $\rho$ satisfying 
$$
0<\rho\le\left (\frac{1}{2(1-p)}\right)^{\frac{1}{2p-1}}\tag 0.3
$$
where $W(\2{\tau},\rho)=(A_0\rho+A_1\rho^{2p}+A_2\rho^{2p-3}e^{2c_2\2{\tau}
\rho})\2{\tau}.
$
Moreover
$$
\lim_{\2{\tau}\to 0^+}\4{V}_p^{\2{\tau}}(\rho)=\left(\frac{\sqrt{\pi}}
{1-p}\right)^n\rho^{\frac{(1-p)n}{2}}.\tag 0.4
$$
Note that when $p=1/2$, one can take $A_0=A_1=A_2=0$ 
and the result reduces to Perelman's monotonicity property for ancient 
$\kappa$-solution of the Ricci flow \cite{P1}. 

When $(M,\2{g})$ is an ancient $\kappa$-solution of the Ricci flow in 
$(-\infty,0)$ with uniformly bounded nonnegative curvature operator, 
then for any $t_0<0$, $p_0\in M$, $0<p<1$, $\tau_2>\tau_1>0$ 
we prove the existence of $\{q_i\}_{i=\infty}^{\infty}\subset M$ and
$\{\2{\tau}_i\}_{i=1}^{\infty}$, $\2{\tau}_i\to\infty$ as $i\to\infty$,
such that the rescaled $l_p$ function $l_p^{\2{\tau}_i}(q,\tau)$ converges 
uniformly on $B_0(q_i,r)\times [\tau_1,\tau_2]$ as $i\to\infty$ for any
$r>0$ where $B_0(q_i,r)$ is the geodesic ball of radius $r>0$ with respect
to the metric $\2{g}(t_0)$. 

We will also prove a conjecture on the reduced distance $l$ and the reduced
volume $\4{V}(\tau)$ used by Perelman without proof in \cite{P1}. Suppose 
$(M,\2{g})$ is an ancient $\kappa$-solution of the Ricci flow with uniformly 
bounded nonnegative curvature operator such that $\2{g}(t)$ is not a flat 
metric for any $t<0$. If $\4{V}(\tau_1)=\4{V}(\tau_2)$ for some $\tau_2>\tau_1
>0$, we prove that the reduced distance $l\in C^{\infty}(M\times [\tau_1,
\tau_2])$ and $g(\tau)=\2{g}(t_0-\tau)$ is a shrinking soliton in $M\times 
[\tau_1,\tau_2]$. 

The plan of the paper is as follows. In section 1 we will use a 
modification of the technique of \cite{P1} to prove the first variation 
formula for the $\Cal{L}_p(q,\2{\tau})$-length. We will also prove the 
existence of a $\Cal{L}_p(q,\2{\tau})$-geodesic minimizer for the 
$L_p(q,\2{\tau})$-length. 
We will prove various properties of the $\Cal{L}_p$-exponential map and 
$\Cal{L}_p$ cut locus in section 2. In section 3 we will prove the second 
variation formula for the $L_p(q,\tau)$-length. We will prove various 
properties of the $L_p(q,\tau)$-length, the generalized reduced distance 
$l_p$ and the generalized reduced volume $\4{V}_p(\tau)$. In section 4 we will 
prove the monotonicity property the generalized reduced volume 
$\4{V}_p(\tau)$. In section 5 we will prove the monotonicity property of 
the rescaled generalized reduced volume $\4{V}_p^{\2{\tau}}(\tau)$ with 
respect to $\2{\tau}$. In section 6 we will prove a conjecture on the 
the reduced distance and the reduced volume used by Perelman without proof 
in \cite{P1}.

We first start with a definition. Let $(M,\2{g})$ be a Riemannian manifold 
with the metric $\2{g}$ evolving by the Ricci flow (0.1) in $M\times (0,T)$.
Let $(p_0,t_0)\in M\times (0,T)$. For any $0<t<t_0$, let 
$\tau=t_0-t$ and 
$$
g(\tau)=\2{g}(t_0-\tau).\tag 0.5
$$
Let $R(q,\tau)$, $R_{ij}(q,\tau)$, $R(X_1,X_2)X_3(q,\tau)$ and $Rm(q,\tau)$ 
be the scalar curvature, Ricci curvature, curvature and Riemannian curvature 
of $M$ at $q$ with respect to the metric $g(\tau)$ and $X_1,X_2,X_3\in T_qM$. 
For any $0<p<1$, $p_0,q\in M$, $\2{\tau}\in (0,t_0)$, and piecewise 
differentiable
curve $\gamma : [0,\2{\tau}]\to M$ joining $p_0$ and $q$ with $\gamma (0)=p_0$
and $\gamma (\2{\tau})=q$, we define the $\Cal{L}_p^{p_0}(q,\2{\tau})$-length of  
the curve $\gamma$ between $p_0$ and $q$ by
$$
\Cal{L}_p^{p_0}(q,\gamma ,\2{\tau})=\int_0^{\2{\tau}}\tau^p(R(\gamma (\tau),\tau)+
|\gamma '(\tau)|^2)\,d\tau
$$
where $|\gamma '(\tau)|=|\gamma '(\tau)|_{g(\tau)}$. Let $\Cal{F}^{p_0}
(q,\2{\tau})$
be the family of all piecewise differentiable curves $\gamma : [0,\2{\tau}]
\to M$ satisfying $\gamma (0)=p_0$ and $\gamma (\2{\tau})=q$, 
$$
L_p^{p_0}(q,\2{\tau})=\inf_{\gamma\in\Cal{F}^{p_0}(q,\2{\tau})}
\Cal{L}_p^{p_0}(q,\gamma ,\2{\tau}),
$$
and let
$$
l_p^{p_0}(q,\tau)=(1-p)\frac{L_p^{p_0}(q,\tau)}{\tau^{1-p}}\tag 0.6
$$
be the generalized reduced distance. Let
$$
\4{V}_p^{p_0}(\tau)=\int_M\tau^{-(1-p)n}e^{-l_p^{p_0}(q,\tau)}dV_{g(\tau)}
$$
be the generalized reduced volume corresponding to the $\Cal{L}_p^{p_0}
(\cdot,\tau)$-length with respect to $(p_0,t_0)$. Then $l^{p_0}(q,\tau)
=l_{\frac{1}{2}}^{p_0}(q,\tau)$ and $\4{V}^{p_0}(\tau)
=\4{V}_{\frac{1}{2}}^{p_0}(\tau)$ are the reduced length and reduced volume of
Perelman \cite{P1}. Let $L^{p_0}(q,\tau)=L_{\frac{1}{2}}^{p_0}(q,\tau)$. 
When there is no ambiguity, we will drop 
the superscript $p_0$.

Let $q_0\in M$ and $0<\tau_0<t_0$. For any $0<p<1$, $q\in M$, 
$\2{\tau}\in (\tau_0,t_0)$, and piecewise differentiable
curve $\gamma : [\tau_0,\2{\tau}]\to M$ joining $q_0$ and $q$ with 
$\gamma (\tau_0)=q_0$ and $\gamma (\2{\tau})=q$, we define the 
$\Cal{L}_{\tau_0,p}^{q_0}(q,\2{\tau})$-length of  
the curve $\gamma$ between $q_0$ and $q$ by
$$
\Cal{L}_{\tau_0,p}^{q_0}(q,\gamma ,\2{\tau})
=\int_{\tau_0}^{\2{\tau}}\tau^p(R(\gamma (\tau),\tau)+
|\gamma '(\tau)|^2)\,d\tau.
$$
Let $\Cal{F}_{\tau_0}^{q_0}(q,\2{\tau})$ be the family of all piecewise 
differentiable curves $\gamma : [\tau_0,\2{\tau}]\to M$ satisfying 
$\gamma (\tau_0)=q_0$ and $\gamma (\2{\tau})=q$ and let
$$
L_{\tau_0,p}^{q_0}(q,\2{\tau})=\inf_{\gamma\in\Cal{F}_{\tau_0}^{q_0}
(q,\2{\tau})}\Cal{L}_{\tau_0,p}^{q_0}(q,\gamma ,\2{\tau}).
$$ 
For any $r>0$, $0\le\tau\le t_0$, we let $B_{\tau}(q,r)$ be the 
geodesic ball of radius $r$ in $M$ around the point $q$ with respect to the 
metric $g(\tau)$. For any $v\in T_{p_0}M$, we let 
$$
\Cal{B}(v,r)=\{v'\in T_{p_0}M:|v-v'|_{g(p_0,0)}<r\}.
$$
We also let $d_{\tau}(q_1,q_2)=d_{g(\tau)}(q_1,q_2)$ be the distance 
between $q_1$ and $q_2$ with respect to the metric $g(\tau)$. For any 
$0<\tau<t_0$ and measurable set $E\subset M$ with respect to the metric 
$g(\tau)$, we let $m_{\tau}(E)$ be the measure of $E$ with respect to the 
metric $g(\tau)$. We let $dV_{g(\tau)}(q)=\sqrt{g(q,\tau)}dq$ be the volume 
form of the metric $g(\tau)$.

Let $\kappa>0$. A Ricci flow $(M,\2{g})$ is said to be $\kappa$-noncollapsing 
at the point $(q',t')$ on the scale $r_0>0$ \cite{P1} if $\forall 0<r\le r_0$,
$$
\text{Vol}_{\2{g}(t')}(B_{\2{g}(t')}(q',r))\ge\kappa r^n
$$
holds whenever
$$
|\2{\text{Rm}}|(q,t)\le r^{-2}\quad\forall d_{\2{g}(t')}(q',q)<r, t'-r^2\le t\le t'
$$
holds where $B_{\2{g}(t')}(q',r)$ is the geodesic ball
of radius $r$ in $M$ around the point $q'$ with respect to the metric 
$\2{g}(t')$. A Ricci flow $(M,\2{g})$  is said to be an ancient 
$\kappa$-solution if it is a solution of the Ricci flow in 
$M\times (-\infty,0]$ such that for each $t\le 0$ the metric $\2{g}(t)$ is
not a flat metric, $(M,\2{g}(t))$ is a complete manifold of nonnegative 
and uniformly bounded curvature, and $(M,\2{g}(t))$ is $\kappa$-noncollapsing 
on all scales at all points of $M\times (-\infty,0]$. 

We will assume that $M$ is complete with respect to $\2{g}(t)$ for any 
$0<t<T$ for the rest of the paper. Unless stated otherwise we will fix
the point $(p_0,t_0)$ and consider the $\Cal{L}_p(q,\tau)$, $l_p(q,\tau)$,
etc. all with respect to this fixed reference point. We also associate the
product manifold $M\times (0,t_0)$ with the product metric $g\,dx^2\oplus 
d\tau^2$.

$$
\text{Section 1}
$$

In this section we will use the technique of \cite{P1} to prove the first 
variation formula for $\Cal{L}_p(q,\gamma,\2{\tau})$ for any curve $\gamma:
[0,\2{\tau}]\to M$ joining $p_0$ and $q$ with $\gamma (0)=p_0$ and $\gamma 
(\2{\tau})=q$. We will prove the non-trivial fact that the $L_p(q,\2{\tau})$ 
length can be realized by some $\Cal{L}_p$-geodesic on $M$. We will let 
$<\cdot,\cdot>_{g(\tau)}$ be the inner product with respect to the metric 
$g(\tau)$. When there is no ambiguity, we will write $<\cdot,\cdot>$ for 
$<\cdot,\cdot>_{g(\tau)}$. 

\proclaim{\bf Lemma 1.1}
Let $\gamma\in\Cal{F}(q,\2{\tau})$ and let $Y:[0,\2{\tau}]\to TM$ be a 
vector field along $\gamma$ with $Y(0)=0$. Suppose $\gamma$ is differentiable 
on $[0,\2{\tau}]$. Then
$$
\delta_Y\Cal{L}_p(q,\gamma ,\2{\tau})=2\2{\tau}^p<X(\2{\tau}),Y(\2{\tau})>
+\int_0^{\2{\tau}}\tau^p<Y,\nabla R-\frac{2p}{\tau}X-2\nabla_XX-4\text{Ric}
(X,\cdot)>\,d\tau\tag 1.1
$$
where $X=X(\tau)=\gamma'(\tau)$ and the inner product in the integral is 
evaluated at $\tau$.
\endproclaim
\demo{Proof}
Let $f:[0,\2{\tau}]\times (-\3,\3)\to M$ be a variation of $\gamma$ 
with respect to the vector field $Y$ such that $f(0,z)=p_0$ for all $z\in
(-\3,\3)$. Then
$$\align
&\frac{d}{dz}\Cal{L}_p(f(\2{\tau},z),f(\cdot,z) ,\2{\tau})\\
=&\frac{d}{dz}\int_0^{\2{\tau}}\tau^p(R(f(\tau,z),\tau)
+|\nabla_{\tau}f|^2)\,d\tau\\ 
=&\int_0^{\2{\tau}}\tau^p(<\nabla_zf,\nabla R>+2<\nabla_{\tau}f,\nabla_z
\nabla_{\tau}f>)\,d\tau\\
=&\int_0^{\2{\tau}}\tau^p(<\nabla_zf,\nabla R>+2<\nabla_{\tau}f,\nabla_\tau
\nabla_zf>)\,d\tau.\tag 1.2\\
=&\int_0^{\2{\tau}}\tau^p(<\nabla_zf,\nabla R>+2\frac{d}{d\tau}
<\nabla_{\tau}f,\nabla_zf>-2<\nabla_{\tau}\nabla_{\tau}f,\nabla_zf>
-4\text{Ric}(\nabla_{\tau}f,\nabla_zf))\,d\tau\\
=&2\2{\tau}^p<\nabla_{\tau}f(\2{\tau},z),\nabla_zf(\2{\tau},z)>\\
&\qquad +\int_0^{\2{\tau}}\tau^p
<\nabla_zf,\nabla R-\frac{2p}{\tau}\nabla_{\tau}f-2\nabla_{\tau}
\nabla_{\tau}f-4\text{Ric}(\nabla_{\tau}f,\cdot)>\,d\tau.\tag 1.3
\endalign
$$
By putting $z=0$ in (1.3), (1.1) follows.
\enddemo

Let $s_0=t_0^{1-p}$. For any $0\le\2{s}<s_0$, $p_0\in M$, 
$\4{\gamma}\in\Cal{F}^{p_0}
(q,\2{s})$, $0\le s\le\2{s}$, let $\4{R}(q,s)=R(q,s^{\frac{1}{1-p}})$,
$$
\4{\Cal{L}}_p^{p_0}(q,\4{\gamma},\2{s})=\frac{1}{1-p}\int_0^{\2{s}}
(s^{\frac{2p}{1-p}}\4{R}(\4{\gamma}(s),s)+(1-p)^2|\4{\gamma}'(s)|^2)\,ds
$$
where $|\4{\gamma}'(s)|=|\4{\gamma}'(s)|_{g(s^{1/(1-p)})}$ and let
$$
\4{L}_p^{p_0}(q,\2{s})=\inf_{\4{\gamma}\in\Cal{F}^{p_0}(q,\2{s})}
\4{\Cal{L}}_p^{p_0}(q,\4{\gamma},\2{s}).
$$
Then by direct computation,
$$
\Cal{L}_p^{p_0}(q,\gamma,\2{\tau})=\4{\Cal{L}}_p^{p_0}(q,\4{\gamma},
\2{\tau}^{1-p})\quad\forall\gamma\in\Cal{F}^{p_0}(q,\2{\tau})\tag 1.4
$$
where
$$
\4{\gamma}(s)=\gamma (\tau),s=\tau^{1-p}.\tag 1.5
$$
Hence
$$
L_p^{p_0}(q,\2{\tau})=\4{L}_p^{p_0}(q,\2{\tau}^{1-p}).\tag 1.6
$$
We will now let $\4{g}(q,s)=g(q,\tau)$, $\4{R}(q,s)=R(q,\tau)$, 
$\4{\text{Ric}}(q,s)=\text{Ric}(q,\tau)$, and $\4{\Gamma}_{ij}^r(q,s)
=\Gamma_{ij}^r(q,s)$ where $s=\tau^{1-p}$ for the rest 
of the paper. When there is no ambiguity, we will drop the superscript $p_0$.

\proclaim{\bf Lemma 1.2}
Let $\4{\gamma}\in\Cal{F}(q,\2{s})$ and let $\4{Y}:[0,\2{s}]\to TM$ be a 
vector field along $\4{\gamma}$ with $\4{Y}(0)=0$. Suppose $\4{\gamma}$ is 
differentiable on $[0,\2{s}]$. Then
$$\align
&\delta_{\4{Y}}\4{\Cal{L}}_p(q,\4{\gamma},\2{s})\\
=&\frac{1}{1-p}\int_0^{\2{s}}(s^{\frac{2p}{1-p}}\4{Y}(\4{R})
+2(1-p)^2<\4{X},\nabla_{\4{X}}\4{Y}>)\,ds\\
=&2(1-p)<\4{X}(\2{s}),\4{Y}(\2{s})>\\
&\qquad +\frac{1}{1-p}\int_0^{\2{s}}<\4{Y},s^{\frac{2p}{1-p}}\nabla\4{R}
-2(1-p)^2\nabla_{\4{X}}\4{X}-4(1-p)s^{\frac{p}{1-p}}
\4{\text{Ric}}(\4{X},\cdot)>\,ds\tag 1.7\endalign
$$
where $\4{X}(s)=\4{\gamma}'(s)$, $\4{R}(s)=\4{R}(\4{\gamma}(s),s)$, 
$\4{\text{Ric}}(s)=\4{\text{Ric}}(\4{\gamma}(s),s)$ and $\nabla
=\nabla^{g(\tau)}$ with $s=\tau^{1-p}$.
\endproclaim
\demo{Proof}
Let $\4{f}:[0,\2{s}]\times (-\3,\3)\to M$ be a variation of $\4{\gamma}$ 
with respect to the vector field $\4{Y}$ such that $\4{f}(0,z)=p_0$ for all 
$z\in (-\3,\3)$. Since
$$
\frac{d}{ds}\4{g}_{ij}(q,s)
=\frac{d}{d\tau}g_{ij}(q,\tau)\cdot\frac{d\tau}{ds}=\frac{2}{1-p}
s^{\frac{p}{1-p}}\4{R}_{ij}(q,s)\tag 1.8
$$
where $s=\tau^{1-p}$, we have
$$\align
&\frac{d}{dz}\4{\Cal{L}}_p(\4{f}(\2{s},z),\4{f}(\cdot,z),\2{s})\\
=&\frac{1}{1-p}\int_0^{\2{s}}(s^{\frac{2p}{1-p}}<\nabla_z\4{f},\nabla\4{R}>
+2(1-p)^2<\nabla_s\4{f},\nabla_z\nabla_s\4{f}>)\,ds\\
=&\frac{1}{1-p}\int_0^{\2{s}}(s^{\frac{2p}{1-p}}<\nabla_z\4{f},\nabla\4{R}>
+2(1-p)^2<\nabla_s\4{f},\nabla_s\nabla_z\4{f}>)\,ds\\
=&\frac{1}{1-p}\int_0^{\2{s}}\biggl\{s^{\frac{2p}{1-p}}
<\nabla_z\4{f},\nabla\4{R}>
+2(1-p)^2\biggl [\frac{d}{ds}<\nabla_s\4{f},\nabla_z\4{f}>
-<\nabla_s\nabla_s\4{f},\nabla_z\4{f}>\\
&\qquad -\frac{2}{1-p}s^{\frac{p}{1-p}}
\4{\text{Ric}}(\nabla_s\4{f},\nabla_z\4{f})\biggr ]\biggr\}\,ds\\
=&2(1-p)<\nabla_s\4{f}(\2{s},z),\nabla_z\4{f}(\2{s},z)>\\
&\qquad +\frac{1}{1-p}\int_0^{\2{s}}<\nabla_z\4{f},s^{\frac{2p}{1-p}}
\nabla\4{R}-2(1-p)^2\nabla_s\nabla_s\4{f}-4(1-p)s^{\frac{p}{1-p}}
\4{\text{Ric}}(\nabla_s\4{f},\cdot)>\,ds\tag 1.9 
\endalign
$$
Putting $z=0$ in (1.9) we get (1.7) and the lemma follows.
\enddemo

From Lemma 1.1 and Lemma 1.2 it is natural to define the following. We say 
that a curve $\4{\gamma}\in\Cal{F}(q,\2{s})$ is a 
$\4{\Cal{L}}_p$-geodesic at $s\in (0,\2{s})$ if
it satisfies
$$ 
\nabla_{\4{X}}\4{X}-\frac{1}{2(1-p)^2}s^{\frac{2p}{1-p}}\nabla\4{R}
+\frac{2}{1-p}s^{\frac{p}{1-p}}\4{\text{Ric}}(\4{X},\cdot)=0\tag 1.10
$$
at $s$ where $\4{X}(s)=\4{\gamma}'(s)$ and $\nabla =\nabla^{g(\tau)}$ with
$\tau=s^{\frac{1}{1-p}}$. We say that it is a $\4{\Cal{L}}_p$-geodesic in 
$(0,\2{s})$ if it satisfies (1.10) in $(0,\2{s})$. Similarly we say a curve 
$\gamma\in\Cal{F}(q,\2{\tau})$ is a $\Cal{L}_p$-geodesic at $\tau\in 
(0,\2{\tau})$ if it satisfies
$$
\nabla_XX-\frac{1}{2}\nabla R+\frac{p}{\tau}X+2\text{Ric}(X,\cdot)=0\tag 1.11
$$
at $\tau$ where $X(\tau)=\gamma '(\tau)$. We say that $\gamma$ is a 
$\Cal{L}_p$-geodesic in $(0,\2{\tau})$ if it satisfies (1.11) in 
$(0,\2{\tau})$. Note that when $p=1/2$, the $\Cal{L}_p$-geodesic is equal to 
the $\Cal{L}$-geodesic defined by Perelman in \cite{P1}.
 
\proclaim{\bf Remark 1.3}
By direct computation $\gamma\in\Cal{F}(q,\2{\tau})$ is a 
$\Cal{L}_p$-geodesic at $\tau\in (0,\2{\tau})$ if and only if $\4{\gamma}
\in\Cal{F}(q,\2{s})$ is a $\4{\Cal{L}}_p$-geodesic at $s\in (0,\2{s})$ where 
$\gamma$, $\4{\gamma}$, $s$ and $\tau$ are related by (1.5) and $\2{s}
=\2{\tau}^{1-p}$. Moreover 
$$
\4{\gamma}'(0)=\frac{1}{1-p}\lim_{\tau\to 0}\tau^p\gamma'(\tau).\tag 1.12
$$
\endproclaim

\proclaim{\bf Lemma 1.4}
For any $\4{v}\in T_{p_0}M$, there exists a unique solution $\4{\gamma}(s)
=\4{\gamma}_{\4{v}}(s)=\4{\gamma}(s;\4{v})$ of (1.10) in $(0,s_0)$ with 
$$\left\{\aligned
&\4{\gamma}(0)=p_0\\
&\4{\gamma}'(0)=\4{v}\endaligned\right.\tag 1.13
$$
for some constant $s_0\in (0,t_0^{1-p}]$ where $(0,s_0)$ is the maximal
interval of existence of the solution. If $s_0<t_0^{1-p}$, then
$$
\lim_{s\to s_0^-}d_0(p_0,\4{\gamma}(s))=\infty.\tag 1.14
$$
If the Ricci curvature of $M$ is uniformly bounded on $(0,t_0]$, then 
$s_0=t_0^{1-p}$. 
\endproclaim
\demo{Proof}
Uniqueness of solution of (1.10) satisfying (1.13) follows by standard O.D.E. 
theory. Hence we only need to prove existence of solution of (1.10) satisfying
(1.13). We will use a continuity argument similar to that of section 17 of 
\cite{KL} to prove the existence of solution of (1.10) satisfying (1.13). We 
first observe that by standard O.D.E. theory there exists a constant $s_0'
\in (0,t_0^{1-p})$ such that (1.10), (1.13), has a unique 
solution $\4{\gamma}(s)$ in $(0,s_0')$. Let $(0,s_0)$ be the maximum interval 
of existence of solution $\4{\gamma}(s)$ of (1.10) and (1.13). Then $s_0\le 
t_0^{1-p}$. If $s_0=t_0^{1-p}$, we are done. So we suppose that
$s_0<t_0^{1-p}$. We claim that (1.14) holds. Suppose not. Then
there exist constants $s_1\in (0,s_0)$ and $C_1>0$ such that
$$
d_0(p_0,\4{\gamma}(s))\le C_1\quad\forall s_1\le s<s_0.\tag 1.15
$$
Let 
$$
r_0=\sup_{0\le s<s_0}d_0(p_0,\4{\gamma}(s)).
$$
By (1.15) $r_0<\infty$. Since $\2{B_0(p_0,r_0)}\times [0,s_0^{1/(1-p)}]$ is 
compact in $(q,\tau)\in M\times [0,t_0)$ when $M$ is equipped with the 
metric $g(0)$, there exists a constant $K_1>0$ such that
$$
|R|+|\nabla R|+|\text{Ric}|\le K_1\tag 1.16
$$
on $(q,\tau)\in\2{B_0(p_0,r_0)}\times [0,s_0^{1/(1-p)}]$.
Then by (1.10) and (1.16),
$$\align
\biggl |\frac{d}{ds}|\4{X}|^2\biggr |=&\biggl |2<\4{X},\nabla_{\4{X}}\4{X}>
+\frac{2}{1-p}s^{\frac{p}{1-p}}\4{\text{Ric}}(\4{X},\4{X})\biggr |\\
=&\biggl |<\4{X},\frac{1}{(1-p)^2}s^{\frac{2p}{1-p}}\nabla\4{R}
-\frac{4}{1-p}s^{\frac{p}{1-p}}\4{\text{Ric}}(\4{X},\cdot)>
+\frac{2}{1-p}s^{\frac{p}{1-p}}\4{\text{Ric}}(\4{X},\4{X})\biggr |\\
\le&\left |\frac{1}{(1-p)^2}s^{\frac{2p}{1-p}}<\4{X},\nabla\4{R}>\right |
+\frac{2}{1-p}s^{\frac{p}{1-p}}\left|\4{\text{Ric}}(\4{X},\4{X})\right|\\
\le&A_1K_1(s^{\frac{2p}{1-p}}|\4{X}|+s^{\frac{p}{1-p}}|\4{X}|^2)\\
\le&C_1(2s^{\frac{p}{1-p}}|\4{X}|^2+(s^{\frac{3p}{1-p}}/4))
\quad\forall 0\le s\le s_0.
\endalign
$$
where $A_1=\max ((1-p)^{-2},2(1-p)^{-1})$ and $C_1=A_1K_1$.
Hence $\forall 0\le s\le s_0$,
$$\align
&\left\{\aligned
&\frac{d}{ds}\biggl (e^{-2C_1(1-p)s^{\frac{1}{1-p}}}|\4{X}|^2\biggr )
\le (C_1/4)e^{-2C_1(1-p)s^{\frac{1}{1-p}}}s^{\frac{3p}{1-p}}\\
&\frac{d}{ds}\biggl (e^{2C_1(1-p)s^{\frac{1}{1-p}}}|\4{X}|^2\biggr )
\ge -(C_1/4)e^{2C_1(1-p)s^{\frac{1}{1-p}}}s^{\frac{3p}{1-p}}
\endaligned\right.\\
\Rightarrow\quad&e^{-C_2s^{\frac{1}{1-p}}}
(|\4{v}|^2-C_3e^{C_2s^{\frac{1}{1-p}}}s^{\frac{1+2p}{1-p}})\le |\4{X}(s)|^2
\le e^{C_2s^{\frac{1}{1-p}}}(|\4{v}|^2+C_3s^{\frac{1+2p}{1-p}})
\tag 1.17\endalign
$$
where $C_2=2(1-p)C_1$ and $C_3=C_2/(8(1+2p))$. By (1.17) and standard 
O.D.E. 
theory there exists a constant $\3_0\in (0,t_0^{1-p}-s_0]$ such that we can 
extend $\4{\gamma}$ to a solution of (1.10), (1.13), on $(0,s_0+\3_0)$. This 
contradicts the maximality of $s_0$. Thus (1.14) holds.

If the Ricci curvature of $M$ is uniformly bounded on $M\times (0,t_0]$, then
by the local estimates for the solutions of Ricci flow \cite{S1} and a similar 
argument as before we will get a contradiction if $s_0<t_0^{1-p}$. Hence 
$s_0=t_0^{1-p}$ and the lemma follows.
\enddemo

By Remark 1.3, Lemma 1.4, and (1.12), we have

\proclaim{\bf Corollary 1.5}
For any $v\in T_{p_0}M$, there exists a unique solution $\gamma (\tau)
=\gamma_v(\tau)=\gamma (\tau ;v)$ of (1.11) in $(0,\tau_0)$ with 
$$\left\{\aligned
&\gamma (0)=p_0\\
&\lim_{\tau\to 0}\tau^p\gamma'(\tau)=v.\endaligned\right.\tag 1.18
$$
for some constant $\tau_0\in (0,t_0^{1-p}]$ where $(0,\tau_0)$ is the maximal
interval of existence of the solution. If $\tau_0<t_0^{1-p}$, 
then
$$
\lim_{\tau\to \tau_0^-}d_0(p_0,\gamma(\tau))=\infty.\tag 1.19
$$
If the Ricci curvature of $M$ is uniformly bounded on $[0,t_0)$, then 
$\tau_0=t_0$.
\endproclaim

We will now prove that the $L_p(q,\2{\tau})$-length can be realized by 
some $\Cal{L}_p(q,\2{\tau})$-geodesic in $M$. We first recall a lemma of 
\cite{Ye1}:

\proclaim{\bf Lemma 1.6}(Lemma 2.1 of \cite{Ye1})
If there exists a constant $c_1>0$ such that 
$$
\text{Ric}(q,\tau)\ge -c_1g(q,\tau)\qquad\text{ on }M\times [0,\2{\tau}],
\tag 1.20
$$
then
$$
e^{-2c_1\tau}g(0)\le g(\tau)\le e^{2c_1(\2{\tau}-\tau)}g(\2{\tau})\quad
\text{ on }M\times [0,\2{\tau}].
$$
If there exists a constant $c_2>0$ such that 
$$
\text{Ric}(q,\tau)\le c_2 g(q,\tau)\quad
\text{ on }M\times [0,\2{\tau}],\tag 1.21
$$
then
$$
e^{2c_2(\tau-\2{\tau})}g(\2{\tau})\le g(\tau)\le e^{2c_2\tau}g(0)\quad
\text{ on }M\times [0,\2{\tau}].
$$
\endproclaim

\proclaim{\bf Lemma 1.7}
Suppose there exists a constant $c_1>0$ such that (1.20) holds. Then
for any $\gamma\in\Cal{F}(q,\2{\tau})$,
$$
\Cal{L}_p(q,\gamma,\2{\tau})\ge -\frac{c_1n}{p+1}\2{\tau}^{p+1}
+\frac{(1-p)e^{-2c_1\2{\tau}}}{\tau_2^{1-p}-\tau_1^{1-p}}d_0(\gamma (\tau_1),
\gamma (\tau_2))^2\quad\forall 0\le\tau_1\le\tau_2\le\2{\tau}.\tag 1.22
$$
\endproclaim
\demo{Proof}
By Lemma 1.6 and the H\"older inequality,
$$\align
\biggl (e^{-c_1\2{\tau}}\int_{\tau_1}^{\tau_2}|\gamma'(\tau)|_{g(0)}\,
d\tau\biggr )^2
\le&\biggl (\int_{\tau_1}^{\tau_2}|\gamma'(\tau)|d\tau\biggr )^2\\
\le&\int_{\tau_1}^{\tau_2}\tau^p|\gamma'(\tau)|^2d\tau\cdot
\int_{\tau_1}^{\tau_2}\tau^{-p}d\tau\\
=&\frac{\tau_2^{1-p}-\tau_1^{1-p}}{1-p}
\int_0^{\2{\tau}}\tau^p|\gamma'(\tau)|^2d\tau
\quad\forall 0\le\tau_1\le\tau_2\le\2{\tau}.\tag 1.23
\endalign
$$
By (1.20),
$$
\int_0^{\2{\tau}}\tau^pR(\gamma (\tau),\tau)d\tau\ge
-\frac{c_1n}{p+1}\2{\tau}^{p+1}.\tag 1.24
$$
By (1.23) and (1.24) the lemma follows.
\enddemo

\proclaim{\bf Lemma 1.8}
Let $r_0>0$. Suppose there exists a constant $c_2>0$ such that (1.21) holds
in $\2{B_0(p_0,r_0)}\times [0,\2{\tau}]$. Then
$$
L_p(q,\tau)\le\frac{c_2n}{p+1}\tau^{p+1}+\frac{e^{2c_2\tau}}{p+1}
\frac{d_0(p_0,q)^2}{\tau^{1-p}}\quad\forall q\in B_0(p_0,r_0),
0<\tau\le\2{\tau}.\tag 1.25
$$
\endproclaim
\demo{Proof}
Let $q\in B_0(p_0,r_0)$, $\tau\in (0,\2{\tau}]$, and let 
$\gamma :[0,\tau]\to M$ be a minimizing geodesic joining $p_0$ and 
$q$ with respect to the metric $g(0)$ with $|\gamma'|_{g(0)}
=d_0(p_0,q)/\tau$ on $[0,\tau]$. Then by Lemma 1.6,
$$\align
L_p(q,\tau)\le&\Cal{L}_p(q,\gamma,\tau)
\le\frac{c_2n}{p+1}\tau^{p+1}
+e^{2c_2\tau}\int_0^{\tau}\rho^p|\gamma'(\rho)|_{g(0)}^2d\rho\\
\le&\frac{c_2n}{p+1}\tau^{p+1}+\frac{e^{2c_2\tau}}{p+1}
\frac{d_0(p_0,q)^2}{\tau^{1-p}}.\endalign
$$
and (1.25) follows.
\enddemo

\proclaim{\bf Lemma 1.9}
Let $\4{\gamma}:[0,\2{s}]\to M$ be a continuous curve satisfying
$$
\int_0^{\2{s}}|\4{\gamma}'(s)|^2\,ds<\infty\tag 1.26
$$
where $|\4{\gamma}'(s)|=|\4{\gamma}'(s)|_{\4{g}(s)}$ and let $\4{Y}(s)\not
\equiv 0$ be a 
smooth vector field along $\4{\gamma}$. Then there exists a variation 
$f:[0,\2{s}]\times [-\3,\3]\to M$ of $\4{\gamma}$ with respect to $\4{Y}(s)$ 
and a constant $C>0$ such that 
$$
\biggl |\frac{\1 f}{\1 s}(s,z)\biggr |^2+\biggl |\nabla_z\biggl (
\frac{\1 f}{\1 s}\biggr )(s,z)
\biggr |^2\le C (|\4{\gamma}'(s)|^2+|\nabla_{\4{X}}\4{Y}(s)|^2)\quad
\forall |z|\le\3\text{ and a.e. }s\in (0,\2{s})\tag 1.27
$$
where $\4{X}(s)=\4{\gamma}'(s)$ and
$$
\lim_{z\to 0}\int_0^{\2{s}}<\frac{\1 f}{\1 s}(s,z),\nabla_z\biggl (
\frac{\1 f}{\1 s}\biggr )
(s,z)>\,ds=\int_0^{\2{s}}<\4{X},\nabla_{\4{X}}\4{Y}>\,ds.\tag 1.28
$$
\endproclaim
\demo{Proof}
For any $s\in [0,\2{s}]$, let $\beta (z,s)
=\beta (z,\4{\gamma}(s),\4{Y}(s))$ be the geodesic with 
respect to the metric $\4{g}(s)$ which satisfies
$$\left\{\aligned
&\beta (0,s)=\4{\gamma}(s)\\
&\frac{\1\beta}{\1 z}(0,s)=\4{Y}(s).\endaligned\right.\tag 1.29
$$
Let $f(s,z)=\beta (z,\4{\gamma}(s),\4{Y}(s))$. By the same argument as the
proof of Proposition 2.2 of Chapter 9 of \cite{C} there exists a constant 
$\3_0>0$ such that for any $0<\3<\3_0$ $f$ is well defined on $[0,\2{s}]\times 
[-\3,\3]$ and is a variation of $\4{\gamma}$ with respect to $\4{Y}(s)$. We 
claim that there exists a constant $0<\3<\3_0$ to be determined later such 
that $f$ satisfies (1.27) and (1.28). Since $f([0,\2{s}]\times [-\3,\3])$ is
compact, there exists a finite family of co-ordinate charts 
$\{(U_i,\phi_i)\}_{i=1}^{i_0}$ such that $f([0,\2{s}]\times [-\3,\3])\subset
\cup_{i=1}^{i_0}U_i$. Without loss of generality we may assume that 
$f([0,\2{s}]\times [-\3,\3])\subset U_1$. We write $\beta=(\beta^1,\beta^2,
\dots,\beta^n)$, $f=(f^1,f^2,\dots,f^n)$, $\4{\gamma}(s)=(a^1(s),a^2(s),\dots,
a^n(s))$, and $\4{Y}(s)=b^i(s)\1/\1 x_i$ in the local coordinates 
$(U_1,\phi_1)$. Then
$$
\frac{\1^2\beta^k}{\1 z^2}+\frac{\1\beta^i}{\1 z}\frac{\1\beta^j}{\1 z}
\4{\Gamma}_{ij}^k(\beta(z,s),s)=0\quad\forall |z|\le\3,0\le s\le\2{s},
k=1,2,\dots,n.\tag 1.30
$$ 
Let
$$
E(z,s)=\4{g}_{kl}(\beta (z,s),s)\frac{\1\beta^k}{\1 z}
\frac{\1\beta^l}{\1 z}.
$$
By (1.29),
$$
E(0,s)=|\4{Y}(s)|^2\le\max_{0\le s\le\2{s}}|\4{Y}(s)|^2=C_1\quad (\text{say})
\tag 1.31
$$
By (1.30),
$$\align
\biggl |\frac{\1 E}{\1 z}\biggr |
=&\biggl |2\4{g}_{kl}\frac{\1\beta^l}{\1 z}\frac{\1^2\beta^k}{\1 z^2}
+\frac{\1\4{g}_{kl}}{\1 x_m}\frac{\1\beta^m}{\1 z}\frac{\1\beta^k}{\1 z}
\frac{\1\beta^l}{\1 z}\biggr |\\
\le&2\biggl |\4{g}_{kl}\frac{\1\beta^l}{\1 z}\frac{\1\beta^i}{\1 z}
\frac{\1\beta^j}{\1 z}\4{\Gamma}_{ij}^k(\beta(z,s),s)\biggr |
+\biggl |\frac{\1\4{g}_{kl}}{\1 x_m}\frac{\1\beta^m}{\1 z}
\frac{\1\beta^k}{\1 z}\frac{\1\beta^l}{\1 z}\biggr |\\
\le&C_2^2E^{\frac{3}{2}}\qquad\qquad\qquad\quad
\forall |z|\le\3,0\le s\le\2{s}\tag 1.32
\endalign
$$
for some constant $C_2>0$. Let $\3=\min (1/(C_1^{\frac{1}{2}}C_2^2),\3_0)$. 
Then by (1.31) and (1.32), 
$$\align
&\biggl |\frac{1}{\sqrt{E(0,s)}}-\frac{1}{\sqrt{E(z,s)}}\biggr |
\le\frac{C_2^2}{2}|z|\le\frac{C_2^2}{2}\3\qquad\qquad\qquad\quad
\forall |z|\le\3,0\le s\le\2{s}\\
\Rightarrow\quad&\frac{1}{\sqrt{E(z,s)}}\ge\frac{1}{\sqrt{E(0,s)}}
-\frac{C_2^2}{2}\3\ge\frac{1}{\sqrt{C_1}}-\frac{C_2^2}{2}\3
\ge\frac{1}{2\sqrt{C_1}}\quad\forall |z|\le\3,0\le s\le\2{s}\\
\Rightarrow\quad&E(z,s)\le 4C_1\qquad\qquad\qquad\qquad\qquad\qquad\qquad
\qquad\qquad\quad\forall |z|\le\3,0\le s\le\2{s}.
\tag 1.33
\endalign
$$
Let 
$$
w(z,s)=\frac{\1\beta}{\1 s}
$$
and let
$$
F(z,s)=|w|^2+|\nabla_zw|^2.
$$
By (1.26), 
$$
\frac{da^i}{ds}\in L^2(0,\2{s})\subset L^1(0,\2{s})\quad\forall i=1,2,\dots,n.
$$
By standard theory on analysis for any $i=1,2,\dots,n$, there exists a set 
$E_i\subset [0,\2{s}]$ of measure zero such that $da^i/ds$ is continuous on 
$(0,\2{s})\setminus
E_i$. Let $E_0=\cup_{i=1}^nE_i$ and $A_0=(0,\2{s})\setminus E_0$. Then 
$|E_0|=0$ and $da^i/ds$ is continuous on $A_0$ for all $i=1,2,\dots,n$. 

We write $\nabla_{\4{X}}\4{Y}(s)=c^i(s)\1/\1 x_i$ in local coordinates. Then
$$
c^k(s)=\frac{db^k}{ds}+b^j\frac{da^i}{ds}\4{\Gamma}_{ij}^k
(\4{\gamma}(s),s)\quad\forall 0\le s\le\2{s},k=1,2,\dots,n.\tag 1.34
$$
By (1.26) and (1.34) $\nabla_{\4{X}}\4{Y}(s)\in L^2(0,\2{s})$ and 
$\nabla_{\4{X}}\4{Y}(s)$ is continuous on $A_0$. Differentiating (1.30) and 
(1.29) with respect to $s\in A_0$,
$$\align
0=&\nabla_s\nabla_z\frac{\1\beta}{\1 z}
=\nabla_z\nabla_s\frac{\1\beta}{\1 z}
+\4{R}(\frac{\1\beta}{\1 s},\frac{\1\beta}{\1 z})\frac{\1\beta}{\1 z}
=\nabla_z\nabla_zw
+\4{R}(w,\frac{\1\beta}{\1 z})\frac{\1\beta}{\1 z}
\tag 1.35
\endalign
$$ 
holds for any $|z|\le\3,s\in A_0,k=1,2,\dots,n$ and
$$\left\{\aligned
&w(0,s)=\4{\gamma}'(s)\qquad\quad\,\,\forall s\in A_0\\
&\nabla_zw(0,s)=\nabla_{\4{X}}\4{Y}(s)\quad\forall s\in A_0.
\endaligned\right.\tag 1.36
$$
Then by (1.33), (1.35), (1.36), and H\"older's inequality, 
$$\align
\biggl |\frac{\1F}{\1 z}\biggr |
=&2\biggl |<w,\nabla_zw>+<\nabla_zw,\nabla_z\nabla_zw>\biggr |\\
\le&F+2\biggl |<\nabla_zw,\4{R}(w,\frac{\1\beta}{\1 z})
\frac{\1\beta}{\1 z}>\biggr |\\
\le&F+|\nabla_zw|^2+\biggl |\4{R}(w,\frac{\1\beta}{\1 z})
\frac{\1\beta}{\1 z}\biggr |^2\\
\le&C_3F\quad\forall |z|\le\3, s\in A_0\\
\Rightarrow\quad F(z,s)\le&e^{C_3z}F(0,s)
\le e^{C_3\3}(|\4{\gamma}'(s)|^2+|\nabla_{\4{X}}\4{Y}|^2)
\endalign
$$
holds for any $|z|\le\3$ and $s\in A_0$ where $C_3>0$ is a constant.
Hence (1.27) follows. By (1.27) and the Cauchy-Schwarz inequality,
$$
\biggl |<\frac{\1 f}{\1 s}(s,z),\nabla_z\biggl (\frac{\1f}{\1 s}\biggr )
(s,z)>\biggr |\le C(|\4{\gamma}'(s)|^2+|\nabla_{\4{X}}\4{Y}|^2)
\quad\forall |z|\le\3, s\in A_0.\tag 1.37
$$
Since both $\gamma'(s)$ and $\nabla_{\4{X}}\4{Y}$ are continuous at $s$ for
any $s\in A_0$, by (1.35), (1.36), and the continuous dependence of solutions
of O.D.E. on initial data, $\nabla_zw$ is continuous at $(z,s)$ for any 
$s\in A_0$ and $|z|\le\3$. Hence
$$
\nabla_z\biggl (\frac{\1 f}{\1 s}\biggr )(s,z)=\nabla_s\biggl (
\frac{\1 f}{\1 z}\biggr )(s,z)\quad\forall |z|\le\3, s\in A_0.\tag 1.38
$$
Since by (1.38) 
$$
<\frac{\1 f}{\1 s}(s,z),\nabla_z\biggl (\frac{\1 f}{\1 s}\biggr )(s,z)>\to
<\4{X},\nabla_{\4{X}}\4{Y}>\quad\text{ as }z\to 0\quad\forall s\in A_0,
$$
by (1.37) and the Lebesgue dominated convergence theorem (1.28) follows.
\enddemo

\proclaim{\bf Theorem 1.10}
Let $\2{\tau}\in (0,t_0)$ and $\2{s}=\2{\tau}^{1-p}$. Suppose $(M,g)$ 
satisfies (1.20) for some constant $c_1>0$. Then for 
any $q\in M$, there exists a $\4{\Cal{L}}_p$-geodesic
$\4{\gamma}\in C^1([0,\2{s}])\cap C^{\infty}((0,\2{s}])$ such that 
$\4{\gamma}(0)=p_0$, $\4{\gamma}(\2{s})=q$, and
$$
\4{L}_p(q,\2{s})=\4{\Cal{L}}_p(q,\4{\gamma},\2{s}).\tag 1.39
$$
\endproclaim
\demo{Proof}
Choose a sequence of curves $\{\4{\gamma}_i\}_{i=1}^{\infty}\subset
\Cal{F}(q,\2{s})$ such that 
$$
\4{\Cal{L}}_p(q,\4{\gamma}_i,\2{s})\le \4{L}_p(q,\2{s})+1\quad\forall i\in
\Bbb{Z}^+\tag 1.40
$$
and
$$
\4{L}_p(q,\2{s})=\lim_{i\to\infty}\4{\Cal{L}}_p(q,\4{\gamma}_i,\2{s}).
\tag 1.41
$$
Let $\gamma_i(\tau)=\4{\gamma}_i(s)$, $\tau=s^{1/(1-p)}$. By Lemma 1.7, (1.4),
(1.6) and (1.41) there exist constants $K=K(\2{\tau},L_p(q,\2{\tau}))>0$
and $C_1>0$ independent of $i\in\Bbb{Z}^+$  
such that
$$\left\{\aligned
&d_0(p_0,\4{\gamma}_i(s))\le K\quad\forall 0\le s\le\2{s},i\in\Bbb{Z}^+\\
&d_0(\4{\gamma}_i(s),\4{\gamma}_i(s'))\le C_1|s-s'|^{1/2}\quad\forall 
s,s'\in [0,\2{s}],i\in\Bbb{Z}^+.\endaligned\right.
\tag 1.42
$$
Hence the sequence of curves $\{\4{\gamma}_i\}_{i=1}^{\infty}$ are uniformly
H\"older continuous on $[0,\2{s}]$. Since $M$ is complete with respect to 
$g(0)$, $\2{B_0(p_0,K)}$ is compact. By the Ascoli Theorem there exists a 
continuous curve $\4{\gamma}:[0,\2{s}]\to\2{B_0(p_0,K)}$ such that 
$\4{\gamma}_i$ converges uniformly to $\4{\gamma}$ on $[0,\2{s}]$ as 
$i\to\infty$. Then $\4{\gamma}(0)=p_0$ and $\4{\gamma}(\2{s})=q$. 
Letting $i\to\infty$ in (1.42),
$$\left\{\aligned
&d_0(p_0,\4{\gamma}(s))\le K\qquad\qquad\qquad\,\,\,\,\forall 0\le s\le\2{s}\\
&d_0(\4{\gamma}(s),\4{\gamma}(s'))\le C_1|s-s'|^{1/2}\quad\forall 
s,s'\in [0,\2{s}].\endaligned\right.
$$
Hence $\4{\gamma}$ is uniformly H\"older continuous on $[0,\2{s}]$. By 
Fatou's Lemma and Lemma 1.6,
$$\align
&\4{\Cal{L}}_p(q,\4{\gamma},\2{s})\le\lim_{i\to\infty}\4{\Cal{L}}_p
(q,\4{\gamma}_i,\2{s})=\4{L}_p(q,\2{s})\\
\Rightarrow\quad&\4{L}_p(q,\2{s})=\4{\Cal{L}}_p(q,\4{\gamma},\2{s})
\quad\text{ and }\quad|\4{\gamma}'|_{g(0)}\in L^2(0,\2{s}).\tag 1.43
\endalign
$$
We now claim that $\4{\gamma}\in C^{\infty}([0,\2{s}])$. Since 
$\2{B_0(p_0,K)}$ is compact, there exists a finite family of co-ordinate 
charts $\{(\phi_k,B_0(q_k,r_k)\}_{k=1}^{k_0}$ such that $\2{B_0(p_0,K)}\subset
\cup_{k=1}^{k_0}B_0(q_k,r_k)$. Let
$$
I_k=I_k(\4{\gamma})=\{s\in [0,\2{s}]:\4{\gamma}(s)\in B_0(q_k,r_k)\}
\quad\forall k=1,2,\dots,k_0.
$$ 
Then $I_k$ is relatively open with respect to the interval $[0,\2{s}]$
for all $k=1,2,\dots,k_0$ and $[0,\2{s}]=\cup_{k=1}^{k_0}I_k$. For any 
$k=1,2,\dots,k_0$, $s\in I_k$, we write 
$$
\phi_k(\4{\gamma}(s))=(a_k^1(s),a_k^2(s),\dots,a_k^n(s))
$$ 
in the local coordinates $(\phi_k,B_0(q_k,r_k))$. When there is no ambiguity, 
we will drop the subscript $k$. To prove the claim we fix one $k\in\{1,2,
\dots,k_0\}$. Then
$$
\4{X}(s)=\4{\gamma}'(s)=\frac{d a^i}{d s}\frac{\1}{\1 x_i}\quad\text{ in }I_k.
$$
By (1.43),
$$
\4{g}_{ij}(\4{\gamma}(s),0)\frac{d a^i}{d s}\frac{d a^j}{d s}\in L^1(I_k)
\quad\Rightarrow\quad\frac{d a^i}{d s}\in L^2(I_k)\quad
\forall i=1,2,\dots,n.\tag 1.44
$$
Let 
$$
\4{Y}(s)=b^j(s)\frac{\1}{\1 x_j}
$$
be a smooth vector field along $\4{\gamma}$ such that $\4{Y}(s)=0$ for any 
$s\not\in I_k$. Since $\4{\gamma}$ is a minimizer of $\4{L}_p(q,\2{s})$ and
by (1.43) $\4{\gamma}$ satisfies (1.26), by Lemma 1.9, Lebesgue dominated
convergence theorem and an argument similar to the proof of Lemma 1.2,
$$\align
&\int_0^{\2{s}}(s^{\frac{2p}{1-p}}\4{Y}(\4{R})
+2(1-p)^2<\4{X},\nabla_{\4{X}}\4{Y}>)\,ds=0\\
\Rightarrow\quad&\int_{I_k}\biggl\{s^{\frac{2p}{1-p}}b^j\frac{\1\4{R}}{\1 x^j}
+2(1-p)^2\4{g}_{lr}\frac{da^l}{ds}\biggl (
\frac{db^r}{ds}+b^j\frac{da^i}{ds}\4{\Gamma}_{ij}^r\biggr )
\biggr\}\,ds=0\tag 1.45\\
\Rightarrow\quad&\biggl |\int_{I_k}\4{g}_{lr}\frac{da^l}{ds}\biggl (
\frac{db^r}{ds}+b^j\frac{da^i}{ds}\4{\Gamma}_{ij}^r\biggr )\,ds\biggl |
\le C\biggl |\int_{I_k}s^{\frac{2p}{1-p}}b^j\frac{\1\4{R}}{\1 x^j}\,ds
\biggr |\\
\Rightarrow\quad&\biggl |\int_{I_k}\4{g}_{lr}\frac{da^l}{ds}
\frac{db^r}{ds}\,ds\biggr |
\le\biggl |\int_{I_k}\4{g}_{lr}\frac{da^l}{ds}b^j\frac{da^i}{ds}
\4{\Gamma}_{ij}^r\,ds\biggr |
+C\sum_{j=1}^n\int_{I_k}|b^j|\,ds.\tag 1.46
\endalign
$$
Since by (1.44) and the Cauchy-Schwarz inequality,
$$
\biggl |\int_{I_k}\4{g}_{lr}\frac{da^l}{ds}b^j\frac{da^i}{ds}
\4{\Gamma}_{ij}^r\,ds\biggr |\le C\max_{1\le j\le n}\|b^j\|_{\infty}
\int_J\4{g}_{lr}\frac{da^l}{ds}\frac{da^r}{ds}\,ds
$$
where $J=\cup_{j=1}^n\text{supp }b^j$, by (1.46) we have
$$
\biggl |\int_{I_k}\4{g}_{lr}\frac{da^l}{ds}\frac{db^r}{ds}\,ds\biggr |
\le C\max_{1\le j\le n}\|b^j\|_{\infty}\int_{J}\4{g}_{lr}\frac{da^l}{ds}
\frac{da^r}{ds}\,ds
+C\sum_{j=1}^n\int_{I_k}|b^j|\,ds.\tag 1.47
$$
We now choose $\phi\in C^{\infty}(\Bbb{R})$, $0\le\phi\le 1$, such that 
$\phi (s)=1$ for all $s\le 0$ and $\phi (s)=0$ for all $s\ge 1$. For any 
$h>0$ let $\phi_{h}(s)=\phi (s/h)$. Since
$$
\int_{\Bbb{R}}\phi '(s)\,ds=-1,
$$
by (1.44) and standard theory in analysis \cite{St} there exists a set
$E_k\subset I_k$ of measure zero such that
$$\left\{\aligned
\lim_{h\to 0^+}\frac{1}{h}\int_{s'-h}^{s'}\frac{da^1}{ds}(s)\phi '
((s'-s)/h)\,ds=-\frac{da^1}{ds}(s')\quad\forall s'\in I_k\setminus E_k\\
\lim_{h\to 0^+}\frac{1}{h}\int_{s'}^{s'+h}\frac{da^1}{ds}(s)\phi '
((s-s')/h)\,ds=-\frac{da^1}{ds}(s')\quad\forall s'\in I_k\setminus E_k.
\endaligned\right.\tag 1.48
$$
Let $s_0'\in I_k\setminus E_k$. Without loss of generality we may assume
that $s_0'\ne 0,\2{s}$. By continuity of $\4{\gamma}$ there exists a constant 
$0<\3\le s_0'/2$ such that
$$
\4{\gamma}(s)\in B_0(q_k,r_k)\quad\forall s\in I_{\3}(s_0')
=(s_0'-\3,s_0'+\3)\subset I_k.
$$ 
Let $s_1,s_2\in I_{\3/2}(s_0')\setminus E_k$ be such that $s_1<s_2$ and let
$\3_1=(s_2-s_1)/3$. Putting 
$$
b^r(s)=\4{g}^{r1}(\4{\gamma}(s),s)\phi_{h_2}(s-s_2)\phi_{h_1}(s_1-s)
\quad\forall r=1,2,\dots,n
$$
in (1.47) where $0<h_1,h_2\le\3_1$ we get
$$\align
&\biggl |\frac{1}{h_2}\int_{s_2}^{s_2+h_2}\frac{da^1}{ds}(s)
\phi '((s-s_2)/h_2)\,ds
-\frac{1}{h_1}\int_{s_1-h_1}^{s_1}\frac{da^1}{ds}(s)\phi '((s_1-s)/h_1)\,ds
\biggr |\\
\le &C\int_{s_1-h_1}^{s_2+h_2}\4{g}_{lr}\frac{da^l}{ds}
\frac{da^r}{ds}\,ds+C(s_2-s_1+h_1+h_2)\\
&\qquad +\biggl |\int_{s_1-h_1}^{s_2+h_2}\4{g}_{lr}\frac{d\4{g}^{r1}}{ds}
\frac{da^l}{ds}\phi_{h_2}(s-s_2)\phi_{h_1}(s_1-s)\,ds\biggr |\\
\le &C\int_{s_1-h_1}^{s_2+h_2}\4{g}_{lr}\frac{da^l}{ds}
\frac{da^r}{ds}\,ds+C(s_2-s_1+h_1+h_2)\\
&\qquad +\biggl |\int_{s_1-h_1}^{s_2+h_2}\4{g}_{lr}\biggl (
\frac{\1\4{g}^{r1}}{\1 s}+\frac{\1\4{g}^{r1}}{\1 x_j}\frac{d a^j}{d s}\biggr )
\frac{da^l}{ds}\phi_{h_2}(s-s_2)\phi_{h_1}(s_1-s)\,ds\biggr |\\
\le &C\int_{s_1-h_1}^{s_2+h_2}\4{g}_{lr}\frac{da^l}{ds}
\frac{da^r}{ds}\,ds+C(s_2-s_1+h_1+h_2)\\
&\qquad +\frac{2}{1-p}\biggl |\int_{s_1-h_1}^{s_2+h_2}s^{\frac{p}{1-p}}
\4{g}_{rl}\4{g}^{ri}\4{R}_{ij}\4{g}^{j1}\frac{da^l}{ds}\phi_{h_2}(s-s_2)
\phi_{h_1}(s_1-s)\,ds\biggr |\quad (\text{by (1.8)})\\
\le &C\int_{s_1-h_1}^{s_2+h_2}\4{g}_{lr}\frac{da^l}{ds}
\frac{da^r}{ds}\,ds+C(s_2-s_1+h_1+h_2)\\
&\qquad +C(s_2-s_1+h_1+h_2)^{1/2}\biggl (\int_{s_1-h_1}^{s_2+h_2}\4{g}_{lr}
\frac{da^1}{ds}\frac{da^r}{ds}\,ds\biggr )^{1/2}\\
\le &C\int_{s_1-h_1}^{s_2+h_2}\4{g}_{lr}\frac{da^l}{ds}
\frac{da^r}{ds}\,ds+C(s_2-s_1+h_1+h_2).\tag 1.49
\endalign
$$
Letting $h_1\to 0$ in (1.49), by (1.26), (1.44), and (1.48),
$$
\biggl |\frac{da^1}{ds}(s_1)\biggr |
\le C\quad\forall s_1\in I_{\3/2}(s_0')\setminus E_k.\tag 1.50
$$ 
for some constant $C>0$. Hence letting $h_1, h_2\to 0$ in (1.49), by (1.26),
(1.48), and (1.50) we have
$$
\biggl |\frac{da^1}{ds}(s_2)-\frac{da^1}{ds}(s_1)\biggr |
\le C|s_2-s_1|\quad\forall s_1,s_2\in I_{\3/2}(s_0')\setminus E_k.\tag 1.51
$$
We now choose $0\le \eta\in C_0^{\infty}(\Bbb{R})$ such that $\eta (s)=0$ for 
any $|s|\ge 1$ and $\int_{\Bbb{R}}\eta\,ds=1$. For any $h>0$, $s_0\in I_k$,
let $\eta_h(s)=\eta(s/h)/h$ and 
$$
a^1\ast\eta_h(s_0)=\int_{\Bbb{R}}a^1(s_0-s)\eta_h(s)\,ds. 
$$
For any $s_0\in I_{\3/2}(s_0')$, we choose a sequence 
$\{s_{0,i}\}_{i=1}^{\infty}\subset I_{\3/2}(s_0')\setminus E_k$ such that 
$\lim_{i\to\infty}s_{0,i}=s_0$. By (1.51) 
$\{da^1(s_{0,i})/ds\}_{i=1}^{\infty}$ is a Cauchy sequence. Hence
$$
\lim_{i\to\infty}\frac{da^1}{ds}(s_{0,i})
$$
exists. Let
$$
f(s_0)=\lim_{i\to\infty}\frac{da^1}{ds}(s_{0,i})\quad\forall s_0\in 
I_{\3}(s_0').
$$
By (1.51) $f$ is well defined on $I_{\3/2}(s_0')$. We now claim that $a^1\in 
C^1(I_{\3/2}(s_0'))$ with $da^1/ds=f$ on $I_{\3/2}(s_0')$. To prove the 
claim we observe that by (1.51)
$$\align
&\biggl |\frac{d}{ds}a^1\ast\eta_h(s_0)-f(s_0)\biggr |\\
=&\biggl |\int_{|\rho|\le h}\frac{da^1}{ds}(s_0-\rho)\eta_h(\rho)\,d\rho
-f(s_0)\biggr |\\
\le&\biggl |\int_{|\rho|\le h}\biggl (\frac{da^1}{ds}(s_0-\rho)-\frac{da^1}{ds}
(s_{0,i})\biggr )\eta_h(\rho)\,d\rho\biggr |+\biggl |\frac{da^1}{ds}(s_{0,i})
-f(s_0)\biggr |\\
\le&C\int_{|\rho|\le h}|s_0-\rho -s_{0,i}|\eta_h(\rho)\,d\rho
+\biggl |\frac{da^1}{ds}(s_{0,i})-f(s_0)\biggr |\\
\le&C(|s_0-s_{0,i}|+h)
+\biggl |\frac{da^1}{ds}(s_{0,i})-f(s_0)\biggr |.\endalign
$$
Letting first $i\to\infty$ and then $h\to 0$ in the above inequality, we get 
that $d(a^1\ast\eta_h)/ds$ converges uniformly to $f$ on $I_{\3/2}(s_0')$ as 
$h\to 0$.
Since $a^1\ast\eta_h$ converges uniformly to $a^1$ on $I_{\3/2}(s_0')$ as 
$h\to 0$. Hence $a^1\in C^1(I_{\3}(s_0'))$ with $da^1/ds=f$ on $I_{\3/2}
(s_0')$.
Since $s_0'$ is arbitrary, $a^1\in C^1([0,\2{s}])$. By a similar 
argument $a^l\in C^1([0,\2{s}])$ for any $l=1,2,\dots,n$. Hence 
$\4{\gamma}\in C^1([0,\2{s}])$ and there exists a constant $C>0$ such that
$$
\biggl |\frac{da_k^l}{ds}(s_2)-\frac{da_k^l}{ds}(s_1)\biggr |
\le C|s_2-s_1|\quad\forall s_1,s_2\in I_k,k=1,2,\dots,k_0,l
=1,2,\dots,n.\tag 1.52
$$
Thus $d^2a_k^l(s)/ds^2$ exists a.e. $s\in I_k$ with
$$
\frac{d^2a_k^l}{ds^2}\in L^{\infty}(I_k)\quad\forall k=1,2,\dots,
k_0,l=1,2,\dots,n.
$$
By (1.45) for any $b^j\in C_0^{\infty}(I_k)$, $j=1,2,\dots, n$,
$$\align
&\int_{I_k}\biggl\{s^{2p/(1-p)}\frac{\1\4{R}}{\1 x_j}
+2(1-p)^2\biggl (-\frac{d }{ds}\biggl (\4{g}_{jl}
\frac{da^l}{ds}\biggr )
+\4{g}_{lr}\frac{da^l}{ds}\frac{da^i}{ds}\4{\Gamma}_{ij}^r\biggr )
\biggr\}b^j\,ds=0\\
\Rightarrow\quad&2(1-p)^2\biggl (-\frac{d }{ds}\biggl (\4{g}_{jl}
\frac{da^l}{ds}\biggr )
+\4{g}_{lr}\frac{da^l}{ds}\frac{da^i}{ds}\4{\Gamma}_{ij}^r\biggr )
+s^{2p/(1-p)}\frac{\1\4{R}}{\1 x_j}=0\quad\text{ a.e. }s\in I_k
\tag 1.53\\
\Rightarrow\quad&2(1-p)^2\biggl (-\4{g}_{jl}\frac{d^2a^l}{ds^2}
-\frac{d\4{g}_{jl}}{ds}\frac{da^l}{ds}
+\4{g}_{lr}\frac{da^l}{ds}\frac{da^i}{ds}\4{\Gamma}_{ij}^r\biggr )
+s^{2p/(1-p)}\frac{\1\4{R}}{\1 x_j}=0\quad\text{ a.e. }s\in I_k\\
\Rightarrow\quad&\frac{d^2a^m}{ds^2}
=-\4{g}^{mj}\frac{d\4{g}_{jl}}{ds}\frac{da^l}{ds}
+\4{g}_{lr}\4{g}^{mj}\frac{da^l}{ds}\frac{da^i}{ds}\4{\Gamma}_{ij}^r
+\frac{1}{(1-p)^2}s^{2p/(1-p)}\4{g}^{mj}\frac{\1\4{R}}{\1 x_j}
\tag 1.54
\endalign
$$
a.e. $s\in I_k$. By (1.52) and (1.54) there exists a set $E_k'$ of
measure zero such that 
$$
\biggl |\frac{d^2a_k^j}{ds^2}(s_2)-\frac{d^2a_k^j}{ds^2}(s_1)\biggr |
\le C|s_2-s_1|\quad\forall s_1,s_2\in I_{\3}(s_0')\setminus E_k',
k=1,2,\dots,k_0,j=1,2,\dots,n.\tag 1.55
$$
By an argument similar to the proof of $\4{\gamma}\in C^1([0,\2{s}])$
but with (1.55) replacing (1.51) in the proof we get that 
$\4{\gamma}\in C^2([0,\2{s}])$ and
$$
\biggl |\frac{d^2a_k^j}{ds^2}(s_2)-\frac{d^2a_k^j}{ds^2}(s_1)\biggr |
\le C|s_2-s_1|\quad\forall s_1,s_2\in I_k,
k=1,2,\dots,k_0,j=1,2,\dots,n
$$
for some constant $C>0$. Hence by (1.53),
$$
2(1-p)^2\biggl (-\frac{d }{ds}\biggl (\4{g}_{jl}
\frac{da^l}{ds}\biggr )
+\4{g}_{lr}\frac{da^l}{ds}\frac{da^i}{ds}\4{\Gamma}_{ij}^r\biggr )
+s^{2p/(1-p)}\frac{\1\4{R}}{\1 x_j}=0\quad\text{ in }I_k
$$
for any $k=1,2,\dots,k_0,j=1,2,\dots,n$. Thus $\4{\gamma}$ satisfies
(1.10) in $(0,\2{s})$. Then by standard O.D.E. theory $\4{\gamma}\in
C^{\infty}((0,\2{s}])$. Hence $\4{\gamma}$ is a $\4{\Cal{L}}_p$-geodesic and 
the theorem follows.
\enddemo

By (1.4), (1.5), (1.6) and Theorem 1.10 we have

\proclaim{\bf Theorem 1.11}
Suppose $(M,g)$ satisfies (1.20) in $[0,\2{\tau}]$ for some constant $c_1>0$. 
Then for any $q\in M$, there exists a $\Cal{L}_p$-geodesic
$\gamma\in C([0,\2{\tau}])\cap C^{\infty}((0,\2{\tau}])$ satisfying (1.18) 
for some $v\in T_{p_0}M$ such that $\gamma (\2{\tau})=q$ and 
$$
L_p(q,\2{\tau})=\Cal{L}_p(q,\gamma,\2{\tau}).
$$
\endproclaim

By putting $p=1/2$ in Theorem 1.11 we obtain a result which is 
conjectured and used without proof in Perelman's paper 
\cite{P1}, \cite{P2}.

\proclaim{\bf Corollary 1.12}
Suppose $(M,g)$ satisfies (1.20) in $[0,\2{\tau}]$ for some constant $c_1>0$. 
Then for any $q\in M$, there exists a $\Cal{L}$-geodesic
$\gamma\in C([0,\2{\tau}])\cap C^{\infty}((0,\2{\tau}])$ satisfying (1.18) 
with $p=1/2$ for some $v\in T_{p_0}M$ such that $\gamma (\2{\tau})=q$ and 
$$
L(q,\2{\tau})=\Cal{L}(q,\gamma,\2{\tau}).
$$
\endproclaim

By an argument similar to the proof of Theorem 1.10 and Theorem 1.11 we have

\proclaim{\bf Theorem 1.13}
Suppose $(M,g)$ satisfies (1.20) in $[0,\2{\tau}]$ for some constant $c_1>0$.
Let $q_0\in M$ and let $\tau_0\in (0,\2{\tau})$. 
Then for any $q\in M$, there exists a $\Cal{L}_p$-geodesic $\gamma\in
C^{\infty}([\tau_0,\2{\tau}])$ satisfying $\gamma (\tau_0)=q_0$, 
$\gamma (\2{\tau})=q$, and 
$$
L_p^{q_0,\tau_0}(q,\2{\tau})=\Cal{L}_p^{q_0,\tau_0}(q,\gamma,\2{\tau}).
$$
\endproclaim

\proclaim{\bf Theorem 1.14}
Let $t_0>0$, $s_0=t_0^{1-p}$, and let $g$ and $\2{g}$ be related by (0.5). 
Suppose $(M,\2{g})$ satisfies (1.20) in $(0,t_0)$ for some constant $c_1>0$.
Then for any $\2{s}\in (0,s_0)$, $q\in M$, $\4{L}_p^{p_0}(q,\2{s})$ is 
locally Lipschitz in $p_0$ with respect to the metric $g(0)=\2{g}(t_0)$.
\endproclaim
\demo{Proof}
Let $r_0>0$, $\2{s}\in (0,s_0)$, $\2{\tau}=\2{s}^{\frac{1}{1-p}}$,
$\2{p}_0\in M$, and let $p_1,p_2\in B_0(\2{p}_0,r_0)$. By Theorem
1.10 for each $i=1,2$, there exists a $\4{\Cal{L}}_p^{p_i}(q,\2{s})$-length 
minimizing $\4{\Cal{L}}_p$-geodesic $\4{\gamma}_i:[0,\2{s}]\to M$ such that
$\4{\gamma}_i(0)=p_i$ and $\4{\gamma}_i(\2{s})=q$. Let $\4{\gamma}:[0,d]\to M$ 
be a normalized minimizing geodesic with respect to the metric $g(0)$ with
$\4{\gamma}(0)=p_1$, $\4{\gamma}(d)=p_2$, $|\4{\gamma}'|=|\4{\gamma}'|_{g(0)}
=1$ on $[0,d]$ with $d=d_0(p_1,p_2)$. Then $\4{\gamma}([0,d])\subset 
B_0(\2{p}_0,3r_0)$. Let $r_1=3r_0+2d_0(\2{p}_0,q)$ and let
$$
K_0=\sup_{\2{B_0(\2{p}_0,r_1)}\times [0,\2{\tau}]}(|R|+|\text{Ric}|).
$$
Let $\gamma_i(\tau)=\4{\gamma}_i(s)$ with $s=\tau^{1-p}$, $i=1,2$.
For $i=1,2$, let $d_i=d_0(p_i,q)$ and let $\2{\gamma}_i:[0,d_i]\to M$
be a normalized minimizing geodesic with respect to the metric $g(0)$ with
$\2{\gamma}_i(0)=p_i$, $\2{\gamma}_i(d_i)=q$. Then $d_i
<r_0+d_0(\2{p}_0,q)$ and $\2{\gamma}_i([0,d_i])\subset B_0(p_0,r_1)$ for
$i=1,2$. Hence by Lemma 1.7 and the proof of Lemma 1.8, there exist constants 
$A_1=A_1(\2{s},r_1,K_0)>0$ and $r_2=r_2(\2{s},r_1,K_0)\ge r_1$ such that
$$\left\{\aligned
&\4{L}_p^{p_i}(q,\2{s})\le A_1\qquad\forall i=1,2\\
&d_0(p_i,\4{\gamma}_i(s))<r_2\quad\forall 0\le s\le\2{s},i=1,2.
\endaligned\right.
$$
Let
$$
K_1=\sup_{\2{B_0(p_0,2r_2)}\times [0,\2{\tau}]}(|R|+|R_t|+|\nabla R|+|\text{Ric}|).
$$
We now assume that $d=d_0(p_1,p_2)<\min (1,\2{s}/4)$. Let 
$$
\beta (s)=\left\{\aligned
&\4{\gamma}(s)\qquad\quad\,\,\,\text{ if }0\le s<d\\
&\4{\gamma}_2(s-d)\quad\,\,\text{ if }d\le s\le\2{s}-d\\
&\4{\gamma}_2(2s-\2{s})\quad\text{ if }\2{s}-d<s\le\2{s}.
\endaligned\right.
$$
Then
$$\align
\4{L}_p^{p_1}(q,\2{s})\le&\4{\Cal{L}}_p^{p_1}(q,\beta ,\2{s})\\
=&\frac{1}{1-p}\int_0^d(s^{\frac{2p}{1-p}}\4{R}(\4{\gamma}(s),s)
+(1-p)^2|\4{\gamma}'(s)|^2)\,ds\\
&\qquad +\frac{1}{1-p}\int_d^{\2{s}-d}(s^{\frac{2p}{1-p}}\4{R}
(\4{\gamma}_2(s'),s)+(1-p)^2|\4{\gamma}_2'(s')|_{\4{g}(s)}^2)\,ds\\
&\qquad +\frac{1}{1-p}\int_{\2{s}-d}^{\2{s}}(s^{\frac{2p}{1-p}}\4{R}(\4{\gamma}_2(s''),s)+4(1-p)^2|\4{\gamma}_2'(s'')|_{\4{g}(s)}^2)\,ds\\
=&I_1+I_2+I_3\tag 1.56
\endalign
$$
where $s'=s-d$ and $s''=2s-\2{s}$.
By Lemma 1.6,
$$
I_1\le\frac{K_1}{1+p}d^{\frac{1+p}{1-p}}+(1-p)e^{2K_1d^{\frac{1}{1-p}}}d
\le\frac{K_1}{1+p}d^{\frac{1+p}{1-p}}+(1-p)e^{2K_1}d,\tag 1.57
$$
$$\align
I_2\le&\frac{1}{1-p}\int_0^{\2{s}-2d}((w+d)^{\frac{2p}{1-p}}
\4{R}(\4{\gamma}_2(w),w+d)+(1-p)^2e^{2K_1d^{\frac{1}{1-p}}}
|\4{\gamma}_2'(w)|^2)\,dw\\
=&\frac{1}{1-p}\int_0^{\2{s}-2d}(w^{\frac{2p}{1-p}}\4{R}(\4{\gamma}_2(w),w)
+(1-p)^2|\4{\gamma}_2'(w)|^2)\,dw\\
&\qquad +\frac{1}{1-p}\int_0^{\2{s}-2d}[(w+d)^{\frac{2p}{1-p}}
-w^{\frac{2p}{1-p}}]\4{R}(\4{\gamma}_2(w),w+d)\,dw\\
&\qquad +\frac{1}{1-p}\int_0^{\2{s}-2d}w^{\frac{2p}{1-p}}(\4{R}
(\4{\gamma}_2(w),w+d)-\4{R}(\4{\gamma}_2(w),w))\,dw\\
&\qquad +(1-p)(e^{2K_1d^{\frac{1}{1-p}}}-1)\int_0^{\2{s}-2d}|
\4{\gamma}_2'(w)|^2\,dw\\
\le&\frac{1}{1-p}\int_0^{\2{s}-2d}(w^{\frac{2p}{1-p}}\4{R}(\4{\gamma}_2(w),w)
+(1-p)^2|\4{\gamma}_2'(w)|^2)\,dw+C_1'd\tag 1.58
\endalign
$$
for some constant $C_1'>0$ and
$$\align
I_3\le&\frac{1}{1-p}K_1\2{s}^{\frac{2p}{1-p}}d+2(1-p)e^{2K_1d^{\frac{1}{1-p}}}
\int_{\2{s}-2d}^{\2{s}}|\4{\gamma}_2'(s)|^2\,ds\\
\le&\frac{1}{1-p}\int_{\2{s}-2d}^{\2{s}}(s^{\frac{2p}{1-p}}\4{R}
(\4{\gamma}_2(s),s)
+(1-p)^2|\4{\gamma}_2'(s)|^2)\,ds +C_1''\int_{\2{s}-2d}^{\2{s}}
|\4{\gamma}_2'(s)|^2\,ds\\
&\qquad +\frac{3K_1}{1-p}\2{s}^{\frac{2p}{1-p}}d.
\tag 1.59
\endalign
$$
Hence by (1.56), (1.57), (1.58) and (1.59),
$$
\4{L}_p^{p_1}(q,\2{s})\le\4{L}_p^{p_2}(q,\2{s})
+C_1''\int_{\2{s}-2d}^{\2{s}}|\4{\gamma}_2'(s)|^2\,ds+C_2'd\tag 1.60
$$
for some constant $C_2'>0$. Let $\4{v}_i=\4{\gamma}_i'(0)$ for $i=1,2$.
By the same argument as the proof of Lemma 1.4, (1.17) holds in 
$(0,\2{s})$ for some constant $C_2>0$, $C_3>0$, 
depending only on $K_1$ with $\4{X}(s)$ and $\4{v}$, being replaced by 
$\4{\gamma}_i'(s)$ and $\4{v}_i=\4{\gamma}_i'(0)$, for $i=1,2$. By (1.17) 
there exist constants $C_4>0$, $C_5>0$ and $C_6>0$ such that
$$
C_4|\4{v}_i|^2-C_5\le |\4{\gamma}_i'(s)|^2\le C_6(1+|\4{v}_i|^2)
\quad\forall 0\le s\le\2{s},i=1,2.\tag 1.61
$$
By (1.61),
$$\align 
C_4\2{s}|\4{v}_2|^2=&C_4\int_0^{\2{s}}|\4{v}_2|^2\,d s\le C_5\2{s}
+\int_0^{\2{s}}|\4{\gamma}_2'(s)|^2\,ds\\
\le&C_5\2{s}+\frac{1}{1-p}\biggl [\4{L}_p^{p_2}(q,\2{s})-\frac{1}{1-p}
\int_0^{\2{s}}s^{\frac{2p}{1-p}}\4{R}(\4{\gamma}_2(s),s)\,ds\biggr ]\\
\le&C_5\2{s}+\frac{1}{1-p}\biggl (\4{L}_p^{p_2}(q,\2{s})+\frac{K_1}{1+p}
\2{s}^{\frac{1+p}{1-p}}\biggr )\\
\le&C_5\2{s}++\frac{1}{1-p}\biggl (A_1+\frac{K_1}{1+p}
\2{s}^{\frac{1+p}{1-p}}\biggr )=C_7\quad (\text{say}).
\tag 1.62\endalign
$$
By (1.61) and (1.62),
$$\align
\int_{\2{s}-2d}^{\2{s}}|\gamma_2'(s)|^2\,ds
\le&2C_6(1+|\4{v}_2|^2)d\le 2C_6(1+(C_7/(C_4\2{s}))d.\tag 1.63
\endalign
$$
By (1.60) and (1.63) there exists a constant 
$C_8=C_8>0$ such that
$$
\4{L}_p^{p_1}(q,\2{s})\le\4{L}_p^{p_2}(q,\2{s})+C_8d.
$$
Interchanging the role of $p_1$ and $p_2$ in the above inequality,
$$
\4{L}_p^{p_2}(q,\2{s})\le\4{L}_p^{p_1}(q,\2{s})+C_8d.
$$
Hence
$$
|\4{L}_p^{p_1}(q,\2{s})-\4{L}_p^{p_2}(q,\2{s})|\le C_8d_0(p_1,p_2)\quad\forall 
p_1,p_2\in B_0(\2{p}_0,r_0),d_0(p_1,p_2)<\min (1,\2{s}/4)
$$
and the theorem follows.
\enddemo

By (1.6) and Theorem 1.14 we have the following theorem.

\proclaim{\bf Theorem 1.15}
Let $t_0>0$ and let $g$ and $\2{g}$ be related by (0.5). Suppose $(M,\2{g})$ 
satisfies (1.20) in $(0,t_0)$ for some constant $c_1>0$.
Then for any $\2{\tau}\in (0,t_0)$, $q\in M$, $L_p^{p_0}(q,\2{\tau})$ is 
locally Lipschitz in $p_0$ with respect to the metric $g(0)=\2{g}(t_0)$.
\endproclaim

$$
\text{Section 2}
$$

In this section we will generalize the $\Cal{L}$-exponential map of Perelman 
\cite{P1} and define the $\Cal{L}_p$-exponential map corresponding to the 
$\Cal{L}_p$-geodesic curve. We will also derive some elementary properties 
of the $\Cal{L}_p$-exponential map. 

We first start will a definition. Let $\2{\tau}>0$. For any $v\in T_{p_0}M$ 
let $\4{v}=v/(1-p)$. By Lemma 1.4 there exists a unique solution 
$\4{\gamma}_{\4{v}}(s)=\4{\gamma}(s;\4{v})$ of (1.10), (1.13), in $(0,s_0)$ 
for some $s_0>0$. Let $\gamma_v(\tau)=\gamma (\tau;v)=\4{\gamma}(s;\4{v})$ 
where $s$ and $\tau$ are related by (1.5). Then $\gamma_v$ is the
unique solution of (1.11), (1.18), in $(0,\tau_0)$ where 
$\tau_0=s_0^{1/(1-p)}$. Similar to \cite{Ye1} for any $\2{\tau}>0$ we let
$$
U_p(\2{\tau})=\{v\in T_{p_0}M:\gamma_v\text{ exists on }(0,\tau_0)
\text{ for some }\tau_0>\2{\tau}\}. 
$$
Then $U_p(\tau_2)\subset U_p(\tau_1)$ for any $0<\tau_1<\tau_2<t_0$.
We define the $\Cal{L}_p$-exponential map 
$\Cal{L}_p\text{-exp}_{p_0}^{\2{\tau}}:U_p(\2{\tau})\to M$ by
$$
\Cal{L}_p\text{-exp}_{p_0}^{\2{\tau}}(v)=\gamma_v(\2{\tau})=\4{\gamma}_{\4{v}}
(\2{\tau}^{1-p}).
$$
By O.D.E. theory and the equivalence of the O.D.E. (1.10), (1.13), and (1.11), 
(1.18), through the transformation (1.5), $U_p(\2{\tau})$ is open in 
$T_{p_0}M$. Note that by Corollary 1.5 if $(M,g)$ has uniformly bounded Ricci 
curvature on $M\times (0,t_0)$, then $U_p(\2{\tau})=T_{p_0}M$ for any 
$0\le\2{\tau}<t_0$.

Let $q_0\in M$, $\tau_0\in (0,t_0)$ and $v\in T_{q_0}M$. By an argument similar
to the proof of Lemma 1.4 and Corollary 1.5 there exists a unique solution
$\gamma_{\tau_0,v}^{q_0}(\tau)=\gamma_{\tau_0}^{q_0}(\tau;v)$ of (1.11) in 
$(\tau_0,\tau_1)$ for some $\tau_1>\tau_0$ such that
$$\left\{\aligned
&\gamma_{\tau_0}^{q_0}(\tau_0;v)=q_0\\
&\tau_0^p\gamma_{\tau_0}^{q_0}{}'(\tau_0;v)=v.\endaligned\right.
$$
For any $\2{\tau}>\tau_0$, let
$$
U_{\tau_0,p}^{q_0}(\2{\tau})=\{v\in T_{q_0}M:\gamma_{\tau_0,v}^{q_0}
\text{ exists on }(\tau_0,\tau_1)\text{ for some }\tau_1>\2{\tau}\}. 
$$
We define the $\Cal{L}_{\tau_0,p}^{q_0}$-exponential map 
$\Cal{L}_{\tau_0,p}^{q_0}\text{-exp}^{\2{\tau}}:U_{\tau_0,p}^{q_0}
(\2{\tau})\to M$ by
$$
\Cal{L}_{\tau_0,p}^{q_0}\text{-exp}^{\2{\tau}}(v)
=\gamma_{\tau_0}^{q_0}(\2{\tau};v).
$$

\proclaim{\bf Lemma 2.1}
Suppose $(M,g)$ satisfies (1.20) in $[0,t_0)$ for some constant $c_1>0$.
Then for any $r_0>0$ and $m_0>0$, there exists a constant $s_1\in 
(0,t_0^{1-p})$ such that for any $\4{v}\in T_{p_0}M$ satisfying 
$|\4{v}|_{g(p_0,0)}\le m_0$ there exists a unique $\4{\Cal{L}}_p$-geodesic 
$\4{\gamma}=\4{\gamma}_{\4{v}}:[0,s_1]\to M$ satisfying (1.13) and
$$
\4{\gamma}(s)\in B_0(p_0,r_0)\quad\forall 0\le s\le s_1.\tag 2.1
$$ 
Hence $\Cal{B}(0,(1-p)m_0)\subset U_p(\tau)$ for any $0<\tau\le\tau_1$ where 
$\tau_1=s_1^{\frac{1}{1-p}}$ and
$$
\bigcup_{0<\tau<t_0}U_p(\tau)=T_{p_0}M.
$$
\endproclaim
\demo{Proof}
We will use an argument similar to the proof of Proposition 2.5 of \cite{Ye1}
to prove the lemma. Let $\4{v}\in T_{p_0}M$ satisfy 
$|\4{v}|_{g(p_0,0)}\le m_0$.
Since $M$ is complete, $\2{B_0(p_0,r_0)}$ is compact. Then there exists a 
constant $K_1>0$ such that (1.16) holds for any $(q,\tau)\in\2{B_0(p_0,r_0)}
\times [0,t_0/2]$. Let $C_2=C_2(K_1)>0$ and $C_3>0$ be as in the proof of 
Lemma 1.4. Let $s_1'=\min (1,(t_0/2)^{1-p},e^{-c_1-(C_2/2)}(m_0^2+C_3)^{-1/2}
r_0/2)$ and $s_1=s_1'/2$. By Lemma 1.4 there exists a maximal interval 
$[0,\2{s})$ such that there exists a unique $\4{\Cal{L}}_p$-geodesic
$\4{\gamma}:[0,\2{s})\to M$ which satisfies (1.13). We claim that $\2{s}\ge 
s_1'$. Suppose not. Then $\2{s}<s_1'$. Let
$$
s_0=\sup\{s'\le\2{s}:\4{\gamma}(s)\in B_0(p_0,r_0)\quad\forall 
0\le s\le s'\}.
$$
Suppose $s_0<\2{s}$. By the same argument as the proof of Lemma 1.4
(1.17) holds in $(0,s_0)$. Hence by (1.17) and Lemma 1.6,
$$\align
&e^{-c_1}d_0(p_0,\4{\gamma}(s_0))\le e^{-c_1}\int_0^{s_0}
|\4{\gamma}'(s)|_{g(0)}\,ds\le\int_0^{s_0}|\4{\gamma}'(s)|\,ds
\le e^{C_2/2}(m_0^2+C_3)^{1/2}s_0\\
\Rightarrow\quad&d_0(p_0,\4{\gamma}(s_0))\le e^{c_1+(C_2/2)}(m_0^2+C_3)^{1/2}
s_1'<r_0.\endalign
$$
By continuity there exists $s_2\in (s_0,\2{s}]$ such that
$$
d_0(p_0,\4{\gamma}(s))<r_0\quad\forall 0\le s\le s_2.
$$
This contradicts the choice of $s_0$. Hence $s_0=\2{s}$. Then (1.17) holds 
on $[0,\2{s}]$. Thus by
(1.17) we can extend $\4{\gamma}$ to a solution of (1.10), (1.13), in 
$(0,\2{s}+\delta)$ for some $\delta\in (0,s_1'-\2{s})$. Contradiction arises.
Hence $\2{s}\ge s_1'$ and the lemma follows. 
\enddemo

\proclaim{\bf Theorem 2.2}
Suppose $(M,g)$ satisfies (1.20) in $[0,t_0)$ for some constant $c_1>0$. Then
there exists a constant $\2{\tau}_0\in (0,t_0)$ such that for any $0<\2{\tau}
\le\2{\tau}_0$ there exist a constant $r_1>0$ and an open set 
$O_1\subset M$ with $p_0\in O_1$ such that 
$\left.\Cal{L}_p\text{-exp}_{p_0}^{\2{\tau}}\right|_{\Cal{B}(0,r_1)}:
\Cal{B}(0,r_1)\to O_1$ is a diffeomorphism.
\endproclaim
\demo{Proof}
We will use a modification of the proof of Proposition 2.6 of \cite{Ye1} to
prove the lemma. Let $(\phi_0,B_0(p_0,r_0))$ be a local normal co-ordinate 
chart around $p_0$. By Lemma 2.1 there exists a constant 
$s_1\in (0,t_0^{1-p})$ such that $\Cal{B}(0,1)\subset U_p(\tau)$ for any 
$0<\tau\le\tau_1=s_1^{\frac{1}{1-p}}$ and (2.1) holds for any 
$\4{\Cal{L}}_p$-geodesic which satisfies (1.13) with $|\4{v}|\le 1$. 
By the inverse function theorem it suffices to check that 
the kernel of $d(\Cal{L}_p\text{-exp}_{p_0}^{\2{\tau}})_0$ is equal to 
zero for sufficiently small $\2{\tau}$. Suppose not. Then there exists 
$\2{\tau}\in (0,\2{\tau}_0)$ and $0\ne v_0\in T_0(T_{p_0}M)=T_{p_0}M$ such that
$$
d(\Cal{L}_p\text{-exp}_{p_0}^{\2{\tau}})_0(v_0)=0\tag 2.2
$$
where $\2{\tau}_0\in (0,\tau_1)$ is some constant to be determined later
in the proof.
Let $\4{v}_0=v_0/(1-p)$ and $\2{s}=\2{\tau}^{1-p}$. By rescaling $v_0$ if 
necessary we may assume 
without loss of generality that $|\4{v}_0|_{g(p_0,0)}=1$. Then
$|v_0|_{g(p_0,0)}=1-p$. For any $0\le z
\le 1$, let
$$
h(s,z)=\4{\gamma}(s;z\4{v_0})
$$
be the solution of (1.10), (1.13), in $[0,s_1]$ with $\4{v}$ being replaced
by $z\4{v_0}$ given by Lemma 2.1. Then $h$ is a variation 
of $\4{\gamma}(s;\4{v_0})$ with $h(s,1)=\4{\gamma}(s;\4{v_0})$.
Let
$$
Y(s,z)=\frac{\1 h}{\1 z}(s,z), \2{Y}(s)=Y(s,0)\text{ and }
\4{X}(s)=\frac{\1 h}{\1 s}(s,0).
$$
We write $\4{X}(s)=a^i(s)\1/\1 x_i$ and $\2{Y}(s)=b^j(s)\1/\1 x_j$ in the 
local co-odinates $(\phi_0,B_0(p_0,r_0))$. By (2.2),
$$\align
&0=\left.\frac{d}{dz}\right |_{z=0}\Cal{L}_p\text{-exp}_{p_0}^{\2{\tau}}(zv_0)
=\left.\frac{d}{dz}\right |_{z=0}\4{\gamma}(\2{s};z\4{v}_0)=Y(\2{s},0)
=\2{Y}(\2{s})\\
\Rightarrow\quad&b^j(\2{s})=0\quad\forall j=1,2,\dots,n.\tag 2.3
\endalign
$$
Note that
$$\align
&h(0,z)=\4{\gamma}(0;z\4{v}_0)=p_0\quad\forall 0\le z\le 1\\
\Rightarrow\quad&Y(0,z)=\frac{\1 h}{\1 z}(0,z)=0\quad\forall 0\le z\le 1\\
\Rightarrow\quad&b^j(0)\frac{\1}{\1 x_j}=\2{Y}(0)=0\\
\Rightarrow\quad&b^j(0)=0\quad\forall j=1,2,\dots,n\tag 2.4
\endalign
$$
and
$$
\nabla_s\2{Y}(0)=\frac{\1^2h}{\1 s\1 z}(0,0)=\frac{\1^2h}{\1 z\1 s}(0,0)
=\left.\frac{\1^2}{\1 z\1 s}\right |_{(0,0)}\4{\gamma}(s;z\4{v}_0)
=\left.\frac{\1 }{\1 z}\right |_{z=0}(z\4{v}_0)=\4{v}_0.\tag 2.5
$$
By an argument similar to the proof of Lemma 1.4 (1.17) holds for some
constants $C_2>0$, $C_3>0$, with $\4{v}=0$. Hence there exists a constant
$C>0$ such that
$$
|\4{X}(s)|\le C\quad\forall 0\le s\le s_1.
$$ 
Since $\4{\gamma}$ satisfies (1.10),
$$\align
&0=\nabla_s\biggl (\frac{\1 h}{\1 s}\biggr )
-\frac{1}{2(1-p)^2}s^{\frac{2p}{1-p}}\nabla\4{R}
+\frac{2}{1-p}s^{\frac{p}{1-p}}\4{\text{Ric}}\biggl (\frac{\1 h}{\1 s}_,\cdot
\biggl )\\
\Rightarrow\quad&0=\nabla_z\nabla_s\left (\frac{\1 h}{\1 s}\right )
-\frac{1}{2(1-p)^2}s^{\frac{2p}{1-p}}\nabla_z(\nabla\4{R})
+\frac{2}{1-p}s^{\frac{p}{1-p}}\nabla_z\left (\4{\text{Ric}}\left (
\frac{\1 h}{\1 s}_,\cdot\right )\right )\\
\Rightarrow\quad&0=\nabla_s\nabla_sY(s,z)
+\4{R}\left (\frac{\1 h}{\1 z}_,\frac{\1 h}{\1 s}\right )\frac{\1 h}{\1 s}
-\frac{1}{2(1-p)^2}s^{\frac{2p}{1-p}}\nabla_z(\nabla\4{R})\\
&\qquad +\frac{2}{1-p}s^{\frac{p}{1-p}}(\nabla_z\4{\text{Ric}})
\left (\frac{\1 h}{\1 s}_,\cdot\right )
+\frac{2}{1-p}s^{\frac{p}{1-p}}\4{\text{Ric}}(\nabla_sY,\cdot)\tag 2.6
\endalign
$$
where $\4{R}(X_1,X_2)X_3(q,s)=R(X_1,X_2)X_3(q,\tau)$ for any $X_1,X_2,X_3\in
T_qM$ with $s=\tau^{1-p}$. Putting $z=0$ we get
$$\align
0=&\nabla_s\nabla_s\2{Y}(s)+\4{R}(\2{Y},\4{X})\4{X}
-\frac{1}{2(1-p)^2}s^{\frac{2p}{1-p}}
\nabla_{\2{Y}(s)}(\nabla\4{R}(h(s,0),s))\\
&\qquad +\frac{2}{1-p}s^{\frac{p}{1-p}}(\nabla_{\2{Y}(s)}\4{\text{Ric}})
(\4{X},\cdot)+\frac{2}{1-p}s^{\frac{p}{1-p}}\4{\text{Ric}}(\nabla_sY,\cdot)
\quad\text{ in }(0,s_1)\\
=&\frac{d^2b^k}{ds^2}+\biggl (b^j\frac{db^i}{ds}
+b^i\frac{db^j}{ds}\biggr )\4{\Gamma}_{ij}^k
+b^ib^jb^m\frac{\1\4{\Gamma}_{ij}^k}{\1 x^m}
+b^m\biggl (\frac{db^r}{ds}+b^ib^j\4{\Gamma}_{ij}^r\biggr )
\4{\Gamma}_{mr}^k\\
&\qquad +g^{jk}b^i<\4{R}(\frac{\1}{\1x_i},\4{X})\4{X},\frac{\1}{\1x_j}>
-\frac{1}{2(1-p)^2}s^{\frac{2p}{1-p}}
g^{kj}b^i<\nabla_i(\nabla\4{R}),\frac{\1}{\1 x_j}>\\
&\qquad +\frac{2}{1-p}s^{\frac{p}{1-p}}g^{jk}a^ib^m\nabla_m\4{R}_{ij}
+\frac{2}{1-p}s^{\frac{p}{1-p}}g^{jk}\biggl (\frac{db^i}{ds}+b^mb^r
\4{\Gamma}_{mr}^i\biggr )\4{R}_{ij}
\tag 2.7\endalign
$$
holds in $(0,s_1)$ for all $k=1,2,\dots,n$. By (2.5),
$$
g_{ij}(p_0,0)\frac{db^i}{ds}(0)\frac{db^j}{ds}(0)=1.
$$
Hence
$$\align
&\lambda_1\sum_{k=1}^n\biggl (\frac{db^k}{ds}(0)\biggr )^2
\le g_{ij}(p_0,0)\frac{db^i}{ds}(0)\frac{\1 b^j}{\1 s}(0)
\le\lambda_2\sum_{k=1}^n\biggl (\frac{db^k}{ds}(0)\biggr )^2\\
\Rightarrow\quad &\lambda_1\sum_{k=1}^n\biggl (\frac{db^k}{ds}(0)
\biggr )^2\le 1
\le\lambda_2\sum_{k=1}^n\biggl (\frac{db^k}{ds}(0)\biggr )^2
\tag 2.8\endalign
$$
for some constant $\lambda_2>\lambda_1 >0$ depending only on $g_{ij}(p_0,0)$.
Let 
$$
E=\sum_{k=1}^n\biggl\{(b^k)^2+\biggl (\frac{db^k}{ds}\biggr )^2\biggr\}
$$ 
and
$$
s_2=\sup\{0<s_1'\le s_1:|b^i(s)|\le 1\quad\forall 0\le s\le s_1',i=1,2,
\dots,n\}.
$$
By (2.4) $s_2>0$. Then by (2.7) and (2.8),
$$\align
\biggl |\frac{dE}{ds}\biggr |=&2\biggl |\sum_{k=1}^n\biggl (
b^k\frac{db^k}{ds}+\frac{db^k}{ds}\frac{d^2b^k}{ds^2}\biggr )\biggr |\\
\le&\sum_{k=1}^n\biggl \{2|b^k|\biggl |\frac{db^k}{ds}\biggr |
+C\biggl |\frac{db^k}{ds}\biggr |
\biggl[\sum_{i=1}^n\left(|b^i|+\left|\frac{db^i}{ds}\right|\right)
+\sum_{i,j=1}^n\biggl (|b^ib^j|+|b^i|\left|\frac{db^j}{ds}\right|
\biggr )\\
&\qquad +\sum_{i,j,m=1}^n|b^ib^jb^m|\biggr]\biggr\}\\
\le&C_4E\quad\forall 0\le s\le s_2\\
\Rightarrow\qquad\qquad&\frac{d}{ds}(e^{-C_4s}E)\le 0\quad\forall 0
\le s\le s_2\\
\Rightarrow\qquad\qquad&E(s)\le e^{C_4s}E(0)
\le e^{C_4s_1}/\lambda_1=C_5\quad\forall 0\le s\le s_2.
\tag 2.9\endalign 
$$
for some constants $C_1>0$ and $C_4>0$ independent of $s_2$. We claim that
$$
s_2\ge\min (1/(2\sqrt{C_5}),s_1).
$$
Suppose not. Then $s_2<\min (1/(2\sqrt{C_5}),s_1)$. By (2.4) and (2.9),
$$
|b^i(s)|=\left|\int_0^s\frac{db^i}{ds}\,ds\right|\le\sqrt{C_5}s\le
\frac{1}{2}\quad\forall 0\le s\le s_2,i=1,2,\dots,n.
$$
Then by continuity there exists a constant $\delta>0$ such that
$$
|b^i(s)|\le 1\quad\forall 0\le s\le s_2+\delta,i=1,2,\dots,n.
$$
This contradicts the maximality of $s_2$. Hence 
$s_2\ge\min (1/(2\sqrt{C_5}),s_1)$. Let 
$$
s_3=\min (1/(2\sqrt{C_5}),s_1).
$$
Then (2.9) holds for all $0\le s\le s_3$.
By (2.8) there exists $i_0\in\{1,2,\dots,n\}$ such that 
$$
\left|\frac{db^{i_0}}{ds}(0)\right|\ge\frac{1}{\sqrt{n\lambda_2}}.
$$
By replacing $\4{v}_0$ by $-\4{v}_0$ if necessary and permutating the indices
we may assume without loss of generality that 
$$
\frac{db^1}{ds}(0)\ge\frac{1}{\sqrt{n\lambda_2}}.\tag 2.10
$$
By (2.7) and (2.9),
$$
\frac{d^2b^1}{ds^2}(s)+C_6\ge 0\quad\forall 0\le s\le s_3\tag 2.11
$$
for some constant $C_6>0$. Let
$$
\2{s}_0=\min (s_1,(2\sqrt{C_5})^{-1},(2\sqrt{n\lambda_2}C_6)^{-1}),
$$
$\2{\tau}_0=\2{s}_0^{1/(1-p)}$, and
$s_4=\sup\{s'\le\2{s}: db^1(s)/ds>0\quad\forall 0\le s<s'\}$. Then
$s_4\le\2{s}\le\2{s}_0$. Integrating (2.11), by (2.10) we 
have
$$
\frac{db^1}{ds}(s)\ge\frac{db^1}{ds}(0)-C_6s
\ge\frac{1}{\sqrt{n\lambda_2}}-C_6\2{s}_0
\ge\frac{1}{2\sqrt{n\lambda_2}}>0\quad\forall 0\le s\le s_4.\tag 2.12
$$
Suppose $s_4<\2{s}$. Then $s_4<\2{s}_0$. By (2.12) and continuity there exists
$\delta\in (0,\2{s}-s_4)$ such 
that $db^1(s)/ds>0$ on $(0,s_4+\delta)$. This contradicts the maximality of 
$s_4$. Hence $s_4=\2{s}$. Integrating (2.12) over $(0,\2{s})$,
$$
b^1(\2{s})>b^1(0)=0.
$$ 
This contradicts (2.3). Hence no such $v_0$ exists. Thus ker$(d(\Cal{L}_p
\text{-exp}_{p_0}^{\2{\tau}})_0)=0$ for any $0<\2{\tau}\le\2{\tau}_0$
and the theorem follows. 
\enddemo

By the proof of Theorem 2.2 it is natural to define the following:

\proclaim{\bf Definition 2.3}
Let $\4{\gamma}(s)$ be a $\4{\Cal{L}}_p$-geodesic in $(s_1,s_2)$. We say 
that a vector field $\4{Y}(s)$ along $\4{\gamma}$ is a 
$\4{\Cal{L}}_p$-Jacobi field in $(s_1,s_2)$ if $\4{Y}(s)$ satisfies
$$
\nabla_s\nabla_s\4{Y}+\4{R}(\4{Y},\4{\gamma}')\4{\gamma}'
-\frac{1}{2(1-p)^2}s^{\frac{2p}{1-p}}\nabla_{\4{Y}}(\nabla\4{R})
+\frac{2}{1-p}s^{\frac{p}{1-p}}\nabla_{\4{Y}}(\4{\text{Ric}}
(\4{\gamma}',\cdot ))=0\tag 2.13
$$
in $(s_1,s_2)$ along $\4{\gamma}$.
\endproclaim

\proclaim{\bf Definition 2.4}
Let $\gamma (\tau)$ be a $\Cal{L}_p$-geodesic in $(\tau_1,\tau_2)$ and let
$\4{\gamma}(s)$ be given by (1.5). We say that a vector field $Y(\tau)$ 
along $\gamma$ is a $\Cal{L}_p$-Jacobi field in $(\tau_1,\tau_2)$ if 
$\4{Y}(s)=Y(s^{1-p})$ is a $\4{\Cal{L}}_p$-Jacobi field in $(s_1,s_2)$
in $(s_1,s_2)$ along $\4{\gamma}$  where $s_i=\tau_i^{1-p}$, $i=1,2$.
\endproclaim

\proclaim{\bf Definition 2.5}
Let $\4{\gamma}(s)$ be a $\4{\Cal{L}}_p$-geodesic on $[0,\2{s}]$. For any 
$0\le s_0<s_1\le\2{s}$, we 
say that $\4{\gamma}(s_1)$ is $\4{\Cal{L}}_p$-conjugate to $\4{\gamma}(s_0)$
along $\left.\4{\gamma}\right|_{[s_0,s_1]}$ if there exists a 
$\4{\Cal{L}}_p$-Jacobi field $\4{Y}(s)\not\equiv 0$ 
along $\left.\4{\gamma}\right|_{[s_0,s_1]}$  such that $\4{Y}(s_0)
=\4{Y}(s_1)=0$.
\endproclaim

\proclaim{\bf Definition 2.6}
Let $\gamma (\tau)$ be a $\Cal{L}_p$-geodesic on 
$[0,\2{\tau}]$. For any $0\le\tau_0<\tau_1\le\2{\tau}$, we say that 
$\gamma(\tau_1)$ is $\Cal{L}_p$-conjugate to $\gamma (\tau_0)$ along 
$\left.\gamma\right|_{[\tau_0,\tau_1]}$ if $\4{\gamma}(s_1)$ is 
$\4{\Cal{L}}_p$-conjugate to $\4{\gamma}(s_0)$ along 
$\left.\4{\gamma}\right|_{[s_0,s_1]}$ where $\4{\gamma}(s)$ is given
by (1.5) with $s=\tau^{1-p}$ and $s_i=\tau_i^{1-p}$ for $i=0,1$.
\endproclaim

\proclaim{\bf Theorem 2.7}
Let $0<\2{\tau}<t_0$, $v,w\in T_{p_0}M$, and let $\gamma=\gamma 
(\tau;v):[0,\2{\tau}]\to M$ be a $\Cal{L}_p$-geodesic which satisfies 
(1.18). Let $\4{\gamma}(s)=\4{\gamma}(s;\4{v})$ be given by (1.5) with
$s=\tau^{1-p}$ where $\4{v}=v/(1-p)$. Suppose $Y(\tau)$ is a 
$\Cal{L}_p$-Jacobi field along $\gamma$ with $Y(0)=0$ and $\nabla_s\4{Y}(0)
=w/(1-p)$ where $\4{Y}(s)=Y(\tau)$ with $s=\tau^{1-p}$. Then
$$
Y(\tau)=d(\Cal{L}_p\text{-exp}_{p_0}^{\tau})_v(w)\quad\forall 0\le\tau
\le\2{\tau}.\tag 2.14
$$
\endproclaim
\demo{Proof}
Let $\alpha :(-\3_1,\3_1)\to T_{p_0}M$ be a curve in $T_{p_0}M$ such that 
$\alpha (0)=v$, $\alpha '(0)=w$. Let $\2{s}=\2{\tau}^{1-p}$. By continuous
dependence of solutions of O.D.E. on the initial data there exist $\3\in 
(0,\3_1)$ such that for any $\rho\in (-\3,\3)$ there exists a unique solution 
$\4{\gamma}(s;\alpha (\rho)/(1-p))$ of (1.10) on $[0,\2{s}]$
which satisfies (1.13) with $\4{v}$ being replaced by $\alpha (\rho)/(1-p)$.
For any $|\rho|<\3$, let $\gamma (\tau;\alpha (\rho))
=\4{\gamma}(s;\alpha (\rho)/(1-p))$ with $s=\tau^{1-p}$. Then $\gamma 
(\tau;\alpha (\rho))$ is the 
$\Cal{L}_p$-geodesic on $[0,\2{\tau}]$ which satisfies 
(1.18) with $v$ being replaced by $\alpha (\rho)$. Let
$$
h(\tau,\rho)=\Cal{L}_p\text{-exp}_{p_0}^{\tau}(\alpha (\rho))
=\gamma (\tau;\alpha (\rho))\quad\forall 0\le\tau\le\2{\tau},|\rho|<\3
$$
and 
$$
\4{Y_1}(s)=Y_1(\tau)=\frac{\1 h}{\1 \rho}(\tau,0)\quad\forall 
0\le\tau\le\2{\tau}, s=\tau^{1-p}.
$$
Then
$$
Y_1(\tau)=d(\Cal{L}_p\text{-exp}_{p_0}^{\tau})_v(w)
\quad\forall 0\le\tau\le\2{\tau}.
$$
Since $h(0,\rho)=p_0$ $\forall\rho\in (-\3,\3)$, 
$$
Y_1(0)=0=Y(0)\quad\Rightarrow\quad\4{Y}_1(0)=\4{Y}(0)=0.\tag 2.15
$$ 
Then 
$$\align
\nabla_s\4{Y}_1(0)=&\frac{\1^2h}{\1 s\1\rho}(0,0)
=\left.\frac{\1}{\1\rho}\right |_{\rho =0}\frac{\1 h}{\1 s}(0,\rho)
=\left.\frac{\1}{\1\rho}\right |_{\rho =0}\4{\gamma}'(0;\alpha (\rho)/(1-p))\\
=&\left.\frac{\1}{\1\rho}\right |_{\rho =0}\frac{\alpha (\rho)}{1-p}
=\frac{\alpha'(0)}{1-p}=\frac{w}{1-p}=\nabla_s\4{Y}(0).\tag 2.16
\endalign
$$
By an argument similar to the proof of Theorem 2.2 $\4{Y}_1(s)$ is a 
$\4{\Cal{L}}_p$-Jacobi field along $\4{\gamma}$. Since both $\4{Y}$ and 
$\4{Y}_1$ satisfies (2.13), (2.15), and (2.16) on $[0,\2{s}]$, by uniqueness 
of O.D.E. $\4{Y}\equiv\4{Y}_1$ on $[0,\2{s}]$.
Hence $Y(\tau)$ satisfies (2.14) and the theorem follows.
\enddemo

By Theorem 2.7 and an argument similar to the proof of Proposition 3.9
in Chapter 5 of \cite{C} we have 

\proclaim{\bf Theorem 2.8}
Let $0<\2{\tau}<t_0$, $v\in T_{p_0}M$, and let $\gamma=\gamma 
(\tau;v):[0,\2{\tau}]\to M$ be a $\Cal{L}_p$-geodesic which satisfies 
(1.18). If $\gamma (\2{\tau})$ is not $\Cal{L}_p$-conjugate to $p_0$, then 
for any $V_0\in T_{\gamma (\2{\tau})}M$ there exists a $\Cal{L}_p$-Jacobi 
field $Y(\tau)$ along $\gamma$ with $Y(0)=0$ and $Y(\2{\tau})=V_0$. 
\endproclaim

\proclaim{\bf Definition 2.9 (cf. Definition 4 of \cite{Ye1})}
For any $\2{\tau}\in (0,t_0)$, we define the injectivity domain 
$\Omega_p(\2{\tau})$ at time $\2{\tau}$ by 
$$\align
\Omega_p(\2{\tau})=\{&q\in M:\exists\text{ a unique }
\Cal{L}_p(q,\2{\tau})\text{-length minimizing }
\Cal{L}_p\text{-geodesic }\gamma:[0,\2{\tau}]\to M\\
&\text{ such that }\gamma (0)=p_0, \gamma (\2{\tau})=q,\text{ and }q
\text{ is not }\Cal{L}_p\text{-conjugate to }p_0\text{ along }\gamma\}. 
\endalign
$$
and we define the cut locus $C_p(\2{\tau})$ at time $\2{\tau}$ by 
$C_p(\2{\tau})=M\setminus\Omega(\2{\tau})$.
\endproclaim

\proclaim{\bf Definition 2.10}
For any $q_0\in M$, $0<\tau_0<\2{\tau}<t_0$, we define the 
$\Cal{L}_{\tau_0,p}^{q_0}$-injectivity domain $\Omega_{\tau_0,p}^{q_0}
(\2{\tau})$ at time $\2{\tau}$ by 
$$\align
\Omega_{\tau_0,p}^{q_0}(\2{\tau})=\{&q\in M:\exists\text{ a unique }
\Cal{L}_{\tau_0,p}^{q_0}(q,\2{\tau})\text{-length minimizing }
\Cal{L}_p\text{-geodesic }\gamma:[\tau_0,\2{\tau}]\to M\\
&\text{ such that }\gamma (\tau_0)=q_0, \gamma (\2{\tau})=q,\text{ and }q
\text{ is not }\Cal{L}_p\text{-conjugate to }q_0\text{ along }\gamma\}. 
\endalign
$$
and we define the $\Cal{L}_{\tau_0,p}^{q_0}$-cut locus $C_{\tau_0,p}^{q_0}
(\2{\tau})$ at time $\2{\tau}$ by $C_{\tau_0,p}^{q_0}(\2{\tau})=M\setminus
\Omega_{\tau_0,p}^{q_0}(\2{\tau})$.
\endproclaim

By the theory of ordinary Riemannian geometry and a similar argument as
the discussion on P.513 of \cite{Ye1} $L_p(q,\tau)$ is a smooth function
in $\cup_{\tau'>0}\Omega_p(\tau')\times\{\tau'\}$ and 
$L_{\tau_0,p}^{q_0}(q,\tau)$ is a smooth function in $\cup_{\tau'>\tau_0}
\Omega_{\tau_0,p}^{q_0}(\tau')\times\{\tau'\}$.

\proclaim{\bf Lemma 2.11}
Let $\2{\tau}\in (0,t_0)$. Suppose $(M,g)$ satisfies (1.20) 
for some constant $c_1>0$. Then for any $0<\rho\le\2{\tau}$, $L_p
(\cdot ,\rho)$ is locally Lipschitz in $M$ with respect to the metric 
$g(\rho)$. 
\endproclaim
\demo{Proof}
We will use a modification of the proof of Proposition 2.12 of \cite{Ye1} and 
the proof of Theorem 1.14 to prove the lemma. Let $0<\rho\le\2{\tau}$, 
$r_0>0$, and let $q_1,q_2\in B_0(p_0,r_0)$. By Theorem 1.11 for $i=1,2$, 
there exists $\Cal{L}_p(q_i,\rho)$-length minimizing $\Cal{L}_p$-geodesic 
$\gamma_i$, with $\gamma_i(0)=p_0$ and $\gamma_i(\rho)=q_i$. 
Let $\gamma :[0,d]\to M$ be a normalized minimizing geodesic with 
respect to the metric $g(0)$ with $\gamma (0)=q_1$, $\gamma (d)=q_2$, 
$|\gamma'|=|\gamma'|_{g(0)}=1$ on $[0,d]$ where $d=d_0(q_1,q_2)$. 
Then $\gamma([0,d])\subset B_0(p_0,2r_0)$. Let
$$
K_0=\sup_{\2{B_0(p_0,2r_0)}\times [0,\2{\tau}]}(|R|+|\text{Ric}|).
$$
Then $K_0<\infty$. For $i=1,2$, let $d_i=d_0(p_0,q_i)$ and let $\2{\gamma}_i:
[0,d_i]\to M$ be a normalized minimizing geodesic with respect to the metric 
$g(0)$ with $\2{\gamma}_i(0)=p_0$, $\2{\gamma}_i(d_i)=q_i$. Then $d_i
<r_0$ and $\2{\gamma}_i([0,d_i])\subset B_0(p_0,2r_0)$ for $i=1,2$. By 
Lemma 1.7 and the proof of Lemma 1.8, there exist constants $A_1
=A_1(\2{\tau},r_0,K_0)>0$ and $r_1=r_1(\2{\tau},r_0,K_0)\ge 2r_0$ such that
$$\left\{\aligned
&L_p(q_i,\rho)\le A_1\qquad\forall i=1,2\\
&d_0(p_0,\gamma_i(\tau))<r_1\quad\forall 0\le\tau\le\rho,i=1,2.
\endaligned\right.
$$
Let
$$
K_1=\sup_{\2{B_0(p_0,2r_1)}\times [0,\2{\tau}]}(|R|+|\nabla R|
+|\text{Ric}|)
$$
and let $\4{\gamma}_i(s)=\gamma_i(\tau)$ with $s=\tau^{1-p}$ for $i=1,2$. Then
$\4{\gamma}_1$, $\4{\gamma}_2$, are $\4{\Cal{L}}_p$-geodesics. We assume now 
$d=d_0(q_1,q_2)<\rho/4$. Similar to the proof of Proposition 2.12 of 
\cite{Ye1} we let 
$$
\beta(\tau)=\left\{\aligned
&\gamma_1(\tau)\qquad\qquad\quad\,\,\text{ if }0\le\tau\le\rho -2d\\
&\gamma_1(2\tau -\rho+2d)\quad\text{ if }\rho-2d\le\tau\le\rho -d\\
&\gamma (\tau -\rho +d)\qquad\,\,\,\text{ if }\rho-d\le\tau\le\rho.
\endaligned\right.
$$
Then by Lemma 1.6,
$$\align
L_p(q_2,\rho)\le&\Cal{L}_p(q_2,\beta ,\rho)\\
\le&L_p(q_1,\rho)-\int_{\rho-2d}^{\rho}\tau^pR(\gamma_1(\tau),\tau)\,d\tau
+\int_{\rho-2d}^{\rho-d}\tau^p[R(\gamma_1(\tau'),\tau)
+4|\gamma_1'(\tau')|_{g(\tau)}^2]\,d\tau\\
&\qquad +\int_{\rho-d}^{\rho}\tau^p[R(\gamma (\tau''),\tau)
+|\gamma'(\tau'')|_{g(\tau)}^2]\,d\tau\\
\le&L_p(q_1,\rho)+(p+1)^{-1}[2K_1(\rho^{p+1}-(\rho-2d)^{p+1})+
e^{2K_1\rho}(\rho^{p+1}-(\rho-d)^{p+1})]\\
&\qquad +2e^{2K_1\rho}\int_{\rho-2d}^{\rho}\tau^p|\gamma_1'(\tau)|^2\,d\tau.
\tag 2.17
\endalign
$$
where $\tau'=2\tau -\rho+2d$ and $\tau''=\tau -\rho +d$. By the same argument 
as the proof of Lemma 1.4, (1.17) holds in $(0,\2{s}_{\rho})$, $\2{s}_{\rho}
=\rho^{1-p}$, for some constant $C_2>0$, $C_3>0$, 
depending only on $K_1$ with $\4{X}(s)$ and $\4{v}$, being replaced by 
$\4{\gamma}_i'(s)$ and $\4{v}_i=\4{\gamma}_i'(0)$, for $i=1,2$. Since 
$\4{\gamma}_i'(s)=\tau^p\gamma_i'(\tau)/(1-p)$, by (1.17) there exist constants
$C_4>0$, $C_5>0$ and $C_6>0$ such that
$$
C_4|v_i|^2-C_5\le\tau^{2p}|\gamma_i'(\tau)|^2\le C_6(1+|v_i|^2)
\quad\forall 0\le\tau\le\rho,i=1,2.\tag 2.18
$$
By (2.18),
$$\align 
C_4\rho |v_i|^2=&C_4\int_0^{\rho}|v_i|^2\,d\tau\le C_5\rho
+\int_0^{\rho}\tau^{2p}|\gamma_i'(\tau)|^2\,d\tau\\
\le&C_5\rho+\rho^p\biggl (L_p(q_i,\rho)+K_1\int_0^{\rho}\tau^{p}\,d\tau
\biggr )\\
\le&C_5\rho+\rho^p(A_1+(K_1\rho^{p+1}/(p+1)))=C_7\quad\forall i=1,2.
\tag 2.19\endalign
$$
By (2.18) and (2.19),
$$\align
\int_{\rho-2d}^{\rho}\tau^p|\gamma_i'(\tau)|^2\,d\tau
\le&C_6(1+|v_i|^2)\int_{\rho-2d}^{\rho}\tau^{-p}\,d\tau
\le 2C_6(1+|v_i|^2)d/(\rho-2d)^p\\
\le&4dC_6(1+(C_7/C_4\rho))/\rho^p\quad\forall i=1,2.\tag 2.20
\endalign
$$
By (2.17) and (2.20) there exists a constant 
$C_8=C_8(\rho,\2{\tau},r_0)>0$ such that
$$
L_p(q_2,\rho)\le L_p(q_1,\rho)+C_8d_0(q_1,q_2)
\quad\forall q_1,q_2\in B_0(p_0,r_0), d_0(q_1,q_2)<\rho/4.
$$
Hence by Lemma 1.6,
$$
L_p(q_2,\rho)\le L_p(q_1,\rho)+C_9d_{\rho}(q_1,q_2)
\quad\forall q_1,q_2\in B_0(p_0,r_0), d_0(q_1,q_2)<\rho/4
$$
for some constant $C_9>0$.
Interchanging the role of $q_1$ and $q_2$ in the above inequality,
$$
L_p(q_1,\rho)\le L_p(q_2,\rho)+C_9d_{\rho}(q_1,q_2)\quad\forall 
q_1,q_2\in B_0(p_0,r_0),d_0(q_1,q_2)<\rho/4.
$$
Hence
$$
|L_p(q_1,\rho)-L_p(q_2,\rho)|\le C_9d_{\rho}(q_1,q_2)\quad\forall 
q_1,q_2\in B_0(p_0,r_0),d_0(q_1,q_2)<\rho/4
$$
and the lemma follows.
\enddemo

By an argument similar to the proof of Proposition 2.13 of \cite{Ye1} we have

\proclaim{\bf Lemma 2.12}
Let $\2{\tau}\in (0,t_0)$. Suppose $(M,g)$ satisfies (1.20) in $(0,\2{\tau})$
for some constant $c_1>0$. Then for any $q\in M$, $L_p(q,\cdot)$ is locally 
Lipschitz in $(0,\2{\tau}]$.
\endproclaim

\proclaim{\bf Lemma 2.13}
Suppose $(M,g)$ satisfies (1.20) in $(0,t_0)$ for some constant $c_1>0$. Then 
$\Omega_p(\2{\tau})$ is open in $M$ for any $\2{\tau}\in (0,t_0)$ and 
$$
\cup_{0<\tau<t_0}\Omega_p(\tau)\times\{\tau\}
$$ 
is open in $M\times (0,t_0)$ with respect to the product metric $g\,dx^2\oplus 
d\tau^2$. Hence $C_p(\2{\tau})$ is close in $M$ for any $\2{\tau}\in (0,t_0)$ 
and $\cup_{0<\tau<t_0}C_p(\tau)\times\{\tau\}$ is closed in $M\times (0,t_0)$ 
with respect to the product metric $g\,dx^2\oplus d\tau^2$.
\endproclaim
\demo{Proof}
Let $\2{\tau}\in (0,t_0)$ and $q\in\Omega_p(\2{\tau})$. Let $\gamma (\tau ;v)$ 
be the minimizing $\Cal{L}_p(q,\2{\tau})$-geodesic given by Theorem 1.11
which satisfies (1.11) and (1.18) for some $v\in T_{p_0}M$.
Since $q$ is not $\Cal{L}_p$-conjugate to $p_0$, by Theorem 2.7,
$$
\text{ker}(d(\Cal{L}_p\text{-exp}_{p_0}^{\2{\tau}})_v)=0\quad
\Rightarrow\quad\text{det}(d(\Cal{L}_p\text{-exp}_{p_0}^{\2{\tau}})_v)\ne 0.
$$
By the inverse function theorem there exist $\3\in (0,\min (\2{\tau},
t_0-\2{\tau})/2)$, $\Cal{B}(v,r_0)$ and an 
open neighbourhood $O(q)$ of $(q,\2{\tau})$ in $M\times (0,t_0)$ such that 
the map 
$$
\phi :\Cal{B}(v,r_0)\times (\2{\tau}-\3,\2{\tau}+\3)\to O(q)
$$
given by $\phi (v,\tau)=(\Cal{L}_p\text{-exp}_{p_0}^{\tau}(v),\tau)$
is a differeomorphism for any $\tau\in (\2{\tau}-\3,\2{\tau}+\3)$ and
$$
\text{det}(d(\Cal{L}_p\text{-exp}_{p_0}^{\tau})_{v'})\ne 0\quad\forall
|\tau -\2{\tau}|<\3,v'\in \Cal{B}(v,r_0).\tag 2.21
$$
We claim that there exists $B_0(q,r_1)\times\{\2{\tau}\}\subset O(q)$ such 
that $B_0(q,r_1)\subset\Omega_p(\2{\tau})$. Suppose not. Then there exists a 
sequence of points $\{q_i\}_{i=1}^{\infty}$, $q_i\not\in\Omega_p(\2{\tau})\quad
\forall i\in\Cal{Z}^+$, such that $q_i\to q$ as $i\to\infty$. By the proof of 
Lemma 2.11, there exist a constant $C_1>0$ and $i_0\in\Bbb{Z}^+$ such that
$$\align
&L_p(q_i,\2{\tau})\le L_p(q,\2{\tau})+C_1d_0(q_i,q)\quad\forall
i\ge i_0\\
\Rightarrow\quad&\exists C_2>0\text{ such that }L_p(q_i,\2{\tau})\le C_2\quad
\forall i\in\Bbb{Z}^+.\tag 2.22
\endalign
$$
Now by Theorem 1.11 for any $i=1,2,\dots$, either 
$$\align
\text{(i) }&q_i\text{ is }\Cal{L}_p\text{-conjugate to }p_0
\text{ along some }\Cal{L}_p(q_i,\2{\tau})\text{-length minimizing }
\Cal{L}_p\text{-geodesic }\gamma (\cdot ;v_i)\\
&\text{satisfying }\gamma (0;v_i)=p\text{ and }\gamma (\2{\tau};v_i)=q_i
\endalign
$$ 
or 
$$\align
\text{(ii) }&\text{there exists two }\Cal{L}_p(q_i,\2{\tau})
\text{-length minimizing }\Cal{L}_p\text{-geodesics }
\gamma (\cdot ;v_i)\text{ and }\gamma (\cdot ;v_i')\text{ satisfying }\\
&\gamma (0;v_i)=\gamma (0;v_i')=p,\gamma (\2{\tau};v_i)=\gamma (\2{\tau};v_i')
=q_i,\text{ with }v_i\ne v_i'
\endalign
$$
where $\gamma (\cdot;w)$ stands for the solution of (1.11) and (1.18)
with $v$ being replaced by $w$.
Then either (i) or (ii) holds for infinitely many $i\in\Bbb{Z}^+$.
We now divide the proof into two cases:

\noindent $\underline{\text{\bf Case 1}}$: (i) holds for infinitely many 
$i\in\Bbb{Z}^+$.  

Without loss of generality we may assume that (i) holds for all 
$i\in\Bbb{Z}^+$. Then by Theorem 2.7,
$$
\text{ker}(d(\Cal{L}_p\text{-exp}_{p_0}^{\2{\tau}})_{v_i})\ne 0\quad\forall
i\in\Bbb{Z}^+\quad
\Rightarrow\quad\text{det}(d(\Cal{L}_p\text{-exp}_{p_0}^{\2{\tau}})_{v_i})=0
\quad\forall i\in\Bbb{Z}^+.\tag 2.23
$$
By (2.22) and an argument similar to the proof of (2.19) there exists
a constant $r_2>0$ such that
$$
|v_i|\le r_2\quad\forall i=1,2,\dots.\tag 2.24
$$
Since the closed ball $\2{\Cal{B}(0,r_2)}$ is compact in $T_{p_0}M$ with 
respect to the metric $g(p_0,0)$,
the sequence $\{v_i\}_{i=1}^{\infty}$ has a subsequence which we may assume 
without loss of generality to be the sequence itself such that $v_i\to v_0$ 
as $i\to\infty$ for some $v_0\in T_{p_0}M$. By continuous dependence of 
solutions of O.D.E. on the initial data $\gamma (\tau ;v_i)$ will 
converge uniformly to a $\Cal{L}_p$-geodesic $\gamma (\tau ;v_0)$ on 
$[0,\2{\tau}]$ as $i\to\infty$. By Fatou's Lemma and (2.22),
$$\align
&L_p(q,\2{\tau})
\le\Cal{L}_p(q,\gamma (\cdot ;v_0),\2{\tau})\le\lim_{i\to\infty}
\Cal{L}_p(q_i,\gamma_i(\cdot;v_i),\2{\tau})
=\lim_{i\to\infty}L_p(q_i,\2{\tau})
\le L_p(q,\2{\tau})\\
\Rightarrow\quad&L_p(q,\2{\tau})=\Cal{L}_p(q,\gamma (\cdot ;v_0),\2{\tau}).
\endalign
$$
Hence $\gamma (\cdot ;v_0)$ is a minimizing $\Cal{L}_p(q,\2{\tau})$-geodesic.
Since $q\in\Omega(\2{\tau})$, the minimizing 
$\Cal{L}_p(q,\2{\tau})$-geodesic is unique. Hence $v_0=v$. Letting $i\to\infty$
in (2.23),
$$
\text{det}(d(\Cal{L}_p\text{-exp}_{p_0}^{\2{\tau}})_{v})=0.
$$
This contradicts (2.21). Hence case 1 does not hold.

\noindent $\underline{\text{\bf Case 2}}$: (ii) holds for infinitely many 
$i\in\Bbb{Z}^+$.  

Without loss of generality we may assume that (ii) holds for all 
$i\in\Bbb{Z}^+$. By the same argument as case 1 there exists $r_2>0$ such that
$v_i,v_i'\in\2{\Cal{B}(0,r_2)}$ for all $i\in\Bbb{Z}^+$. Then as in case 1
by choosing a subsequence if necessary we may assume without loss of 
generality that $v_i\to v$ and $v_i'\to v$ as $i\to\infty$. Then there 
exists $i_0\in\Bbb{Z}^+$ such that $v_i,v_i'\in\Cal{B}(v,r_0)$ for all
$i\ge i_0$. Since the map $\phi$ is a differeomorphism,
$$
\Cal{L}_p\text{-exp}_{p_0}^{\2{\tau}}(v_i)
=\Cal{L}_p\text{-exp}_{p_0}^{\2{\tau}}(v_i')=q_i\quad\forall i\ge i_0\qquad
\Rightarrow\quad v_i=v_i'\quad\forall i\ge i_0.
$$
Contradiction arise. Hence case 2 does not hold. Thus no such sequence 
$\{q_i\}_{i=1}^{\infty}$ exists. Hence there exists $B_0(q,r_1)\times
\{\2{\tau}\}\subset O(q)$ such that $B_0(q,r_1)\subset\Omega_p(\2{\tau})$. 
Therefore $\Omega_p(\2{\tau})$ is open. By a similar argument 
$\cup_{0<\tau<t_0}\Omega_p(\tau)\times\{\tau\}$ is open in $M\times (0,t_0)$ 
and the lemma follows.
\enddemo

By Theorem 1.11, Lemma 2.13, and an argument similar to the proof of 
Proposition 2.16 of \cite{Ye1} but with Lemma 2.11 replacing Proposition 2.12 
in the proof there we have

\proclaim{\bf Lemma 2.14}
Suppose $(M,g)$ satisfies (1.20) in $(0,t_0)$ for some constant $c_1>0$. 
Then $C_p(\2{\tau})$ is a closed set of measure zero for any $\2{\tau}\in 
(0,t_0)$.
\endproclaim

\proclaim{\bf Lemma 2.15}
Suppose $(M,g)$ satisfies (1.20) in $(0,t_0)$ for some constant $c_1>0$. Then 
$\cup_{0<\tau<t_0}C_p(\tau)\times\{\tau\}$ is a closed set of measure zero in 
$M\times (0,t_0)$ with respect to the product metric $g\,dx^2\oplus d\tau^2$.
\endproclaim
\demo{Proof}
This result for the case $p=1/2$ is stated in \cite{Ye1}. We
will give a proof of it here for any $0<p<1$. By Lemma 2.13 we know that
$\cup_{0<\tau<t_0}C_p(\tau)\times\{\tau\}$ is closed in $M\times (0,t_0)$.
It suffices to show that $\cup_{\tau_1\le\tau\le\tau_2}(C_p(\tau)\cap
B_0(p_0,r_0))\times\{\tau\}$ has measure zero in $M\times (0,t_0)$ 
for any $0<\tau_1<\tau_2<t_0$ and $r_0>0$. 

Let $\2{\tau}\in [\tau_1,\tau_2]$, $\delta>0$, $\Omega_p(\tau,r_0)=\Omega_p 
(\tau)\cap B_0(p_0,r_0)$, $D_p(r_0)=\cup_{0<\tau<t_0}\Omega_p(\tau,r_0)
\times\{\tau\}$, and
$$
C_p(\tau_1,\tau_2,r_0)
=\cup_{\tau_1\le\tau\le\tau_2}(C_p(\tau)\cap B_0(p_0,r_0))\times\{\tau\}.
$$ 
We choose a compact set $K(\2{\tau})\subset\Omega_p(\2{\tau},r_0)$ such 
that 
$$
m_{\2{\tau}}(\Omega (\2{\tau},r_0)\setminus K(\2{\tau}))<\delta.
$$ 
Then by Lemma 2.14,
$$
m_{\2{\tau}}(B_0(p_0,r_0)\setminus K(\2{\tau}))<\delta.\tag 2.25
$$
Since $D_p(r_0)$ is open, for any $q\in K(\2{\tau})$ there exists $\3_q>0$
and an open ball $O_{\2{\tau}}(q)\subset\Omega (\2{\tau})$ containing $q$ 
such that $O_{\2{\tau}}(q)\times [\2{\tau}-\3_q,\2{\tau}+\3_q]\subset 
D_p(r_0)$. Since $K(\2{\tau})$ is compact, 
$$
K(\2{\tau})\subset\cup_{i=1}^{n(\2{\tau})}O_{\2{\tau}}(q_i)
$$
for some $q_1,\dots,q_{n(\2{\tau})}\in K(\2{\tau})$. Let $\3_{\2{\tau}}
=\min_{1\le i\le n_0}\3_{q_i}$. Since $[\tau_1,\tau_2]$ is compact, there 
exists $\tau_1\le\2{\tau}_1<\2{\tau}_2<\dots<\2{\tau}_{k_0}\le\tau_2$ such 
that
$$
[\tau_1,\tau_2]\subset\cup_{k=1}^{k_0}
(\2{\tau}_k-\3_k,\2{\tau}_k+\3_k)
$$ 
where $\3_k=\3_{\2{\tau}_k}$ for all $k=1,2,\dots,k_0$. Let $I_1=
(\2{\tau}_1-\3_1,\2{\tau}_1+\3_1)\cap [\tau_1,\tau_2]$ and 
$I_k=((\2{\tau}_k-\3_k,\2{\tau}_k+\3_k)\setminus\cup_{j=1}^{k-1}I_j)\cap 
[\tau_1,\tau_2]$ for all $k=2,3,\dots,k_0$. Then
$$
[\tau_1,\tau_2]=\cup_{k=1}^{k_0}I_k
$$ 
For any $k=1,2,\dots,k_0$, let $E_k=(\cup_{i=1}^{n_k}O_{\2{\tau}_k}(q_i))
\times I_k$ where $n_k=n(\2{\tau}_k)$. Then $\cup_{k=1}^{k_0}E_k\subset 
D_p(r_0)$ and
$$\align
C_p(\tau_1,\tau_2,r_0)\subset&(B_0(p_0,r_0)\times [\tau_1,\tau_2])\setminus
\cup_{k=1}^{k_0}E_k\subset\cup_{k=1}^{k_0}( B_0(p_0,r_0)\setminus
K(\2{\tau}_k))\times I_k\\
\Rightarrow\quad m(C_p(\tau_1,\tau_2,r_0))
\le&m\left (\cup_{k=1}^{k_0}( B_0(p_0,r_0)\setminus
K(\2{\tau}_k))\times I_k\right)\\
\le&\sum_{k=1}^{k_0}m\left((B_0(p_0,r_0)\setminus K(\2{\tau}_k))\times 
I_k\right).\tag 2.26
\endalign
$$
Let $C_1=\sup_{\2{B_0(p_0,r_0)}\times [\tau_1,\tau_2]}R(q,\tau)$. Then 
$C_1<\infty$. Note that by \cite{H5} the volume form $\sqrt{g}$ of $M$
satisfies
$$\align
&\left|\frac{d\sqrt{g}}{d\tau}\right|=|R\sqrt{g}|\le C_1\sqrt{g}\quad
\text{ in }\2{B_0(p_0,r_0)}\times [\tau_1,\tau_2]\tag 2.27\\
\Rightarrow\quad&\sqrt{g}(q,\tau)\le e^{C_1\3_k}\sqrt{g}(q,\2{\tau}_k)
\quad\forall q\in\2{B_0(p_0,r_0)},\tau\in I_k,k=1,2,\dots,n.\tag 2.28
\endalign
$$
By (2.25), (2.26), and (2.28),
$$
m(C_p(\tau_1,\tau_2,r_0))\le\sum_{k=1}^{k_0}\left\{e^{C_1\3_k}|I_k|
m_{\2{\tau}_k}(B_0(p_0,r_0)\setminus K(\2{\tau}_k))\right\}
\le e^{C_1}\delta\sum_{k=1}^{k_0}|I_k|\le e^{C_1}(\tau_2-\tau_1)\delta.
$$
Letting $\delta\to 0$,
$$
m(C_p(\tau_1,\tau_2,r_0))=0
$$
and the lemma follows.
\enddemo

By Lemma 2.11 and Lemma 2.12 we have

\proclaim{\bf Lemma 2.16}
Suppose $(M,g)$ satisfies (1.20) in $(0,t_0)$ for some constant $c_1>0$. 
Then $|\nabla L_p(q,\tau)|$ and $\1 L_p(q,\tau)/\1\tau$ are locally bounded 
measurable functions on $M\times (0,t_0)$. 
\endproclaim

$$
\text{Section 3}
$$

In this section we will prove the second variation formula for 
$L_p(q,\2{\tau})$. We will prove various properties of the 
$L_p(q,\tau)$-length, the generalized reduced distance $l_p$, and the 
generalized reduced volume $\4{V}_p(\tau)$. We will now assume that $(M,g)$ 
satisfies (1.20) in $(0,t_0)$ for some constant $c_1>0$ for the
rest of the paper. For any $\2{\tau}\in (0,t_0)$,
let  
$$\align
U_p'(\2{\tau})=&\{v\in U_p(\2{\tau}):\gamma_v(\2{\tau})\in\Omega_p(\2{\tau})
\text{ where }\gamma_v(\cdot)=\gamma (\cdot;v):[0,\2{\tau}]\to M\text{ is the }
\Cal{L}_p\text{-geodesic}\\
&\text{ that satisfies (1.18)}\}.
\endalign
$$
Note that by the definition of $\Omega_p(\2{\tau})$ and Theorem 1.11,
$$
\left.\Cal{L}_p\text{-exp}_{p_0}^{\2{\tau}}\right |_{U_p'(\2{\tau})}:
U_p'(\2{\tau})\to\Omega_p(\2{\tau})
$$
is a diffeomorphism. For any $v\in U_p(\2{\tau})$, let $J_p(v,\2{\tau})$ be 
the Jacobian of the $\Cal{L}_p\text{-exp}_{p_0}^{\2{\tau}}$ map at $v$. Let
$$
\Omega_p=\cup_{0<\tau<t_0}\Omega_p(\tau)\times\{\tau\}.
$$
By the same argument as the discussion on P.518 of \cite{Ye1} $L_p(q,\tau)$ is 
a smooth function in $\Omega_p$. If $\tau\in (0,t_0)$ and $q\in \Omega_p
(\tau)$, then there exists a unique $\Cal{L}_p(q,\tau)$-length minimizing 
$\Cal{L}_p$-geodesic $\gamma$ satisfying $\gamma (0)=p_0$, $\gamma (\tau)=q$,
such that $q$ is not $\Cal{L}_p$-conjugate to $p_0$. Then by Lemma 1.1,
$$
\nabla L_p(q,\tau)=2\tau^p\gamma '(\tau).\tag 3.1
$$

\proclaim{\bf Lemma 3.1}
Let $(q,\2{\tau})\in\Omega_p$ and let $\gamma$ be a $\Cal{L}_p$-geodesic 
satisfying $\gamma (0)=p_0$, $\gamma (\2{\tau})=q$, which minimizes the 
$\Cal{L}_p(q,\2{\tau})$-length. Suppose $Y$ is as in Lemma 1.1. Then
$$\align
&\delta_Y^2 L_p(q,\2{\tau})\\
\le &2\2{\tau}^p<X(\2{\tau}),\nabla_YY(\2{\tau})>\\
&\qquad +\int_0^{\2{\tau}}\tau^p\{\text{Hess}_R(Y,Y)+2<R(Y,X)Y,X>
+2|\nabla_XY|^2+2\nabla_XRic(Y,Y)\\
&\qquad -4\nabla_YRic(X,Y)\}\,d\tau\tag 3.2
\endalign
$$
where $X=X(\tau)=\gamma '(\tau)$.
\endproclaim
\demo{Proof}
We will use a modification of the argument of section 7 of \cite{P1} to 
prove the theorem. Let $f:[0,\2{\tau}]\times (-\3,\3)\to M$  be as in the 
proof of Lemma 1.1. Since 
$$\left\{\aligned
&L_p(f(\2{\tau},z),\2{\tau})\le\Cal{L}_p(f(\2{\tau},z),f(\cdot,z),\2{\tau})
\quad\forall |z|<\3\\
&L_p(f(\2{\tau},0),\2{\tau})=\Cal{L}_p(f(\2{\tau},0),f(\cdot,0),\2{\tau}),
\endaligned\right.
$$
differentiating (1.2) with respect to $z$ and putting $z=0$,
$$\align
&\delta_Y^2 L_p(q,\2{\tau})\\
=&\left.\frac{d^2}{dz^2}\right|_{z=0}L_p(f(\2{\tau},z),\2{\tau})
\le\left.\frac{d^2}{dz^2}\right|_{z=0}\Cal{L}_p(f(\2{\tau},z),f(\cdot,z),
\2{\tau})\\
\le &\int_0^{\2{\tau}}\tau^p(Y(Y(R))+2|\nabla_XY|^2+2<X,\nabla_Y\nabla_XY>)\,
d\tau\\
\le&\int_0^{\2{\tau}}\tau^p\{Y(Y(R))+2|\nabla_XY|^2+2<X,\nabla_X\nabla_YY>+
2<R(Y,X)Y,X>\}\,d\tau.\tag 3.3\endalign
$$
Since
$$\align
<\nabla_YY,X>=&Y(<Y,X>)-<Y,\nabla_YX>=Y(<Y,X>)-<Y,\nabla_XY>\\
=&Y(<Y,X>)-\frac{1}{2}X(<Y,Y>),\endalign
$$
we have
$$\align
&\frac{d}{d\tau}<\nabla_YY,X>\\
=&X(<\nabla_YY,X>)+\frac{\1}{\1\tau}<\nabla_YY,X>\\
=&<\nabla_X\nabla_YY,X>+<\nabla_YY,\nabla_XX>+Y(\frac{\1 g}{\1\tau}(Y,X))
-\frac{1}{2}X(\frac{\1 g}{\1\tau}(Y,Y))\\
=&<\nabla_X\nabla_YY,X>+<\nabla_YY,\nabla_XX>+2Y(\text{Ric}(Y,X))
-X(\text{Ric}(Y,Y))\\
=&<\nabla_X\nabla_YY,X>+<\nabla_YY,\nabla_XX>
+2\text{Ric}(\nabla_YY,X)+2\nabla_Y\text{Ric}(X,Y)\\
&\qquad -(\nabla_X\text{Ric})(Y,Y).\tag 3.4
\endalign
$$
Note that (3.4) is stated in section 7 of \cite{P1} but there is no proof of 
it in \cite{P1}. We refer the reader to \cite{KL} for another proof of (3.4)
by B.~Klein and J.~Lott. By (3.4),
$$\align
&2\int_0^{\2{\tau}}\tau^p<X,\nabla_X\nabla_YY>\,d\tau\\
=&2\int_0^{\2{\tau}}\tau^p\biggl\{\frac{d}{d\tau}<X,\nabla_YY>
-<\nabla_XX,\nabla_YY>-2\text{Ric}(\nabla_YY,X)-2\nabla_Y\text{Ric}(X,Y)\\
&\qquad +(\nabla_X\text{Ric})(Y,Y)\biggr\}\,d\tau\\
=&2\2{\tau}^p<X(\2{\tau}),\nabla_YY(\2{\tau})>
-2\int_0^{\2{\tau}}\tau^p\biggl\{\frac{p}{\tau}<\nabla_YY,X>
+<\nabla_YY,\frac{1}{2}\nabla R-\frac{p}{\tau}X-2\text{Ric} (X,\cdot)>\\
&\qquad +2\text{Ric}(\nabla_YY,X)+2\nabla_Y\text{Ric}(X,Y)
-\nabla_X\text{Ric}(Y,Y)\biggr\}\,d\tau\\
=&2\2{\tau}^p<X(\2{\tau}),\nabla_YY(\2{\tau})>+\int_0^{\2{\tau}}\tau^p\{
-(\nabla_YY)R
-4\nabla_Y\text{Ric}(X,Y)+2\nabla_X\text{Ric}(Y,Y)\}\,d\tau.\tag 3.5
\endalign
$$
By (3.3) and (3.5), (3.2) follows.
\enddemo

\proclaim{\bf Lemma 3.2}
Let $(q,\2{\tau})\in\Omega_p$ and let $\gamma$, $X$, be as in Lemma 3.1. 
Let $b>(1-p)/2$ be a constant and let $Y(\tau)$ be a vector field along 
$\gamma$ such that $|Y(\2{\tau})|=1$ and $Y(\tau)$ solves the O.D.E.
$$
\nabla_XY=-\text{Ric}(Y,\cdot)+\frac{b}{\tau}Y\quad\text{ in }(0,\2{\tau}).
\tag 3.6
$$
Then
$$\align
\text{Hess}_{L_p(q,\2{\tau})}(Y(\2{\tau}),Y(\2{\tau}))
\le&-2\2{\tau}^p\text{Ric}(q,\2{\tau})(Y(\2{\tau}),Y(\2{\tau}))
+\frac{2b^2}{(p+2b-1)\2{\tau}^{1-p}}\\
&\qquad +(2p-1)\int_0^{\2{\tau}}\tau^{p-1}\text{Ric}(Y,Y)\,d\tau
-\int_0^{\2{\tau}}\tau^pH(X,Y)\,d\tau\tag 3.7
\endalign
$$
and
$$\align
\Delta L_p(q,\2{\tau})\le&-2\2{\tau}^pR(q,\2{\tau})+\frac{2nb^2}{(p+2b-1)
\2{\tau}^{1-p}}+\frac{2p-1}{\2{\tau}^{2b}}\int_0^{\2{\tau}}\tau^{p+2b-1}R\,
d\tau\\
&\qquad -\frac{1}{\2{\tau}^{2b}}\int_0^{\2{\tau}}\tau^{p+2b}H(X)\,d\tau
\tag 3.8
\endalign
$$
where 
$$\align
H(X,Y)=&-\text{Hess}_R(Y,Y)-2<R(Y,X)Y,X>-4(\nabla_X\text{Ric}(Y,Y)
-\nabla_Y\text{Ric}(Y,X))\\
&\qquad -2\text{Ric}_{\tau}(Y,Y)+2|\text{Ric}(Y,\cdot)|^2
-\frac{1}{\tau}\text{Ric}(Y,Y)
\endalign
$$
and
$$
H(X)=-R_{\tau}-\frac{1}{\tau}R-2<X,\nabla R>+2\text{Ric}(X,X)\tag 3.9
$$
is the Hamilton's expressions for the matrix Harnack inequality and
the trace Harnack inequality respectively (with time equal to $-\tau$). 
\endproclaim
\demo{Proof}
We will use a modification of the argument of section 7 of \cite{P1} to 
prove the lemma. By (3.6),
$$\align
&\frac{d}{d\tau}|Y|^2=2\,\text{Ric}(Y,Y)+2<\nabla_XY,Y>=\frac{2b}{\tau}|Y|^2
\quad\forall 0\le\tau\le\2{\tau}\\
\Rightarrow\quad&|Y(\tau)|^2=\left (\frac{\tau}{\2{\tau}}\right )^{2b}
\quad\forall 0\le\tau\le\2{\tau}\quad\text{and}\quad Y(0)=0.\tag 3.10\\
\endalign
$$
Let
$$\align
I_1=&\int_0^{\2{\tau}}\tau^p\{\text{Hess}_R(Y,Y)+2<R(Y,X)Y,X>+2|\nabla_XY|^2
+2\nabla_X\text{Ric}(Y,Y)\\
&\qquad -4\nabla_Y\text{Ric}(X,Y)\}\,d\tau.\endalign
$$
Then
$$
I_1=-\int_0^{\2{\tau}}\tau^pH(X,Y)\,d\tau -I_2\tag 3.11
$$
where
$$
I_2=\int_0^{\2{\tau}}\tau^p\left\{2\,\text{Ric}_{\tau}(Y,Y)-
2|\text{Ric}(Y,\cdot)|^2+\frac{1}{\tau}\text{Ric}(Y,Y)
+2 \nabla_X\text{Ric}(Y,Y)-2|\nabla_XY|^2\right\}\,d\tau.\tag 3.12
$$
By P.17 of \cite{P1} and (3.6),
$$\align
\frac{d}{d\tau}\text{Ric}(Y,Y)
=&\text{Ric}_{\tau}(Y,Y)+\nabla_X\text{Ric}(Y,Y)+2\,\text{Ric}(\nabla_XY,Y)\\
=&\text{Ric}_{\tau}(Y,Y)+\nabla_X\text{Ric}(Y,Y)-2|\text{Ric}(Y,\cdot)|^2
+\frac{2b}{\tau}\text{Ric}(Y,Y).\tag 3.13
\endalign
$$
By (3.6), (3.10), (3.12), and (3.13),
$$\align
I_2=&\int_0^{\2{\tau}}\tau^p\biggl\{2\frac{d}{d\tau}\text{Ric}(Y,Y)
-2\nabla_X\text{Ric}(Y,Y)+4|\text{Ric}(Y,\cdot)|^2
-\frac{4b}{\tau}\text{Ric}(Y,Y)-2|\text{Ric}(Y,\cdot)|^2\\
&\qquad +\frac{1}{\tau}\text{Ric}(Y,Y)
+2\nabla_X\text{Ric}(Y,Y)-2\left [|\text{Ric}(Y,\cdot)|^2
-\frac{2b}{\tau}\text{Ric}(Y,Y)+\frac{b^2}{\tau^2}|Y|^2\right ]
\biggr \}\,d\tau\\
=&2\2{\tau}^p\text{Ric}(q,\2{\tau})(Y(\2{\tau}),Y(\2{\tau}))
-\frac{2b^2}{(p+2b-1)\2{\tau}^{1-p}}
+(1-2p)\int_0^{\2{\tau}}\tau^{p-1}\text{Ric}(Y,Y)\,d\tau.\tag 3.14
\endalign
$$
By Lemma 1.1,
$$
\delta_YL_p(q,\2{\tau})=2\2{\tau}^p<X,Y>\quad\Rightarrow\quad
\delta_{\nabla_YY}L_p(q,\2{\tau})=2\2{\tau}^p<X,\nabla_YY>.\tag 3.15
$$
By (3.11), (3.14), (3.15), and Lemma 3.1, (3.7) follows.
Let $\{V_i\}_{i=1}^n$ be an orthonormal basis of $T_{\gamma (\2{\tau})}M$
with respect to the metric $g(\gamma (\2{\tau}),\2{\tau})$. For any $i=1,2,
\dots,n$, let $\4{Y}_i$ the solution of (3.6) with $\4{Y}_i(\2{\tau})=V_i$. 
By an argument similar to the proof of (3.10), 
$$
<\4{Y}_i(\tau),\4{Y}_j(\tau)>=\left (\frac{\tau}{\2{\tau}}\right )^{2b}
\delta_{ij}\quad\forall i,j=1,2,\dots,n.
$$
Let $e_i=\4{Y}_i/|\4{Y}_i|$. Then $\4{Y}_i(\tau)=(\tau/\2{\tau})^b
e_i(\tau)$. 
By putting $Y=\4{Y}_i$ in (3.7) and summing over $i=1,2,\dots,n$, by an 
argument similar to 
the derivation of (7.10) of \cite{P1} in \cite{KL}, we get (3.8) and the 
lemma follows.
\enddemo

\proclaim{\bf Lemma 3.3}
Let $(q,\2{\tau})\in\Omega_p$ and let $\gamma$, $X$, be as in Lemma 3.1.
Suppose $Y(\tau)$ is a $\Cal{L}_p$-Jacobi field along $\gamma$ with $Y(0)=0$.
Then
$$
\text{Hess}_{L_p(q,\2{\tau})}(Y(\2{\tau}),Y(\2{\tau}))=2\2{\tau}^p
<\nabla_{Y(\2{\tau})}X(\2{\tau}),Y(\2{\tau})>.
$$
\endproclaim
\demo{Proof}
Let $\alpha :(-\3,\3)\to M$ be a curve in $M$ such that $\alpha (0)=q$,
$\alpha'(0)=Y(\2{\tau})$. Since $\Omega_p(\2{\tau})$ is open, by choosing 
$\3>0$ sufficiently small we may assume without loss of generality that 
$\alpha (-\3,\3)\subset\Omega_p(\2{\tau})$. Then for any $z\in (-\3,\3)$ 
there exists a unique $\Cal{L}_p(\alpha (z),\2{\tau})$-length minimizing 
geodesic $\gamma_z:[0,\2{\tau}]\to M$ which satisfies $\gamma_z(0)=p_0$
and  $\gamma_z(\2{\tau})=\alpha (z)$. Let $f:[0,\2{\tau}]\times (-\3,\3)\to 
M$ be given by $f(\tau,z)=\gamma_z(\tau)$ and let
$$
\2{Y}(\tau)=\frac{\1 f}{\1 z}(\tau,0).
$$
Then $\2{Y}$ is a $\Cal{L}_p$-Jacobi field along $\gamma$ with $\2{Y}(0)=0$
and $\2{Y}(\2{\tau})=Y(\2{\tau})$. By uniqueness of solution of the O.D.E.
for $\Cal{L}_p$-Jacobi field, $\2{Y}(\tau)=Y(\tau)$ for any 
$0\le\tau\le\2{\tau}$. By Lemma 1.1, (1.3) holds. Since
$$
L_p(f(\2{\tau},z),\2{\tau})=\Cal{L}_p(f(\2{\tau},z),f(\cdot,z),\2{\tau})
\quad\forall |z|<\3,
$$
differentiating (1.3) with respect 
to $z$,
$$\align
&\frac{d^2}{dz^2}L_p(f(\2{\tau},z),\2{\tau})\\
=&\frac{d^2}{dz^2}\Cal{L}_p(f(\2{\tau},z),f(\cdot,z),\2{\tau})\\
=&2\2{\tau}^p<\nabla_z\nabla_{\tau}f(\2{\tau},z),\nabla_zf(\2{\tau},z)>
+2\2{\tau}^p<\nabla_{\tau}f(\2{\tau},z),\nabla_z\nabla_zf(\2{\tau},z)>\\
&\qquad +\int_0^{\2{\tau}}\tau^p
<\nabla_z\nabla_zf,\nabla R-\frac{2p}{\tau}\nabla_{\tau}f-2\nabla_{\tau}
\nabla_{\tau}f-4\text{Ric}(\nabla_{\tau}f,\cdot)>\,d\tau\\
&\qquad +\int_0^{\2{\tau}}\tau^p
<\nabla_zf,\nabla_z\{\nabla R-(2p/\tau)\nabla_{\tau}f
-2\nabla_{\tau}\nabla_{\tau}f-4\text{Ric}(\nabla_{\tau}f,\cdot)\}>\,
d\tau\\
=&2\2{\tau}^p<\nabla_z\nabla_{\tau}f(\2{\tau},z),\nabla_zf(\2{\tau},z)>
+2\2{\tau}^p<\nabla_{\tau}f(\2{\tau},s),\nabla_z\nabla_zf(\2{\tau},z)>
+I_1+I_2.\tag 3.16
\endalign
$$
Note that since $\gamma$ is a $\Cal{L}_p$-geodesic, $I_1$ vanishes when
$z=0$. Since $Y(\tau)$ is a $\Cal{L}_p$-Jacobi field along $\gamma$, by
the derivation of the $\Cal{L}_p$-Jacobi equation in the proof of Theorem
2.2, $I_2$ also vanishes when $z=0$. Hence by putting $z=0$ in (3.16), by 
(3.15),
$$\align  
&\delta_Y^2L_p(q,\2{\tau})=\delta_{\nabla_YY}L_p(q,\2{\tau})
+2\2{\tau}^p<\nabla_YX,Y>\\
\Rightarrow\quad&\text{Hess}_{L_p(q,\2{\tau})}(Y(\2{\tau}),Y(\2{\tau}))
=\delta_Y^2L_p(q,\2{\tau})-\delta_{\nabla_YY}L_p(q,\2{\tau})
=2\2{\tau}^p<\nabla_YX,Y>.
\endalign
$$
\enddemo

\proclaim{\bf Lemma 3.4}
Let $\2{\tau}\in (0,t_0)$, $v\in U_p'(\2{\tau})$, and let $b>(1-p)/2$ be
a constant. Then 
$$
\frac{d}{d\tau}\log J_p(v,\2{\tau})\le\frac{b^2n}{(p+2b-1)\2{\tau}}
+\frac{2p-1}{2\2{\tau}^{p+2b}}\int_0^{\2{\tau}}\tau^{p+2b-1}R\,d\tau
-\frac{1}{2\2{\tau}^{p+2b}}\int_0^{\2{\tau}}\tau^{p+2b}H(X)\,d\tau\tag 3.17
$$
where the integration is along the $\Cal{L}_p$-geodesic $\gamma_v(\tau)$ which 
satisfies (1.11) and (1.18), $X(\tau)=\gamma_v'(\tau)$, and $H(X)$ is given 
by (3.9).
\endproclaim
\demo{Proof}
We will use a modification of the proof of a similar result for the case 
$p=b=1/2$ in \cite{P1} to prove the lemma. Let $\gamma_v:[0,\2{\tau}]\to M$ 
be the unique $\Cal{L}_p$-geodesic which satisfies (1.11) and (1.18). Let 
$\{V_i\}_{i=1}^n$ be an orthonormal basis of $T_{\gamma_v(\2{\tau})}M$
with respect to the metric $g(\gamma_v(\2{\tau}),\2{\tau})$. By Theorem 2.8 
for any $i=1,2,\dots,n$, there exists a $\Cal{L}_p$-Jacobi field $Y_i(\tau)$ 
along $\gamma_v$ with $Y_i(0)=0$ and $Y_i(\2{\tau})=V_i$. Then by Lemma 3.3,
$$\align
\frac{d}{d\tau}|Y_i|^2(\2{\tau})=&2\,\text{Ric}(Y_i(\2{\tau}),Y_i(\2{\tau}))
+2<\nabla_{X(\2{\tau})}Y_i(\2{\tau}),Y_i(\2{\tau})>\\
=&2\,\text{Ric}(Y_i(\2{\tau}),Y_i(\2{\tau}))+\frac{1}{\2{\tau}^p}
\text{Hess}_{L_p(q,\2{\tau})}(Y_i(\2{\tau}),Y_i(\2{\tau}))\quad\forall 
i=1,2,\dots,n.\tag 3.18
\endalign
$$
For any $i=1,2,\dots,n$, let $\4{Y}_i(\tau)$ be the solution of (3.6) with
$\4{Y}_i(\2{\tau})=V_i$. Then by Lemma 3.2 and (3.18),
$$
\frac{d}{d\tau}|Y_i|^2(\2{\tau})\le\frac{2b^2}{(p+2b-1)\2{\tau}}
+\frac{2p-1}{\2{\tau}^p}\int_0^{\2{\tau}}\tau^{p-1}\text{Ric}
(\4{Y}_i,\4{Y}_i)\,d\tau
-\frac{1}{\2{\tau}^p}\int_0^{\2{\tau}}\tau^pH(X,\4{Y}_i)\,d\tau.
\tag 3.19
$$
Summing (3.19) over $i=1,2,\dots,n$, similar to the proof of Lemma 3.2 we
have
$$
\sum_{i=1}^n\frac{d}{d\tau}|Y_i|^2(\2{\tau})\le\frac{2nb^2}{(p+2b-1)\2{\tau}}
+\frac{2p-1}{\2{\tau}^{p+2b}}\int_0^{\2{\tau}}\tau^{p+2b-1}R\,d\tau
-\frac{1}{\2{\tau}^{p+2b}}\int_0^{\2{\tau}}\tau^{p+2b}H(X)\,d\tau.\tag 3.20
$$
Now
$$
\frac{d}{d\tau}\log J_p(v,\2{\tau})
=\frac{1}{2|Y_i(\2{\tau})|^2}\sum_{i=1}^n\frac{d}{d\tau}|Y_i|^2(\2{\tau})
=\frac{1}{2}\sum_{i=1}^n\frac{d}{d\tau}|Y_i|^2(\2{\tau}).\tag 3.21
$$
Hence by (3.20) and (3.21) the lemma follows.
\enddemo

By putting $b=(1-p)$ in (3.17) we have

\proclaim{\bf Corollary 3.5}
Let $\2{\tau}\in (0,t_0)$ and $v\in U_p'(\2{\tau})$. Then 
$$
\frac{d}{d\tau}\log J_p(v,\2{\tau})\le\frac{(1-p)n}{\2{\tau}}
+\frac{2p-1}{2\2{\tau}^{2-p}}\int_0^{\2{\tau}}\tau^{1-p}R\,d\tau
-\frac{1}{2\2{\tau}^{2-p}}\int_0^{\2{\tau}}\tau^{2-p}H(X)\,d\tau\tag 3.22
$$
where the integration is along the $\Cal{L}_p$-geodesic $\gamma_v(\tau)$ 
which satisfies (1.11) and (1.18), $X(\tau)=\gamma_v'(\tau)$, and $H(X)$ is 
given by (3.9).
\endproclaim

\proclaim{\bf Lemma 3.6}
Let $q\in M$ and let $\4{\gamma}:[0,\2{s}]\to M$ be a $\4{\Cal{L}}_p$-geodesic
satisfying $\4{\gamma}(0)=p_0$ and $\4{\gamma}(\2{s})=q$. 
Suppose there exists $s_0\in (0,\2{s})$ such that $\4{\gamma}(s_0)$ is 
$\4{\Cal{L}}_p$-conjugate to $p_0$ along $\4{\gamma}$. Then there exists
a vector field $\4{Y}_1$ along $\4{\gamma}$ such that
$$
\delta_{\4{Y}_1}^2\4{\Cal{L}}_p(q,\4{\gamma},\2{s})<0.
$$
\endproclaim
\demo{Proof}
Let $\4{X}(s)=\4{\gamma}'(s)$. 
Since $\4{\gamma}(s_0)$ is $\4{\Cal{L}}_p$-conjugate to $p_0$ along 
$\4{\gamma}$, there exists a $\4{\Cal{L}}_p$-Jacobi field $\4{Y}:[0,s_0]\to M$
along $\4{\gamma}|_{[0,s_0]}$, $\4{Y}\not\equiv 0$, such that $\4{Y}(0)=0$ 
and $\4{Y}(s_0)=0$. Since $\4{Y}\not\equiv 0$, $\nabla_{\4{X}}\4{Y}(s_0)\ne 0$.
Let $W$ be a parallel vector field along $\4{\gamma}$ with respect to the
metric $\4{g}(s_0)=g(s_0^{1/(1-p)})$ such that $W(s_0)=\nabla_{\4{X}}\4{Y}
(s_0)$. Let
$$
\4{Y}_0(s)=\left\{\aligned
&\4{Y}(s)\quad\forall 0\le s\le s_0\\
&0\qquad\,\,\,\forall s_0<s\le\2{s}.\endaligned\right.
$$
Let $h\in (0,\min (s_0,\2{s}-s_0)/2)$ be a constant to be determined later.
We choose $\phi\in C^{\infty}(\Bbb{R})$, $0\le\phi\le 1$ on $[0,\2{s}]$, 
such  that $\phi(s)=0$ for all $|s-s_0|\ge h$ and $\phi (s_0)=1$. 
Let $\4{Y}_1(s)=\4{Y}_0(s)+\lambda\phi W(s)$ where $\lambda\in\Bbb{R}$ is some 
constant to determined later. Let $\4{f}:[0,\2{s}]\times (-\3,\3)$ be a 
variation of $\4{\gamma}$ with respect to $\4{Y}_1$ such that $\4{f}(0,z)
=p_0$ on $(-\3,\3)$ and $\4{f}(s,0)=\4{\gamma}(s)$ for any $0\le s
\le\2{s}$ given by Proposition 2.2 of Chapter 9 of \cite{C}. By the same 
argument as the proof of Lemma 1.2 (1.9) 
holds. Differentiating (1.9) with respect to $z$, by the same argument as 
the proof of (2.6),
$$\align
&\frac{d^2}{dz^2}\4{\Cal{L}}_p(\4{f}(\2{s},z),\4{f}(\cdot,z),\2{s})\\
=&\frac{1}{1-p}\int_0^{\2{s}}<\nabla_z\nabla_z\4{f},s^{\frac{2p}{1-p}}
\nabla\4{R}-2(1-p)^2\nabla_s\nabla_s\4{f}-4(1-p)s^{\frac{p}{1-p}}
\4{\text{Ric}}(\nabla_s\4{f},\cdot)>\,ds\\
&\qquad -2(1-p)\int_0^{\2{s}}<\nabla_z\4{f},\nabla_z\biggl\{
\nabla_s\nabla_s\4{f}
-\frac{1}{2(1-p)^2}s^{\frac{2p}{1-p}}\nabla\4{R}\\
&\qquad +\frac{2}{1-p}
s^{\frac{p}{1-p}}\4{\text{Ric}}(\nabla_s\4{f},\cdot)\biggr\}>\,ds\\
=&I_1(z)-2(1-p)I_2(z)
\tag 3.23
\endalign
$$
where 
$$
I_1(z)=\frac{1}{1-p}\int_0^{\2{s}}<\nabla_z\nabla_z\4{f},s^{\frac{2p}{1-p}}
\nabla\4{R}-2(1-p)^2\nabla_s\nabla_s\4{f}-4(1-p)s^{\frac{p}{1-p}}
\4{\text{Ric}}(\nabla_s\4{f},\cdot)>\,ds
$$
and
$$\align
I_2(z)=&\int_0^{\2{s}}<\nabla_z\4{f},\nabla_s\nabla_s\nabla_z\4{f}
+\4{R}(\nabla_z\4{f},\nabla_s\4{f})\nabla_s\4{f}-\frac{1}{2(1-p)^2}
s^{\frac{2p}{1-p}}\nabla_z(\nabla\4{R})\\
&\qquad +\frac{2}{1-p}
s^{\frac{p}{1-p}}\nabla_z\left (\4{\text{Ric}}(\nabla_s\4{f},\cdot)\right)>
\,ds.
\endalign
$$
Since $\4{\gamma}$ is a $\4{\Cal{L}}_p$-geodesic, 
$$
I_1(0)=0.\tag 3.24
$$
Since $\4{Y}$
is a $\4{\Cal{L}}_p$-Jacobi field on $[0,s_0]$,
$$\align
I_2(0)=&\int_0^{\2{s}}<\4{Y}_0+\lambda\phi W,\nabla_{\4{X}}\nabla_{\4{X}}
(\4{Y}_0+\lambda\phi W)
+\4{R}(\4{Y}_0+\lambda\phi W,\4{X})\4{X}\\
&\qquad -\frac{1}{2(1-p)^2}s^{\frac{2p}{1-p}}\nabla_{\4{Y}_0+\lambda\phi W}
(\nabla\4{R})+\frac{2}{1-p}s^{\frac{p}{1-p}}\nabla_{\4{Y}_0+\lambda\phi W}
(\4{\text{Ric}}(\4{X},\cdot))>\,ds\\
=&\lambda\int_0^{\2{s}}<\4{Y}_0+\lambda\phi W,\nabla_{\4{X}}\nabla_{\4{X}}
(\phi W)+\phi\4{R}(W,\4{X})\4{X}-\frac{1}{2(1-p)^2}s^{\frac{2p}{1-p}}
\phi\nabla_{W}(\nabla\4{R})\\
&\qquad +\frac{2}{1-p}s^{\frac{p}{1-p}}\phi\nabla_{W}
(\4{\text{Ric}}(\4{X},\cdot))>\,ds\\
=&\lambda I_{2,1}+\lambda^2I_{2,2}\tag 3.25
\endalign
$$
where 
$$\align
I_{2,1}=&\int_0^{s_0}<\4{Y},\nabla_{\4{X}}\nabla_{\4{X}}(\phi W)
+\phi\4{R}(W,\4{X})\4{X}-\frac{1}{2(1-p)^2}s^{\frac{2p}{1-p}}
\phi\nabla_{W}(\nabla\4{R})\\
&\qquad +\frac{2}{1-p}s^{\frac{p}{1-p}}\phi\nabla_{W}
(\4{\text{Ric}}(\4{X},\cdot))>\,ds\tag 3.26
\endalign
$$
and
$$\align
I_{2,2}=&\int_0^{\2{s}}<\phi W,\nabla_{\4{X}}\nabla_{\4{X}}(\phi W)
+\phi\4{R}(W,\4{X})\4{X}-\frac{1}{2(1-p)^2}s^{\frac{2p}{1-p}}
\phi\nabla_{W}(\nabla\4{R})\\
&\qquad +\frac{2}{1-p}s^{\frac{p}{1-p}}\phi\nabla_{W}
(\4{\text{Ric}}(\4{X},\cdot))>\,ds.\tag 3.27
\endalign
$$
Now by (1.8),
$$\align
&\int_0^{s_0}<\4{Y},\nabla_{\4{X}}\nabla_{\4{X}}(\phi W)>\,ds\\
=&\int_0^{s_0}\left\{\frac{d}{ds}<\4{Y},\nabla_{\4{X}}(\phi W)>
-<\nabla_{\4{X}}\4{Y},\nabla_{\4{X}}(\phi W)>
-\frac{2}{1-p}s^{\frac{p}{1-p}}\4{\text{Ric}}
(\4{Y},\nabla_{\4{X}}(\phi W))\right\}\,ds\\
=&-\int_0^{s_0}\biggl\{<\nabla_{\4{X}}\4{Y},\nabla_{\4{X}}(\phi W)>
+\frac{2}{1-p}s^{\frac{p}{1-p}}\4{\text{Ric}}(\4{Y},\nabla_{\4{X}}
(\phi W))\biggr\}\,ds\\
=&-\int_0^{s_0}\biggl\{\frac{d}{ds}<\nabla_{\4{X}}\4{Y},\phi W>
-<\nabla_{\4{X}}\nabla_{\4{X}}\4{Y},\phi W>-\frac{2}{1-p}s^{\frac{p}{1-p}}\phi
\4{\text{Ric}}(\nabla_{\4{X}}\4{Y},W)\\
&\qquad +\frac{2}{1-p}
s^{\frac{p}{1-p}}\4{\text{Ric}}(\4{Y},\nabla_{\4{X}}(\phi W))\biggr\}\,ds\\
=&-|\nabla_{\4{X}}\4{Y}(s_0)|^2+\int_0^{s_0}\biggl\{\phi<\nabla_{\4{X}}
\nabla_{\4{X}}\4{Y},W>
+\frac{2}{1-p}s^{\frac{p}{1-p}}\phi\,\4{\text{Ric}}(\nabla_{\4{X}}\4{Y},W)\\
&\qquad -\frac{2}{1-p}s^{\frac{p}{1-p}}\4{\text{Ric}}(\4{Y},\nabla_{\4{X}}
(\phi W))\biggr \}\,ds.\tag 3.28
\endalign
$$
Since
$$\align
&s^{\frac{p}{1-p}}\4{\text{Ric}}(\4{Y},\nabla_{\4{X}}(\phi W))\\
=&\frac{d}{ds}(s^{\frac{p}{1-p}}\4{\text{Ric}}(\4{Y},\phi W))
-\frac{p}{1-p}s^{\frac{2p-1}{1-p}}\4{\text{Ric}}(\4{Y},\phi W)
-s^{\frac{p}{1-p}}\frac{\1}{\1 s}(\4{\text{Ric}})(\4{Y},\phi W)\\
&\qquad -s^{\frac{p}{1-p}}\nabla_{\4{X}}(\4{\text{Ric}})(\4{Y},\phi W)
-s^{\frac{p}{1-p}}\4{\text{Ric}}(\nabla_{\4{X}}\4{Y},\phi W)\tag 3.29
\endalign
$$
By (3.26), (3.28), and (3.29),
$$
I_{2,1}=-|\nabla_{\4{X}}\4{Y}(s_0)|^2+\int_{s_0-h}^{s_0}\phi 
G(W,\4{X},\4{Y})\,ds
\tag 3.30
$$
where
$$\align
&G(W,\4{X},\4{Y})\\
=&<\4{Y},\4{R}(W,\4{X})\4{X}-\frac{1}{2(1-p)^2}s^{\frac{2p}{1-p}}
\nabla_{W}(\nabla\4{R})+\frac{2}{1-p}s^{\frac{p}{1-p}}
\nabla_{W}(\4{\text{Ric}}(\4{X},\cdot))>\\
&\qquad+<\nabla_{\4{X}}\nabla_{\4{X}}\4{Y},W>
+\frac{2}{1-p}s^{\frac{p}{1-p}}\4{\text{Ric}}(\nabla_{\4{X}}\4{Y},W)
+\frac{2}{1-p}\biggl \{\frac{p}{1-p}s^{\frac{2p-1}{1-p}}\4{\text{Ric}}
(\4{Y},W)\\
&\qquad +s^{\frac{p}{1-p}}\frac{\1}{\1 s}(\4{\text{Ric}})(\4{Y}, W)
 +s^{\frac{p}{1-p}}\nabla_{\4{X}}(\4{\text{Ric}})(\4{Y},W)
+s^{\frac{p}{1-p}}\4{\text{Ric}}(\nabla_{\4{X}}\4{Y}, W)\biggr\}.
\endalign
$$
Let $C_1=1+\max_{0\le s\le\2{s}}|G(W,\4{X},\4{Y})|$ and let
$$
h=\frac{|\nabla_{\4{X}}\4{Y}(s_0)|^2}{2C_1}.
$$
Then by (3.30),
$$
I_{2,1}\le -\frac{1}{2}|\nabla_{\4{X}}\4{Y}(s_0)|^2.\tag 3.31
$$
We now choose $\lambda<0$ such that $0>\lambda>-|\nabla_{\4{X}}
\4{Y}(s_0)|^2/[4(1+|I_{2,2}|)]$. Then by putting $z=0$ in (3.23), by 
(3.24), (3.25), and (3.31),
$$
\delta_{\4{Y}_1}^2\4{\Cal{L}}_p(q,\4{\gamma},\2{s})\le |\lambda |(1-p)
(-|\nabla_{\4{X}}\4{Y}(s_0)|^2+2\lambda I_{2,2})\le\lambda (1-p)
|\nabla_{\4{X}}\4{Y}(s_0)|^2/2<0
$$ 
and the lemma follows.
\enddemo

As a consequence of Lemma 3.6 and the equivalence of the $\Cal{L}_p$-geodesic 
and $\4{\Cal{L}}_p$-geodesic by relations (1.4), (1.5), we have

\proclaim{\bf Corollary 3.7}
Let $q\in M$ and let $\gamma :[0,\2{\tau}]\to M$ be a $\Cal{L}_p$-geodesic
satisfying $\gamma(0)=p_0$ and $\gamma (\2{\tau})=q$. 
Suppose there exists $\tau_0\in (0,\2{\tau})$ such that $\gamma (\tau_0)$ is 
$\Cal{L}_p$-conjugate to $p_0$ along $\gamma$. Then there exists
a vector field $Y_1$ along $\gamma$ such that
$$
\delta_{Y_1}^2\Cal{L}_p(q,\gamma,\2{\tau})<0.
$$
\endproclaim

By an argument similar to the proof of Lemma 3.6 and Corollary 3.7
we have

\proclaim{\bf Lemma 3.8}
Let $q\in M$and let 
$\gamma :[0,\2{\tau}]\to M$ be a $\Cal{L}_p$-geodesic satisfying 
$\gamma(0)=p_0$ and $\gamma (\2{\tau})=q$. Suppose there exists $\tau_0
\in (0,\2{\tau})$ such that $q$ is $\Cal{L}_p$-conjugate to $\gamma (\tau_0)$ 
along $\left.\gamma\right |_{[\tau_0,\2{\tau}]}$. Then there exists a vector 
field $Y_1$ along $\gamma$ such that
$$
\delta_{Y_1}^2\Cal{L}_p(q,\gamma,\2{\tau})<0.
$$
\endproclaim

\proclaim{\bf Corollary 3.9}
$$
U_p'(\tau_2)\subset U_p'(\tau_1)\quad\forall 0<\tau_1<\tau_2<t_0.
$$ 
\endproclaim
\demo{Proof}
Let $0<\tau_1<\tau_2<t_0$. Let $v\in U_p'(\tau_2)$. Then $\gamma_v(\tau_2)
\in\Omega_p(\tau_2)$ where $\gamma_v(\cdot)=\gamma (\cdot;v):
[0,\tau_2]\to M$ is the 
$\Cal{L}_p$-geodesic that satisfies (1.18). By the definition of 
$\Omega_p(\tau_2)$, $\gamma (\cdot;v)$ is the unique $\Cal{L}_p
(\gamma_v(\tau_2),\tau_2)$-length minimizing $\Cal{L}_p$-geodesic joining 
$p_0$ and $\gamma_v(\tau_2)$ and $\gamma_v(\tau_2)$ is not 
$\Cal{L}_p$-conjugate to $p_0$ along $\gamma_v$. By an argument similar to 
the proof of Proposition 2.2 of chapter 13 of \cite{C}, $\left.\gamma (\cdot;v)
\right|_{[0,\tau_1]}$ is the unique $\Cal{L}_p(\gamma_v(\tau_1),
\tau_1)$-length minimizing $\Cal{L}_p$-geodesic joining $p_0$ and 
$\gamma_v(\tau_1)$. Suppose $\gamma_v(\tau_1)$ is $\Cal{L}_p$-conjugate to 
$p_0$ along $\left.\gamma_v(\cdot)\right|_{[0,\tau_1]}$. Then by Corollary 
3.7 $\gamma_v$ is not a $\Cal{L}_p(\gamma_v(\tau_2),\tau_2)$-length 
minimizing $\Cal{L}_p$-geodesic.
Contradiction arises. Hence $\gamma_v(\tau_1)$ is not $\Cal{L}_p$-conjugate 
to $p_0$ along $\left.\gamma (\cdot;v)\right|_{[0,\tau_1]}$. Thus $v\in 
U_p'(\tau_1)$ and the lemma follows.
\enddemo

By Lemma 3.8 and an argument similar to the proof of Corollary 3.9 we have

\proclaim{\bf Corollary 3.10}
Let $\2{\tau}\in (0,t_0)$ and $q\in M$. Suppose $\gamma :[0,\2{\tau}]\to M$ 
is the $\Cal{L}_p(q,\2{\tau})$-length minimizing $\Cal{L}_p$-geodesic 
which satisfies $\gamma (0)=p_0$, $\gamma (\2{\tau})=q$, given by Theorem 1.11.
Then $q\in\Omega_{\tau_1,p}^{\gamma (\tau_1)}(\2{\tau})$ for any $0<\tau_1
<\2{\tau}$. 
\endproclaim

\proclaim{\bf Remark 3.11}
By Corollary 3.10 and an argument similar to the proof of Proposition 2.15
of \cite{Ye1}, (3.2), (3.7), (3.8), (3.17), (3.22), etc. in this section holds 
in $M\times (0,t_0)$ in the barrier sense of Perelman \cite{P1}.
\endproclaim

$$
\text{Section 4}
$$

In this section we will prove the monotonicity property of the generalized
reduced volume $\4{V}_p(\tau)$ for $1/2\le p<1$. We first start with a lemma.

\proclaim{\bf Lemma 4.1}
Suppose $M$ has nonnegative curvature operator in $(0,T)$. Then for any 
$\2{\tau}\in (0,t_0)$, $v\in U_p'(\2{\tau})$, there exists a constant 
$C_1=C_1(v,\2{\tau})>0$ such that 
$$
\frac{d}{d\tau}\biggl (\tau^{-(1-p)n}e^{-C_1\tau}J_p(v,\tau)\biggr )\le 0
\quad\forall 0<\tau\le\2{\tau}\tag 4.1
$$
and
$$
\lim_{\tau\to 0^+}\tau^{-(1-p)n}J_p(v,\tau)=(1-p)^{-n}.\tag 4.2
$$
Hence
$$
\tau^{-(1-p)n}e^{-C_1\tau}J_p(v,\tau)\le (1-p)^{-n}\quad\forall 0<\tau
\le\2{\tau}.\tag 4.3
$$
If $M$ also has uniformly bounded scalar curvature on $(0,T)$, then we can 
take
$$
C_1=\biggl (\frac{(2p-1)_+}{2(2-p)}+\frac{t_0}{2(2-p)(t_0-\2{\tau})}
\biggr )\|R\|_{L^{\infty}}.\tag 4.4
$$  
\endproclaim
\demo{Proof}
Let $v\in U_p'(\2{\tau})$ and let $\gamma_v:[0,\2{\tau}]\to M$ be the 
unique $\Cal{L}_p$-geodesic which satisfies (1.18). We extend $\gamma_v$ to
a $\Cal{L}_p$-geodesic on $[0,\2{\tau}+\3)$ for some constant $\3\in (0,t_0
-\2{\tau})$. Let
$$
r_1=\sup_{0\le\tau\le\2{\tau}+\3}d_0(p_0,\gamma_v(\tau))
$$
and let 
$$
C_1=\biggl\{\frac{(2p-1)_+}{2(2-p)}+\frac{\2{\tau}+\3}{2\3(2-p)}\biggr\}
\sup\Sb q\in\2{B_0(p_0,r_1)}\\
0\le\tau\le\2{\tau}+\3\endSb|R(q,\tau)|.
$$
If $M$ also has uniformly bounded scalar curvature on $(0,T)$, by Corollary 
1.5 we can choose $\3=t_0-\2{\tau}$ and let $C_1$ be given by (4.4). 
Since $M$ is complete with respect to the metric $g(\tau)$ for any $\tau\in 
(0,t_0)$, $\2{B_0(p_0,r_1)}\times [0,\2{\tau}]$ is compact. Hence $C_1
<\infty$. Let $H(X)$ be given by (3.9). Since $M$ has nonnegative 
curvature operator in $(0,T)$, as observed by Perelman \cite{P1} by 
Hamilton's Harnack inequality for the solutions of Ricci flow \cite{H4},
$$
H(X(\tau))\ge -\biggl (\frac{1}{\tau}+\frac{1}{\2{\tau}+\3-\tau}\biggr )
R(\gamma (\tau),\tau)\ge -\frac{\2{\tau}+\3}{\tau (\2{\tau}+\3-\tau)}
R(\gamma (\tau),\tau)\quad\forall 0<\tau\le\2{\tau}.\tag 4.5
$$
Since $v\in U_p'(\2{\tau})$, by Corollary 3.9 $v\in U_p'(\tau)$ for 
any $0<\tau\le\2{\tau}$. Hence by Corollary 3.5 and (4.5),
$$\align
&\frac{d}{d\tau}\log J_p(v,\tau)\le\frac{(1-p)n}{\tau}
+\frac{2p-1}{2\tau^{2-p}}\int_0^{\tau}\rho^{1-p}R\,d\rho
-\frac{1}{2\tau^{2-p}}\int_0^{\tau}\rho^{2-p}H(X)\,d\rho\\
&\qquad\qquad\qquad\,\,\le\frac{(1-p)n}{\tau}
+\frac{2p-1}{2\tau^{2-p}}\int_0^{\tau}\rho^{1-p}R\,d\rho
+\frac{\2{\tau}+\3}{2\tau^{2-p}}\int_0^{\tau}
\frac{\rho^{1-p}}{\2{\tau}+\3-\rho}R\,d\rho\\
&\qquad\qquad\qquad\,\,\le\frac{(1-p)n}{\tau}+C_1\quad\forall 
0<\tau\le\2{\tau}\\
\Rightarrow\quad&\frac{d}{d\tau}\log \biggl (\tau^{-(1-p)n}e^{-C_1\tau}
J_p(v,\tau)\biggr )\le 0\quad\forall 0<\tau\le\2{\tau}
\endalign
$$
and (4.1) follows. 
Let $\4{\gamma}_{\4{v}}(s)=\gamma_v(\tau)$ and let $\4{J}_p(v,s)=J_p(v,\tau)$ 
where $\4{v}=v/(1-p)$ and $s=\tau^{1-p}$. Then $\4{\gamma}_{\4{v}}(s)$ 
satisfies (1.10) and (1.13) in $(0,\2{s})$ where $\2{s}=\2{\tau}^{1-p}$. 
We write
$$
\4{\gamma}_{\4{v}}(s)=(\4{\gamma}_{\4{v}}^1(s),\4{\gamma}_{\4{v}}^2(s),
\dots,\4{\gamma}_{\4{v}}^n(s))
$$
and 
$$
v=(v^1,v^2,\dots,v^n)
$$
in the normal coordinate system around $p_0$ with respect to the metric 
$g(p_0,0)$. Differentiating (1.13) with respect to $v^j$, $j=1,2,\dots,n$,
$$\left\{\aligned
&\frac{\1\4{\gamma}_{\4{v}}^i}{\1v^j}(0)=0\qquad\qquad\quad\,\,\,
\forall i,j=1,2,\dots,n\\
&\frac{d}{ds}\biggl (\frac{\1\4{\gamma}_{\4{v}}^i}{\1v^j}\biggr )(0)
=\frac{\delta^i_j}{1-p}\quad\forall i,j=1,2,\dots,n.\endaligned\right.
$$
Hence there exists $s_0\in (0,\2{s})$ and functions $\3^i_j(s)$ such
that $\3^i_j(s)\to 0$ as $s\to 0$ for all $i,j=1,2,\dots,n$ and
$$
\frac{\1\4{\gamma}_{\4{v}}^i}{\1v^j}(s)
=\frac{\delta^i_j+\3^i_j(s)}{1-p}s\quad\forall i,j=1,2,\dots,n.\tag 4.6
$$
Since $\sqrt{g(p_0,0)}=1$ in the normal coordinates around $p_0$, by (4.6)
$$\align
&\4{J}_p(v,s)
=\sqrt{g(p_0,s^{\frac{1}{1-p}})}\,\text{det}\biggl (
\frac{\1\4{\gamma}_{\4{v}}^i}{\1v^j}(s)\biggr )
=\frac{s^n}{(1-p)^n}\sqrt{g(p_0,s^{\frac{1}{1-p}})}\,
\text{det}(\delta^i_j+\3^i_j(s))\\
\Rightarrow\quad&\lim_{s\to 0^+}s^{-n}\4{J}_p(v,s)
=\frac{1}{(1-p)^n}\lim_{s\to 0^+}\text{det}(\delta^i_j+\3^i_j(s))=(1-p)^{-n}
\endalign
$$
and (4.2) follows. By (4.1) and (4.2), (4.3) follows.
\enddemo

\proclaim{\bf Lemma 4.2}
Let $\2{\tau}\in (0,t_0)$ and  $v\in U_p'(\2{\tau})$. Suppose 
$\gamma_v:[0,\2{\tau}]\to M$ is the unique $\Cal{L}_p$-geodesic that satisfies 
(1.18). Then
$$
lim_{\tau\to 0^+}l_p(\gamma_v(\tau),\tau)=|v|^2.\tag 4.7
$$
\endproclaim
\demo{Proof}
Let $\2{s}=\2{\tau}^{1-p}$. Let $r_1>0$, $\4{v}$, and $\4{\gamma}_{\4{v}}$ 
be as in the proof of Lemma 4.1 and let
$$
K_1=\sup\Sb q\in\2{B_0(p_0,r_1)}\\
0\le\tau\le\2{\tau}\endSb (|R|+|\nabla R|+|\text{Ric}|).
$$
Now
$$
l_p(\gamma_v(\tau),\tau)=\frac{1-p}{\tau^{1-p}}\int_0^{\tau}\rho^p R(\gamma_v
(\rho),\rho)\,d\rho+\frac{1-p}{\tau^{1-p}}\int_0^{\tau}\rho^p|\gamma_v'
(\rho)|^2\,d\rho=I_1(\tau)+I_2(\tau)\tag 4.8
$$
where
$$
|I_1(\tau)|\le\frac{1-p}{1+p}K_1\tau^{2p}\to 0\quad\text{ as }\tau\to 0.
\tag 4.9
$$
By the same argument as the proof of Lemma 1.4, there exist constants 
$C_2>0$, $C_3>0$, such that (1.17) holds on $[0,\2{s}]$. Then by (1.17),
$$
\int_0^{\tau}\frac{e^{-C_2\rho}}{\rho^p}
(|v|^2-C_3'e^{C_2\rho}\rho^{1+2p})\,d\rho
\le\int_0^{\tau}\rho^p|\gamma_v'(\rho)|^2\,d\rho
\le\int_0^{\tau}\frac{e^{C_2\rho}}{\rho^p}
(|v|^2+C_3'\rho^{1+2p})\,d\rho\tag 4.10
$$
where $C_3'=(1-p)^2C_3$. Now
$$
\int_0^{\tau}\frac{e^{C_2\rho}}{\rho^p}
(|v|^2+C_3'\rho^{1+2p})\,d\rho
\le e^{C_2\tau}\biggl (\frac{\tau^{1-p}}{1-p}|v|^2+C_3'
\frac{\tau^{2+p}}{2+p}\biggr )\tag 4.11
$$
and
$$
\int_0^{\tau}\frac{e^{-C_2\rho}}{\rho^p}
(|v|^2-C_3'e^{C_2\rho}\rho^{1+2p})\,d\rho
\ge\frac{\tau^{1-p}}{1-p}e^{-C_2\tau}|v|^2-C_3'\frac{\tau^{2+p}}{2+p}
\tag 4.12
$$
By (4.10), (4.11), and (4.12),
$$\align
&e^{-C_2\tau}|v|^2-C_3'\frac{1-p}{2+p}\tau^{1+2p}
\le\frac{1-p}{\tau^{1-p}}\int_0^{\tau}\rho^p|\gamma_v'(\rho)|^2\,d\rho
\le e^{C_2\tau}\biggl (|v|^2+C_3'\frac{1-p}{2+p}\tau^{1+2p}\biggr ).\qquad
\tag 4.13
\endalign
$$
Letting $\tau\to 0$ in (4.13),
$$
\lim_{\tau\to 0^+}I_2(\tau)=|v|^2.\tag 4.14
$$
By (4.8), (4.9), and (4.14), we get (4.7).
\enddemo

\proclaim{\bf Theorem 4.3}
Suppose $M$ has nonnegative curvature operator with respect to the metric 
$g(\tau)$ for any $\tau\in [0,T)$. Suppose $M$ also has uniformly bounded 
scalar curvature on $M\times (0,T)$ when $1/2<p<1$. Let $A_0=0$ if $p=1/2$. 
For any $1/2<p<1$ and $0<c<1$, let
$$
A_0=\biggl (\frac{(2p-1)_+}{2(2-p)}+\frac{1}{2(2-p)c}\biggr )
\|R\|_{L^{\infty}}\tag 4.15
$$
and let $\tau_0$ be given by (0.2). Then 
$$
e^{-A_0\tau_2}\4{V}_p(\tau_2)\le e^{-A_0\tau_1}\4{V}_p(\tau_1)
\le (\sqrt{\pi}/(1-p))^n\quad\forall 0<\tau_1\le\tau_2<\2{\tau}_1,1/2\le p
<1\tag 4.16
$$
where $\2{\tau}_1=(1-c)\tau_0$ if $1/2<p<1$ and $\2{\tau}_1=t_0$ if $p=1/2$.
\endproclaim
\demo{Proof}
Let $p\in [1/2,1)$, $0<\tau_2<\2{\tau}_1$, and $v\in U_p'(\tau_2)$. Let 
$\gamma=\gamma_v:[0,\tau_2]\to M$ be the unique $\Cal{L}_p$-geodesic that 
satisfies (1.18) and let $X(\tau)=\gamma '(\tau)$. By Corollary 3.9 and its 
proof, $v\in U_p'(\tau)$ for any $0<\tau\le\tau_2$ and $\left.\gamma
\right|_{[0,\tau]}$ is the unique $\Cal{L}_p(\gamma_v(\tau),\tau)$-length 
minimizing 
$\Cal{L}_p$-geodesic between $p_0$ and $\gamma_v(\tau)$ for any $0<\tau
\le\tau_2$. Hence $L_p(\gamma (\tau),\tau)=\Cal{L}_p(\gamma (\tau),\gamma,
\tau)$ and
$$
\frac{dL_p}{d\tau}(\gamma (\tau),\tau)=\frac{d}{d\tau}\Cal{L}_p
(\gamma (\tau),\gamma,\tau)
=\tau^p(R(\gamma (\tau),\tau)+|X(\tau)|^2)\quad\forall 0<\tau\le\tau_2.
\tag 4.17
$$
When there is no ambiguity we will write $R$, $X$, $L_p$, and $l_p$ for 
$R(\gamma (\tau),\tau)$, $X(\tau)$, $L_p(\gamma (\tau),\tau)$, and
$l_p(\gamma (\tau),\tau)$. Then
$$\align
\frac{dl_p}{d\tau}(\gamma (\tau),\tau)=&\frac{(1-p)}{\tau^{1-p}}
\frac{dL_p}{d\tau}(\gamma (\tau),\tau)-\frac{(1-p)^2}{\tau^{2-p}}
L_p(\gamma (\tau),\tau)\\
=&\frac{(1-p)}{\tau^{2-p}}[\tau^{p+1}(R+|X|^2)-(1-p)L_p]\quad\forall
0<\tau\le\tau_2.\tag 4.18
\endalign
$$
Now by (1.11),
$$\align
\frac{d}{d\tau}(R+|X|^2)=&R_{\tau}+<X,\nabla R>+2<X,\nabla_XX>
+2\,\text{Ric}(X,X)\\
=&R_{\tau}+<X,\nabla R>+<X,\nabla R-\frac{2p}{\tau}X-4\,\text{Ric}(X,\cdot)>
+2\,\text{Ric}(X,X)\\
=&R_{\tau}+2<X,\nabla R>-2\,\text{Ric}(X,X)-\frac{2p}{\tau}|X|^2\\
=&-H(X)-\frac{2p}{\tau}(R+|X|^2)+\frac{2p-1}{\tau}R\tag 4.19
\endalign
$$
where $H(X)$ is given by (3.9). Hence
$$\align
\frac{d}{d\tau}[\tau^{p+1}(R+|X|^2)]=&\tau^{p+1}\biggl\{\frac{d}{d\tau}
(R+|X|^2)+\frac{(p+1)}{\tau}(R+|X|^2)\biggr\}\\
=&\tau^{p+1}\biggl\{-H(X)-\frac{2p}{\tau}(R+|X|^2)+\frac{2p-1}{\tau}R
+\frac{(p+1)}{\tau}(R+|X|^2)\biggr\}\\
=&-\tau^{p+1}H(X)+(1-p)\tau^p(R+|X|^2)+(2p-1)\tau^pR\quad\forall
0<\tau\le\tau_2.\tag 4.20
\endalign
$$
Since
$$
\lim_{\tau\to 0}\tau^{p+1}|X|^2=0\tag 4.21
$$
by (1.18), integrating (4.20) over $(0,\tau)$,
$$\align
&\tau^{p+1}(R+|X|^2)=-\int_0^{\tau}\rho^{p+1}H(X)\,d\rho +(1-p)L_p+(2p-1)
\int_0^{\tau}\rho^pR\,d\rho\quad\forall 0<\tau\le\tau_2\\
\Rightarrow\quad&(1-p)L_p-\tau^{p+1}(R+|X|^2)
=\int_0^{\tau}\rho^{p+1}H(X)\,d\rho -(2p-1)\int_0^{\tau}\rho^pR\,d\rho
\quad\forall 0<\tau\le\tau_2.\tag 4.22
\endalign
$$
Let 
$$
Z_p(v,\tau)=\tau^{-(1-p)n}e^{-l_p(\gamma_v(\tau),\tau)}e^{-A_0\tau}
J_p(v,\tau).
$$
By Corollary 3.5, (4.18), and (4.22), $\forall 0<\tau\le\tau_2$,
$$\align
\frac{d}{d\tau}\log Z_p(v,\tau)=&\frac{2p-1}{2\tau^{2-p}}\int_0^{\tau}
\rho^{1-p}R\,d\rho
-\frac{1}{2\tau^{2-p}}\int_0^{\tau}\rho^{2-p}H(X)\,d\rho -A_0\\
&\qquad +\frac{(1-p)}{\tau^{2-p}}\biggl\{\int_0^{\tau}\rho^{p+1}H(X)\,d\rho 
-(2p-1)\int_0^{\tau}\rho^pR\,d\rho\biggr\}\\
\le&\frac{2p-1}{2\tau^{2-p}}\int_0^{\tau}\rho^{1-p}R\,d\rho
-A_0-\frac{1}{2\tau^{2-p}}\int_0^{\tau}\rho^{2-p}(1-2(1-p)\rho^{2p-1})
H(X)\,d\rho.\tag 4.23
\endalign
$$
When $p=1/2$, the right hand side of (4.23) is $\le 0$. When $1/2<p<1$, 
by Corollary 1.5 we can extend $\gamma$ to a $\Cal{L}_p$-geodesic on 
$(0,(1-c)\tau_0)$. Since $\tau_2<(1-c)\tau_0$, by the Hamilton's Harnack 
inequality \cite{H4} and an argument similar to the proof of (4.5),
$$
H(X(\tau))\ge -\frac{1}{c\tau}R(\gamma (\gamma (\tau),\tau)
\quad\forall 0<\tau<(1-c)\tau_0\tag 4.24
$$
By (4.24) when $1/2<p<1$, the right hand side of (4.23) is bounded above by
$$
\le\frac{2p-1}{2\tau^{2-p}}\int_0^{\tau}\rho^{1-p}R\,d\rho
-A_0+\frac{1}{2c\tau^{2-p}}\int_0^{\tau}\rho^{1-p}(1-2(1-p)\rho^{2p-1})R
\,d\rho\le 0
$$
for any $0<\tau<(1-c)\tau_0$. Hence by Corollary 3.9 and (4.23),
$$\align
\quad Z_p(v,\tau_2)\le&Z_p(v,\tau_1)\quad\forall 0<\tau_1\le\tau_2,v\in 
U_p'(\tau_2)\\
\Rightarrow\quad \int_{U_p'(\tau_2)}Z_p(v,\tau_2)\,dv\le&\int_{U_p'(\tau_2)}
Z_p(v,\tau_1)\,dv\le\int_{U_p'(\tau_1)}Z_p(v,\tau_1)\,dv\quad\forall 
0<\tau_1\le\tau_2<\2{\tau}_1\\
\Rightarrow\qquad\quad e^{-A_0\tau_2}\4{V}_p(\tau_2)\le&e^{-A_0\tau_1}
\4{V}_p(\tau_1)\quad\forall 0<\tau_1\le\tau_2<\2{\tau}_1.\tag 4.25
\endalign
$$
By (4.25), Lemma 4.1, Lemma 4.2, and the monotone convergence theorem,
$$
e^{-A_0\tau_2}\4{V}_p(\tau_2)\le\int_{U_p'(\tau_2)}\lim_{\tau\to 0}Z_p
(v,\tau)\,dv\le (1-p)^{-n}\int_{\Bbb{R}^n}e^{-|v|^2}dv
=(\sqrt{\pi}/(1-p))^n\tag 4.26
$$
holds for any $0<\tau_2<\tau_0$. By (4.25) and (4.26) we get (4.16) and the 
lemma follows.
\enddemo

$$
\text{Section 5}
$$

In this section we will assume that $(M,\2{g})$ is an ancient solution of
the Ricci flow in $(-\infty,0)$. We will fix a point $(p_0,t_0)\in M\times 
(-\infty,0)$ and let $g$ and $\2{g}$ be related by (0.5). Unless stated
otherwise we will also assume that $M$ has nonnegative curvature operator 
with respect to $g(\tau)$ for any $0\le\tau<\infty$. We will consider the 
$\Cal{L}_p$-length, $L_p$-distance, $\4{V}_p(\tau)$, etc. all with respect 
to this point $(p_0,t_0)$. We will derive various scaling properties of
these geometric quantities in this section. For any $\2{\tau}>0$, let
$$
g^{\2{\tau}}(\tau)=\frac{1}{\2{\tau}}g(\2{\tau}\tau)\quad\forall\tau\ge 0.
$$
and let $R^{\2{\tau}}(q,\tau)$ be the scalar curvature of M at $q$ with 
respect to the metric $g^{\2{\tau}}(\tau)$. We also let 
$\Cal{L}_p^{\2{\tau}}$,
$L_p^{\2{\tau}}$, $l_p^{\2{\tau}}$, $V_p^{\2{\tau}}$ be the corresponding
$\Cal{L}_p$, $L_p$, $l_p$, $V_p$ functions with respect to the metric 
$g^{\2{\tau}}$. Note that
$$
R^{\2{\tau}}(q,\tau)=\2{\tau}R(q,\2{\tau}\tau)\quad\forall q\in M,\tau,
\2{\tau}>0.
$$
For any curve $\gamma$ in $M$ we let $\gamma^{\2{\tau}}$ be the curve in $M$
given by $\gamma^{\2{\tau}}(\rho)=\gamma (\2{\tau}\rho)$.

\proclaim{\bf Lemma 5.1}
Let $\gamma :(\tau_1,\tau_2)\to M$ be a $\Cal{L}_p$-geodesic in $(\tau_1,
\tau_2)$. Then for any $\2{\tau}>0$, $\gamma^{\2{\tau}}$ is a 
$\Cal{L}_p$-geodesic with respect to $g^{\2{\tau}}$ in $(\tau_1/\2{\tau},
\tau_2/\2{\tau})$. If $\tau_1=0$ and $\gamma$ satisfies (1.18) for some 
$v\in T_{p_0}M$, then $\gamma^{\2{\tau}}$ satisfies (1.18) with $v$ being 
replaced by $\2{\tau}^{1-p}v$.
\endproclaim
\demo{Proof}
Let $X(\tau)=X(\gamma(\tau))=\gamma'(\tau)$ and let $X^{\2{\tau}}(\rho)
=X^{\2{\tau}}(\gamma^{\2{\tau}}(\rho))=\gamma^{\2{\tau}}{}'
(\rho)$. Then $X^{\2{\tau}}(\rho)=\2{\tau}X(\2{\tau}\rho)$,
and $\forall \tau_1<\rho<\tau_2$, $Y\in T_{\gamma (\2{\tau}\rho)}M$,
$$\align
&(\nabla_XX)(\2{\tau}\rho)=\frac{D}{d\tau}X(\2{\tau}\rho)
=\frac{1}{\2{\tau}}\frac{D}{d\rho}X(\2{\tau}\rho)
=\frac{1}{\2{\tau}^2}(\nabla_{X^{\2{\tau}}}X^{\2{\tau}})(\rho)\\
\Rightarrow\quad&<\nabla_XX(\2{\tau}\rho),Y>_{g(\2{\tau}\rho)}
=\frac{1}{\2{\tau}}
<\nabla_{X^{\2{\tau}}}X^{\2{\tau}}(\rho),Y>_{g^{\2{\tau}}
(\rho)},\tag 5.1
\endalign
$$
$$\left\{\aligned
&\text{Ric}_{g(\2{\tau}\rho)}(X(\2{\tau}\rho),Y)
=\frac{1}{\2{\tau}}\text{Ric}_{g^{\2{\tau}}(\rho)}(X^{\2{\tau}}(\rho),Y)\\
&<\frac{1}{2\2{\tau}\rho}X(\2{\tau}\rho),Y>_{g(\2{\tau}\rho)}
=\frac{1}{\2{\tau}}<\frac{1}{2\rho}X^{\2{\tau}}(\rho),
Y>_{g^{\2{\tau}}(\rho)}\\
&<\nabla R(\gamma (\2{\tau}\rho),\2{\tau}\rho),Y>_{g(\2{\tau}
\rho)}
=\frac{1}{\2{\tau}}<\nabla R^{\2{\tau}}(\gamma (\2{\tau}\rho),\rho),
Y>_{g^{\2{\tau}}(\rho)}.\endaligned\right.
\tag 5.2
$$
Since $\gamma$ satisfies (1.11) in $(\tau_1,\tau_2)$, by (1.11), (5.1),
and (5.2), $\forall \tau_1<\rho<\tau_2, Y\in T_{\gamma (\2{\tau}\rho)}M$,
$$\align
&<\nabla_{X^{\2{\tau}}}X^{\2{\tau}}-\frac{1}{2}\nabla R^{\2{\tau}}
+\frac{p}{\rho}X^{\2{\tau}}+2\,\text{Ric}_{g^{\2{\tau}}(\rho)}
(X^{\2{\tau}},\cdot),Y>_{g^{\2{\tau}}(\rho)}=0\\
\Rightarrow\quad&\nabla_{X^{\2{\tau}}}X^{\2{\tau}}-\frac{1}{2}\nabla 
R^{\2{\tau}}+\frac{p}{\rho}X^{\2{\tau}}+2\,\text{Ric}_{g^{\2{\tau}}(\rho)}
(X^{\2{\tau}},\cdot)=0\quad\forall \tau_1<\rho<\tau_2.
\endalign
$$ 
If $\tau_1=0$ and $\gamma$ satisfies (1.18) for some $v\in T_{p_0}M$, then
$$
\lim_{\rho\to 0}\rho^pX^{\2{\tau}}(\rho)=\lim_{\rho\to 0}\rho^p\2{\tau}
X(\2{\tau}\rho)=\2{\tau}^{1-p}\lim_{\rho\to 0}(\2{\tau}\rho)^p
X(\2{\tau}\rho)=\2{\tau}^{1-p}v
$$
and the lemma follows.
\enddemo

\proclaim{\bf Lemma 5.2}
For any $q\in M$, $\2{\tau}>0$, $\tau>0$, the following holds.
$$\align
\text{(i) }&L_p^{\2{\tau}}(q,\tau)=\frac{L_p(q,\2{\tau}\tau)}{\2{\tau}^p}
\qquad\qquad\qquad\qquad\qquad\qquad\qquad\qquad\qquad\qquad\qquad\qquad\\
\text{(ii) }&l_p^{\2{\tau}}(q,\tau)=\2{\tau}^{1-2p}l_p(q,\2{\tau}\tau).
\qquad\qquad\qquad\qquad\qquad\qquad\qquad\qquad\qquad\qquad\qquad\qquad
\endalign
$$
\endproclaim
\demo{Proof}
Let $\gamma\in\Cal{F}(q,\2{\tau}\tau)$. Then $\gamma^{\2{\tau}}\in\Cal{F}
(q,\tau)$. Let $X$ and $X^{\2{\tau}}$ be as in the proof of Lemma 5.1. Then
$$\align
\Cal{L}_p^{\2{\tau}}(q,\gamma^{\2{\tau}},\tau)=&\int_0^{\tau}\rho^p(
R^{\2{\tau}}(\gamma^{\2{\tau}}(\rho),\rho)
+|X^{\2{\tau}}(\rho)|^2_{g^{\2{\tau}}(\rho)})\,d\rho\\
=&\int_0^{\tau}\rho^p(\2{\tau}R(\gamma (\2{\tau}\rho),\2{\tau}\rho)
+\2{\tau}|X(\2{\tau}\rho)|^2_{g(\2{\tau}\rho)})\,d\rho\\
=&\frac{1}{\2{\tau}^p}\int_0^{\2{\tau}\tau}z^p(R(\gamma (z),z)
+|X(z)|^2_{g(z)})\,dz\\
=&\frac{1}{\2{\tau}^p}\Cal{L}_p(q,\gamma ,\2{\tau}\tau)\tag 5.3
\endalign
$$
Since
$$
\gamma\in\Cal{F}(q,\2{\tau}\tau)\quad\Leftrightarrow\quad
\gamma^{\2{\tau}}\in\Cal{F}(q,\tau),
$$
by taking infimum in (5.3) over $\gamma\in\Cal{F}(q,\2{\tau}\tau)$,
(i) follows. By (0.6) and (i), (ii) follows.
\enddemo

\proclaim{\bf Lemma 5.3}
There exists a constant $C>0$ such that
$$
|\nabla^{g^{\2{\tau}}(\tau)}l_p^{\2{\tau}}|^2(q,\tau)
+4\frac{(1-p)^2}{\tau^{2(1-2p)}}R^{\2{\tau}}(q,\tau)
\le\frac{C}{\tau^{2(1-p)}}l_p^{\2{\tau}}(q,\tau)\quad\forall \2{\tau}>0,
\tau>0,q\in\Omega_p(\2{\tau}\tau).\tag 5.4
$$ 
\endproclaim
\demo{Proof}
We will use a modification of the proof of this result for the case $p=1/2$
and $\2{\tau}=1$ in \cite{P1} to prove the lemma. We will first prove (5.4) 
for the case
$\2{\tau}=1$. Let $\tau>0$, and $q\in\Omega_p(\tau)$. By Theorem 1.11 
there exists a unique $\Cal{L}_p(q,\tau)$-length minimizing 
$\Cal{L}_p$-geodesic $\gamma$ satisfying $\gamma (0)=p_0$, $\gamma (\tau)=q$. 
Let $X(\rho)=\gamma'(\rho)$ for any $\rho\in [0,\tau]$. Choose $\tau_0>2\tau$.
Then $\tau<\tau_0/2$ and (4.24) holds with $c=1/2$. Hence by (4.22) and 
(4.24),
$$\align
(R(q,\tau)+|X(\tau)|^2)\le&\frac{1}{\tau^{p+1}}\biggl (
C_1L_p(q,\tau)-\int_0^{\tau}\rho^{p+1}H(X)d\rho \biggr )\\
\le&\frac{1}{\tau^{p+1}}\biggl (
C_1L_p(q,\tau)+2\int_0^{\tau}\rho^pRd\rho \biggr )\\
\le&\frac{C_1+2}{\tau^{p+1}}L_p(q,\tau)\tag 5.5
\endalign
$$ 
where $H(X)$ is given by (3.9) and 
$$
C_1=\left\{\aligned
&(1-p)\quad\text{ if }0<p<1/2\\
&p\qquad\quad\,\,\,\text{ if }1/2\le p<1.\endaligned\right.
$$
Then by (0.6), (3.1), (5.5),
$$\align
|\nabla l_p(q,\tau)|^2=&\frac{(1-p)^2}{\tau^{2(1-p)}}
|\nabla L_p(q,\tau)|^2\\
=&\frac{(1-p)^2}{\tau^{2(1-p)}}[4\tau^{2p}(|X(\tau)|^2+R(q,\tau))
-4\tau^{2p}R(q,\tau)]\\
\le&\frac{(1-p)^2}{\tau^{2(1-p)}}\biggl (\frac{4C_1'}{\tau^{1-p}}L_p(q,\tau)
-4\tau^{2p}R(q,\tau)\biggr )\\
\le&\frac{4(1-p)C_1'}{\tau^{2(1-p)}}l_p(q,\tau)
-\frac{4(1-p)^2}{\tau^{2(1-2p)}}R(q,\tau)\quad\forall\tau>0,
q\in\Omega_p(\tau)\tag 5.6
\endalign
$$
where $C_1'=C_1+2$.
Hence by (5.6) and Lemma 5.2 for any $\2{\tau}>0$, $\tau>0$, $q\in
\Omega_p(\2{\tau}\tau)$,
$$\align
|\nabla^{g^{\2{\tau}}(\tau)}l_p^{\2{\tau}}(q,\tau)|^2
=&\2{\tau}\cdot\2{\tau}^{2(1-2p)}|\nabla l_p(q,\2{\tau}\tau)|^2\\
\le&\2{\tau}^{3-4p}\biggl (\frac{4(1-p)C_1'}{(\2{\tau}\tau)^{2(1-p)}}
l_p(q,\2{\tau}\tau)
-\frac{4(1-p)^2}{(\2{\tau}\tau)^{2(1-2p)}}R(q,\2{\tau}
\tau)\biggr )\\
\le&\frac{4(1-p)C_1'}{\tau^{2(1-p)}}
l_p^{\2{\tau}}(q,\tau)
-\frac{4(1-p)^2}{\tau^{2(1-2p)}}R^{\2{\tau}}(q,\tau)
\endalign
$$
and (5.4) follows.
\enddemo

\proclaim{\bf Lemma 5.4}
There exists a constant $C>0$ such that
$$
\left |\frac{\1 l_p^{\2{\tau}}}{\1\tau}(q,\tau)\right |
\le C\frac{l_p^{\2{\tau}}(q,\tau)}{\tau}\quad\forall\2{\tau}>0,\tau>0,q\in 
\Omega_p(\2{\tau}\tau).\tag 5.7
$$ 
\endproclaim
\demo{Proof}
We will use a modification of the proof of this result for the case $p=1/2$
and $\2{\tau}=1$ in \cite{P1} to prove the lemma. We will first prove (5.7) 
for the case
$\2{\tau}=1$. Let $\tau>0$, and $q\in\Omega_p(\tau)$. Let $\gamma$ and $X$
be as in the proof of Lemma 5.3. Then by (3.1) and (5.5),
$$\align
\frac{\1 L_p}{\1\tau}(q,\tau)=&\frac{dL_p}{d\tau}(q,\tau)
-\nabla L_p\cdot X=\tau^p(R+|X|^2)-2\tau^p|X|^2\\
=&2\tau^pR-\tau^p(R+|X|^2)\tag 5.8\\
\Rightarrow\quad\left|\frac{\1 L_p}{\1\tau}(q,\tau)\right|
\le&\tau^p(R+|X|^2)\le\frac{C_1'}{\tau}L_p(q,\tau)\tag 5.9
\endalign
$$
where $C_1'>0$ is as in the proof of Lemma 5.3. By (0.6) and (5.9),
$$
\left|\frac{\1 l_p}{\1\tau}(q,\tau)\right|=\frac{(1-p)}{\tau^{1-p}}\left|
\frac{\1 L_p}{\1\tau}(q,\tau)-\frac{(1-p)}{\tau}L_p(q,\tau)\right|
\le C\frac{l_p(q,\tau)}{\tau}\quad\forall\tau>0,q\in\Omega_p(\tau).
\tag 5.10
$$
For the general case we let $\2{\tau}>0$, $\rho>0$, and $q\in\Omega_p
(\2{\tau}\rho)$. Then by Lemma 5.2 and (5.10),
$$
\left|\frac{\1 l_p^{\2{\tau}}}{\1\rho}(q,\rho)\right|=\left|\frac{\1}{\1\rho}
[\2{\tau}^{1-2p}l_p(q,\2{\tau}\rho)]\right|
=\2{\tau}^{2(1-p)}\left|\frac{\1 l_p}{\1\tau}(q,\2{\tau}\rho)\right|
\le C\2{\tau}^{2(1-p)}\frac{l_p(q,\2{\tau}\rho)}{\2{\tau}\rho}
=C\frac{l_p^{\2{\tau}}(q,\rho)}{\rho}
$$
and the lemma follows.
\enddemo

For any $\tau\ge 0$, $q\in M$, we let
$$
\4{L}_p(q,\tau)=\tau^{1-p}L_p(q,\tau)\tag 5.11
$$
and
$$\left\{\aligned
&G_p(q,\tau)=\4{L}_p(q,\tau)-\frac{n}{2p}\tau\\
&G_p(\tau)=\min_{q\in M}G_p(q,\tau).
\endaligned\right.\tag 5.12
$$
Note that by Lemma 1.7 $G_p(\tau)$ is well-defined. 

\proclaim{\bf Lemma 5.5}
$G_p(\tau)$ is a decreasing function of $\tau\ge 0$. 
\endproclaim 
\demo{Proof}
Let $\tau>0$ and $q\in\Omega_p(\tau)$. Let $\gamma$ and $X$ be as in the
proof of Lemma 5.3. By (4.22) and (5.8),
$$
\frac{\1 L_p}{\1\tau}(q,\tau)=2\tau^pR-\frac{(1-p)}{\tau}L_p(q,\tau)
+\frac{1}{\tau}\int_0^{\tau}\rho^{p+1}
H(X)\,d\rho-\frac{(2p-1)}{\tau}\int_0^{\tau}\rho^pR\,d\rho\tag 5.13
$$
where the integration is along the curve $\gamma$. Putting $b=1/2$ in (3.8),
$$
\Delta L_p(q,\tau)\le -2\tau^pR(q,\tau)+\frac{n}{2p\tau^{1-p}}
+\frac{2p-1}{\tau}\int_0^{\tau}\rho^pR\,
d\rho-\frac{1}{\tau}\int_0^{\tau}\rho^{p+1}H(X)\,d\rho
\tag 5.14
$$
By (5.13) and (5.14),
$$\align
&\frac{\1 L_p}{\1\tau}(q,\tau)+\Delta L_p(q,\tau)\le\frac{n}{2p\tau^{1-p}}
-\frac{(1-p)}{\tau}L_p(q,\tau)\quad\forall\tau>0,q\in\Omega_p(\tau)\\
\Rightarrow\quad&\frac{\1\4{L}_p}{\1\tau}(q,\tau)+\Delta\4{L}_p(q,\tau)
\le\tau^{1-p}\biggl (\frac{\1 L_p}{\1\tau}(q,\tau)+\Delta L_p(q,\tau)
+\frac{(1-p)}{\tau}L_p(q,\tau)\biggr )=\frac{n}{2p}\\
\Rightarrow\quad&\frac{\1 G_p}{\1\tau}(q,\tau)+\Delta G_p(q,\tau)\le 0
\quad\forall\tau>0,q\in\Omega_p(\tau).\tag 5.15
\endalign
$$
By (5.15), Corollary 3.10, Remark 3.11, and an argument similar to the proof of
Proposition 2.15 of \cite{Ye1},
$$
\frac{\1 G_p}{\1\tau}(q,\tau)+\Delta G_p(q,\tau)\le 0\quad\text{ in }M\times 
(0,\infty)\tag 5.16
$$
in the barrier sense of Perelman \cite{P1}. By the same argument as the proof 
of Lemma 3.1 of \cite{Ye1} but with (5.16) replacing (3.1) in the proof there 
the lemma follows. 
\enddemo

\proclaim{\bf Corollary 5.6}
$$
\min_{q\in M}l_p^{\2{\tau}}(q,\tau)\le\frac{n(1-p)}{2p\tau^{1-2p}}\quad
\forall\2{\tau}>0,\tau>0.\tag 5.17
$$
\endproclaim
\demo{Proof}
Let $\2{\tau}>0$ and $\tau>0$. By Lemma 5.5,
$$\align
&0=G_p(0)\ge G_p(\tau)=\min_{q\in M}\biggl (\tau^{1-p}L_p(q,\tau)
-\frac{n}{2p}\tau\biggr )
=\min_{q\in M}\biggl (\frac{\tau^{2(1-p)}}{1-p}l_p(q,\tau)-\frac{n}{2p}\tau
\biggr )\\
\Rightarrow\quad&\min_{q\in M}l_p(q,\tau)\le\frac{n(1-p)}{2p\tau^{1-2p}}.
\tag 5.18
\endalign
$$
By Lemma 5.2 and (5.18), (5.17) follows.
\enddemo

By Corollary 5.6 for any $\2{\tau}>0$ there exists $q(\2{\tau})\in M$ such that
$$
l_p^{\2{\tau}}(q(\2{\tau}),1)=\min_{q\in M}l_p^{\2{\tau}}(q,1)
\le\frac{n(1-p)}{2p}.\tag 5.19
$$

\proclaim{\bf Lemma 5.7}
For any $r_0>0$, $\tau_2>\tau_1>0$, there exists a constant $C_1=C_1(r_0,
\tau_1,\tau_2)>0$ such that 
$$
R^{\2{\tau}}(q,\tau)+l_p^{\2{\tau}}(q,\tau)\le C_1\tag 5.20
$$
holds for any $\tau_1\le\tau\le\tau_2$, $\2{\tau}>0$ and $q\in M$ satisfying
$$
d_{g^{\2{\tau}}(1)}(q(\2{\tau}),q)\le r_0.\tag 5.21
$$
\endproclaim
\demo{Proof}
Let $\tau_1\le\tau\le\tau_2$, $\2{\tau}>0$, and $q\in M$ satisfy (5.21). Let 
$\gamma:[0,d]\to M$ be a minimal normalized geodesic joining $q$ and 
$q(\2{\tau})$ with respect to the metric $g^{\2{\tau}}(1)$ where
$d=d_{g^{\2{\tau}}(1)}(q(\2{\tau}),q)$. Then by Lemma 5.3,
$$\align
|l_p^{\2{\tau}}(q,1)^{\frac{1}{2}}-l_p^{\2{\tau}}(q(\2{\tau}),
1)^{\frac{1}{2}}|
=&\left|\int_0^{d}\frac{\1}{\1\rho}l_p^{\2{\tau}}(\gamma (\rho),
1)^{\frac{1}{2}}\,d\rho\right|\\
=&\left|\int_0^{d}<\nabla^{g^{\2{\tau}}(1)}l_p^{\2{\tau}}
(\gamma (\rho),1)^{\frac{1}{2}},\gamma'(\rho)>_{g^{\2{\tau}}(1)}\,d\rho
\right|\\
\le&\int_0^{d}|\nabla^{g^{\2{\tau}}(1)}(l_p^{\2{\tau}}
(\gamma (\rho),1)^{\frac{1}{2}})|\,d\rho\\
\le&Cd\\
\le&Cr_0\\
\Rightarrow\qquad\qquad\qquad\qquad l_p^{\2{\tau}}(q,1)
\le&(l_p^{\2{\tau}}(q(\2{\tau}),1)^{\frac{1}{2}}+Cr_0)^2
\quad\forall\tau_1\le\tau\le\tau_2.\tag 5.22
\endalign
$$
By Lemma 5.4 there exists a constant $a>0$ such that $\forall\2{\tau}>0,
\tau_1\le\tau\le\tau_2$,
$$
l_p^{\2{\tau}}(q,\tau)\le l_p^{\2{\tau}}(q,1)\biggl (\tau^a
+\frac{1}{\tau^a}\biggr )\tag 5.23
$$
holds for any $q\in\Omega_p(\2{\tau}\tau_2)$. Since $\Omega_p(\2{\tau}\tau_2)$
is dense in $M$ and $l_p^{\2{\tau}}(q,\tau)$ is a continuous function, (5.23) 
holds for any $q\in M$. Hence by (5.19), (5.22), and (5.23),
$$
l_p^{\2{\tau}}(q,\tau)\le (\sqrt{n(1-p)/(2p)}+Cr_0)^2\biggl (\tau_2^a
+\frac{1}{\tau_1^a}\biggr )\quad\forall\tau_1\le\tau\le\tau_2.\tag 5.24
$$
By (5.24) and Lemma 5.3, (5.20) follows.
\enddemo

Since
$$
\frac{\1}{\1\tau}g_{ij}^{\2{\tau}}=2R_{ij}^{\2{\tau}},
$$
by (5.20), Lemma 1.6, Lemma 5.3, Lemma 5.4 and an argument similar to the 
proof of Lemma 5.7 we have the following lemma.

\proclaim{\bf Lemma 5.8}
For any $r_0>0$, $\tau_2>\tau_1>0$, there exists a constant $C>0$ 
such that
$$\left\{\aligned
&|l_p^{\2{\tau}}(q_1,\tau)-l_p^{\2{\tau}}(q_2,\tau)|\le Cd_{g^{\2{\tau}}(1)}
(q_1,q_2)\quad\forall \2{\tau}>0,\tau_1\le\tau\le\tau_2,\text{ and }q_1,q_2
\text{ satisfying (5.21)}\\
&|l_p^{\2{\tau}}(q,\rho_1)-l_p^{\2{\tau}}(q,\rho_2)|\le C|\rho_1-\rho_2|
\qquad\,\,\,\forall\2{\tau}>0, \rho_1,\rho_2\in [\tau_1,\tau_2]\text{ and }
q\text{ satisfying (5.21)}.\endaligned\right.
$$
\endproclaim

By Lemma 5.7, Lemma 5.8, and an argument similar to the sketch of proof of 
Proposition 11.2 of \cite{P1} and a diagonalization argument we have

\proclaim{\bf Corollary 5.9}
Suppose $(M,\2{g})$ has nonnegative curvature operator in $(-\infty,0)$. 
If $(M,\2{g})$ is $\kappa$-noncollapsing on all scales, then there exist a 
sequence $\{q_i\}_{i=1}^{\infty}\subset M$ and a sequence 
$\{\2{\tau}_i\}_{i=1}^{\infty}\subset\Bbb{R}^+$, $\2{\tau}_i\to\infty$ as 
$i\to\infty$, such that $l^{\2{\tau}_i}(q,\tau)$ converges uniformly on
$$
d_{g^{\2{\tau}_i}(1)}(q_i,q)\le r_0,\tau_1\le\tau\le\tau_2
$$ 
as $i\to\infty$ for any $r_0>0$ and $\tau_2>\tau_1>0$.
\endproclaim

For any $\2{\tau}>0$, $\rho>0$, $0<p<1$, let
$$
\4{V}_p^{\2{\tau}}(\rho)=\int_M\rho^{-\frac{(1-p)n}{2}}e^{-l_p^{\2{\tau}}
(q,\rho)}\,dV_{g^{\2{\tau}}(\rho)}(q).
$$

\proclaim{\bf Theorem 5.10}
Suppose $(M,\2{g})$ is an ancient $\kappa$-solution of the Ricci flow.
Let $g$ and $\2{g}$ be related by (0.5) for some constant $t_0<0$. 
Let $\2{\tau}_0>0$ for $1/2<p<1$ and $\2{\tau}_0=0$ for $p=1/2$.
When $1/2<p<1$, suppose also that $(M,g(\tau))$ is compact and satisfies 
(1.21) in $M\times (0,\infty)$ for some constant $c_2>0$. Let $A_0=0$ for 
$p=1/2$ and $A_0$ be given by (4.15) with $c=1$ for $1/2<p<1$. Then for any 
$1/2\le p<1$ there exist constants $A_1\ge 0$, $A_2\ge 0$, such that 
$e^{-W(\2{\tau},\rho)}\4{V}_p^{\2{\tau}}(\rho)$ is a monotone decreasing 
function of $\2{\tau}>\2{\tau}_0$ for any $\rho$ satisfying (0.3) where
$$
W(\2{\tau},\rho)=(A_0\rho+A_1\rho^{2p}+A_2\rho^{2p-3}
e^{2c_2\2{\tau}\rho})\2{\tau}\tag 5.25
$$ 
with $A_1=A_2=0$ for $p=1/2$. Moreover (0.4) holds for any $1/2\le p<1$. 
\endproclaim
\demo{Proof}
Let $\rho$ satisfy (0.3), $\2{\tau}_1>\2{\tau}_0$, and $v\in U_p'(\2{\tau}_1
\rho)$. Let $\gamma_v :[0,\2{\tau}_1\rho]\to M$ be the 
unique $\Cal{L}_p(\gamma_v(\2{\tau}_1\rho),\2{\tau}_1\rho)$-length minimizing 
$\Cal{L}_p$-geodesic given by Theorem 1.11 which satisfies (1.18). 
By Corollary 3.9 and an argument similar to the proof of Theorem 4.3, 
$v\in U_p'(\2{\tau}\rho)$ and $L_p(\gamma (\2{\tau}\rho),\2{\tau}\rho)
=\Cal{L}_p(\gamma (\2{\tau}\rho),\gamma,\2{\tau}\rho)$ for any
$0<\2{\tau}\le\2{\tau}_1$. Let
$$
Z_p^{\2{\tau}}(v,\rho)=\2{\tau}^{-\frac{n}{2}}\rho^{-\frac{(1-p)n}{2}}
e^{-W(\2{\tau},\rho)}e^{-\2{\tau}^{1-2p}l_p(\gamma_v(\2{\tau}\rho),\2{\tau}
\rho)}J_p(v,\2{\tau}\rho)
$$
where $W(\2{\tau},\rho)$ is given by (5.25), $A_1\ge 0$, $A_2\ge 0$, are 
constants to be determined later for $1/2<p<1$ and $A_1=A_2=0$ for $p=1/2$. 
By Lemma 5.2, Corollary 3.5, (4.18) and (4.22), for any $\2{\tau}_0<\2{\tau}
\le\2{\tau}_1$,
$$\align
&\frac{d}{d\2{\tau}}\log Z_p^{\2{\tau}}(v,\rho)\\
\le&-\frac{n}{2\2{\tau}}-\frac{(1-2p)}{\2{\tau}^{2p}}l_p(
\gamma_v(\2{\tau}\rho),
\2{\tau}\rho)-\2{\tau}^{(1-2p)}\rho\frac{d}{d\tau}(l_p
(\gamma_v(\2{\tau}\rho),\2{\tau}\rho))
+\rho\frac{d}{d\tau}\log J_p(v,\2{\tau}\rho)\\
&\qquad -(A_0\rho+A_1\rho^{2p})-2c_2A_2\rho^{2p-2}\2{\tau}
e^{2c_2\2{\tau}\rho}\\
\le&-\frac{n}{2\2{\tau}}-\frac{(1-2p)}{\2{\tau}^{2p}}l_p
+\frac{(1-p)}{\2{\tau}^{1+p}\rho^{1-p}}\biggl (
\int_0^{\2{\tau}\rho}w^{p+1}H(X)\,dw-(2p-1)\int_0^{\2{\tau}\rho}w^pR\,dw
\biggr )\\
&\qquad +\frac{(1-p)n}{\2{\tau}}
+\frac{2p-1}{2\2{\tau}^{2-p}\rho^{1-p}}\int_0^{\2{\tau}\rho}w^{1-p}R\,dw
-\frac{1}{2\2{\tau}^{2-p}\rho^{1-p}}\int_0^{\2{\tau}\rho}w^{2-p}H(X)
\,dw\\
&\qquad -(A_0\rho+A_1\rho^{2p})-2c_2A_2\rho^{2p-2}\2{\tau}
e^{2c_2\2{\tau}\rho}\\
\le&-(2p-1)\frac{n}{2\2{\tau}}+\frac{(2p-1)}{\2{\tau}^{2p}}l_p
+\frac{(2p-1)}{2\2{\tau}^{2-p}\rho^{1-p}}\int_0^{\2{\tau}\rho}w^{1-p}R\,dw
-(A_0\rho+A_1\rho^{2p})\\
&\qquad -2c_2A_2\rho^{2p-2}\2{\tau}e^{2c_2\2{\tau}\rho}
-\frac{1}{2\2{\tau}^{2-p}\rho^{1-p}}
\int_0^{\2{\tau}\rho}\biggl (1-2(1-p)\biggl (\frac{w}{\2{\tau}}
\biggr )^{2p-1}\biggr )w^{2-p}H(X)\,dw\tag 5.26
\endalign
$$
where the integration is along the curve $\gamma_v$. We now divide the
proof into two cases.

\noindent $\underline{\text{\bf Case 1}}$:$p=1/2$.

Then the right hand side of (5.26) is $\le 0$ for any $\2{\tau}_0<\2{\tau}
\le\2{\tau}_1$. 

\noindent $\underline{\text{\bf Case 2}}$:$1/2<p<1$. 

Since $M$ has uniformly bounded Ricci curvature, by Corollary 1.5 we can 
extend $\gamma_v$ to a $\Cal{L}_p$-geodesic on $(0,\infty)$. For any $0<c<1$,
$\tau>0$, choose $\tau_0$ such that $\tau<(1-c)\tau_0$. Then by the Hamilton
Harnack inequality \cite{H4} (4.24) holds. Letting $\tau_0\to\infty$ and 
$c\to 1$ in (4.24),
$$
H(X(\tau))\ge -\frac{1}{\tau}R(\gamma (\gamma (\tau),\tau)
\quad\forall\tau>0.\tag 5.27
$$
By (5.27) the right hand side of (5.26) is bounded above by
$$\align
\le&-(2p-1)\frac{n}{2\2{\tau}}+\frac{(2p-1)}{\2{\tau}^{2p}}l_p
+\frac{(2p-1)}{2\2{\tau}^{2-p}\rho^{1-p}}\int_0^{\2{\tau}\rho}w^{1-p}R\,dw
-(A_0\rho+A_1\rho^{2p})\\
&\qquad -2c_2A_2\rho^{2p-2}\2{\tau}e^{2c_2\2{\tau}\rho}
+\frac{1}{2\2{\tau}^{2-p}\rho^{1-p}}
\int_0^{\2{\tau}\rho}\biggl (1-2(1-p)\biggl (\frac{w}{\2{\tau}}
\biggr )^{2p-1}\biggr )w^{1-p}R\,dw\\
\le&-(2p-1)\frac{n}{2\2{\tau}}
+\frac{(2p-1)}{\2{\tau}^{2p}}l_p-A_1\rho^{2p}-2c_2A_2\rho^{2p-2}\2{\tau}
e^{2c_2\2{\tau}\rho}\tag 5.28
\endalign
$$
Since $M$ is compact and satisfies (1.21), by (0.6) and Lemma 1.8
there exists a constant $C_0>0$ such that
$$
l_p(\gamma_v(\2{\tau}\rho),\2{\tau}\rho)\le C_0\biggl ((\2{\tau}\rho)^{2p}
+\frac{e^{2c_2\2{\tau}\rho}}{(\2{\tau}\rho)^{2-2p}}\biggr ).\tag 5.29
$$
Let $A_1=(2p-1)C_0$ and $A_2=(2p-1)C_0/(2c_2\2{\tau}_0^3)$. Then by (5.28)
and (5.29) the right hand side of (5.26) is $\le 0$ for any 
$\2{\tau}_0<\2{\tau}\le\2{\tau}_1$.

By case 1 and case 2 and an argument similar to the proof of Theorem 4.3
we get that $e^{-W(\2{\tau},\rho)}\4{V}_p^{\2{\tau}}(\rho)$ is a monotone 
decreasing function of $\2{\tau}>\2{\tau}_0$. Hence when $p=1/2$, 
$\4{V}_p^{\2{\tau}}(\rho)$ is a monotone decreasing function of $\2{\tau}>0$.
We now write
$$
Z_p^{\2{\tau}}(v,\rho)=e^{-W(\2{\tau},\rho)}
[(\2{\tau}\rho)^{-(1-p)n}J_p(v,\2{\tau}\rho)][\2{\tau}^{-\frac{n}{2}(2p-1)}
e^{-\2{\tau}^{1-2p}l_p(\gamma_v(\2{\tau}\rho),\2{\tau}\rho)}]
\rho^{\frac{(1-p)n}{2}}.
$$
By Theorem 2.2 and Corollary 3.9 there exist constants $\2{\tau}_0>0$ and 
$r_1>0$ such that $\2{\Cal{B}(0,r_1)}\subset  U_p'(\2{\tau}\rho)$ for all 
$0<\2{\tau}\le\2{\tau}_0$. Since $\2{\Cal{B}(0,r_1)}$ is compact and the 
solution of a $\Cal{L}_p$-geodesic depends continuously on the initial data, 
by an argument similar to the proof of Lemma 4.1, for any $\3>0$, there exists 
$\2{\tau}_1\in (0,\2{\tau}_0)$ and a constant $C_1>0$ such that $\forall 
v\in\2{\Cal{B}(0,r_1)}$, 
$$
(1-p)^{-n}-\3\le (\2{\tau}\rho)^{-(1-p)n}e^{-C_1\2{\tau}\rho}J_p(v,\2{\tau}
\rho)\le (1-p)^{-n}+\3\quad\forall 0<\2{\tau}\le\2{\tau}_1.
$$
Hence $\forall v\in\2{\Cal{B}(0,r_1)}, 0<\2{\tau}\le\2{\tau}_1$, 
$$\left\{\aligned
&Z_p^{\2{\tau}}(v,\rho)\le e^{C_1\2{\tau}\rho -W(\2{\tau},\rho)}
((1-p)^{-n}+\3)\rho^{\frac{(1-p)n}{2}}[\2{\tau}^{-\frac{n}{2}(2p-1)}
e^{-\2{\tau}^{1-2p}l_p(\gamma_v(\2{\tau}\rho),\2{\tau}\rho)}]\\
&Z_p^{\2{\tau}}(v,\rho)\ge e^{C_1\2{\tau}\rho -W(\2{\tau},\rho)}
((1-p)^{-n}-\3)\rho^{\frac{(1-p)n}{2}}[\2{\tau}^{-\frac{n}{2}(2p-1)}
e^{-\2{\tau}^{1-2p}l_p(\gamma_v(\2{\tau}\rho),\2{\tau}\rho)}].
\endaligned\right.\tag 5.30
$$
Let
$$
g(v,\2{\tau},\rho)=\2{\tau}^{-\frac{n}{2}(2p-1)}
e^{-\2{\tau}^{1-2p}l_p(\gamma_v(\2{\tau}\rho),\2{\tau}\rho)}.
$$
By the proof of Lemma 4.2 there exists constants $C_2>0$, $C_3'>0$, $K_1>0$
such that (4.8), (4.9), and (4.13) holds. Let $C_4=(1-p)K_1/(1+p)$ and
$C_5=C_3'(1-p)/(2+p)$. Then by (4.8), (4.9), and (4.13), 
$$\align
g(v,\2{\tau},\rho)\le&\2{\tau}^{-\frac{n}{2}(2p-1)}e^{-\2{\tau}^{1-2p}
[-C_4(\2{\tau}\rho)^{2p}+e^{-C_2\2{\tau}\rho}|v|^2-C_5(\2{\tau}\rho)^{1+2p}]}\\
\le&\2{\tau}^{-\frac{n}{2}(2p-1)}e^{-e^{-C_2\2{\tau}\rho}
(|v|/\2{\tau}^{(2p-1)/2})^2}e^{C_4\2{\tau}\rho^{2p}+C_5\2{\tau}^2\rho^{1+2p}}
\quad\forall v\in\Omega_p(\2{\tau}\rho).
\tag 5.31
\endalign
$$
Similarly 
$$
g(v,\2{\tau},\rho)\ge\2{\tau}^{-\frac{n}{2}(2p-1)}e^{-e^{C_2\2{\tau}\rho}
(|v|/\2{\tau}^{(2p-1)/2})^2}e^{-C_4\2{\tau}\rho^{2p}
-C_5\2{\tau}^2\rho^{1+2p}e^{C_2\2{\tau}\rho}}\quad\forall v\in
\Omega_p(\2{\tau}\rho).\tag 5.32
$$
Hence
$$\align
&\int_{U_p'(\2{\tau}\rho)}g(v,\2{\tau},\rho)\,dv\le e^{C_4\2{\tau}\rho^{2p}
+C_5\2{\tau}^2\rho^{1+2p}}\2{\tau}^{-\frac{n}{2}(2p-1)}\int_{T_{p_0}M}
e^{-e^{-C_2\2{\tau}\rho}(|v|/\2{\tau}^{(2p-1)/2})^2}\,dv\\
&\qquad\qquad\qquad\qquad\,\,\le e^{C_4\2{\tau}\rho^{2p}
+C_5\2{\tau}^2\rho^{1+2p}+(nC_2/2)\2{\tau}\rho}\int_{\Bbb{R}^n}
e^{-|v'|^2}\,dv'\\
&\qquad\qquad\qquad\qquad\,\,\le e^{C_4\2{\tau}\rho^{2p}+C_5\2{\tau}^2
\rho^{1+2p}+(nC_2/2)\2{\tau}\rho}\pi^{\frac{n}{2}}\\
\Rightarrow\quad&\limsup_{\2{\tau}\to 0^+}\int_{U_p'(\2{\tau}\rho)}
g(v,\2{\tau},\rho)\,dv\le\pi^{\frac{n}{2}}.\tag 5.33
\endalign
$$
By (5.30) and (5.33),
$$\align
&\limsup_{\2{\tau}\to 0^+}\4{V}_p^{\2{\tau}}(\rho)\le [(1-p)^{-n}+\3]
\pi^{\frac{n}{2}}\\
\Rightarrow\quad&\limsup_{\2{\tau}\to 0^+}\4{V}_p^{\2{\tau}}(\rho)
\le (1-p)^{-n}\pi^{\frac{n}{2}}\quad\text{ as }\3\to 0.\tag 5.34
\endalign
$$
Similarly by (5.32),
$$
\int_{\Cal{B}(0,r_1)}g(v,\2{\tau},\rho)\,dv\ge e^{-C_4\2{\tau}\rho^{2p}
-C_5\2{\tau}^2\rho^{1+2p}-(nC_2/2)\2{\tau}\rho}\int_{\Cal{B}(0,r_2)}
e^{-|v'|^2}\,dv'.\tag 5.35
$$
where $r_2=\2{\tau}^{-\frac{2p-1}{2}}e^{\frac{C_2\2{\tau}\rho}{2}}r_1$.
Since $r_2\to\infty$ as $\2{\tau}\to 0$, letting $\2{\tau}\to 0$ in (5.35),
$$
\liminf_{\2{\tau}\to 0^+}\int_{\Cal{B}(0,r_1)}
g(v,\2{\tau},\rho)\,dv\ge\pi^{\frac{n}{2}}.\tag 5.36
$$
By (5.30) and (5.36),
$$\align
&\liminf_{\2{\tau}\to 0^+}\4{V}_p(\rho)\ge[(1-p)^{-n}-\3]\pi^{\frac{n}{2}}\\
\Rightarrow\quad&\liminf_{\2{\tau}\to 0^+}\4{V}_p(\rho)\ge (1-p)^{-n}
\pi^{\frac{n}{2}}\quad\text{ as }\3\to 0.\tag 5.37
\endalign
$$
By (5.34) and (5.37), (0.4) follows.
\enddemo

$$
\text{Section 6}
$$

In this section we will prove a conjecture on the reduced distance $l$ and  
the reduced volume $\4{V}(\tau)$ used by Perelman in \cite{P1}. This result 
was used in the proof of Proposition 11.2 of \cite{P1} but no proof was given 
by Perelman in \cite{P1}. We will assume that $(M,\2{g})$ is an ancient 
$\kappa$-solution with $g$ and $\2{g}$ being related by (0.5) for some fixed 
$t_0<0$. For any $\tau>0$, let $\Omega (\tau)=\Omega_{\frac{1}{2}}(\tau)$.
We also fix a point $p_0\in M$ and have $L(q,\tau)$, $\4{V}(\tau)$, 
etc. all defined with respect to the reference point $(p_0,t_0)$.

By an argument similar to the proof of Theorem 6 of \cite{H1} we have

\proclaim{\bf Lemma 6.1}
Let $\tau_2>\tau_1>0$. For any $r_0>0$ there exists a unique solution
$0\le f\in C^{\infty}(\2{B_0(p_0,r_0)}\times [\tau_1,\tau_2])$ of 
$$\left\{\aligned
&f_{\tau}=\Delta f-Rf\qquad\qquad\,\,\,
\text{ in }B_0(p_0,r_0)\times (\tau_1,\tau_2)\\
&f(q,\tau)=\tau^{-\frac{n}{2}}e^{-l(q,\tau)}\quad\,\,\,\text{ on }
\1 B_0(p_0,r_0)\times (\tau_1,\tau_2)\\
&f(q,\tau_1)=\tau_1^{-\frac{n}{2}}e^{-l(q,\tau_1)}\quad\text{ in }B_0(p_0,r_0).
\endaligned\right.\tag 6.1
$$
\endproclaim

We now state and prove a conjecture of Perelman (cf. Proposition 11.2 of 
\cite{P1}). Note that this conjecture was used implicitly by Perelman in his 
proof of Proposition 11.2 but no proof of it was given in \cite{P1}.

\proclaim{\bf Theorem 6.2}
Suppose $\4{V}(\tau_1)=\4{V}(\tau_2)$ for some $\tau_2>\tau_1>0$. Then
$l(q,\tau)\in C^{\infty}(M\times (\tau_1,\tau_2))$ and satisfies
$$
l_{\tau}-\Delta l+|\nabla l|^2-R+\frac{n}{2\tau}=0\tag 6.2
$$
in $M\times (\tau_1,\tau_2)$ in the classical sense.
\endproclaim
\demo{Proof}
Suppose $\4{V}(\tau_1)=\4{V}(\tau_1)$ for some $\tau_2>\tau_1>0$. Let
$r_0>0$ and let $f$ be the solution of (6.1) given by Lemma 6.1. 
Let
$$
Q(\phi)=\int_{\tau_1}^{\tau_2}\int_M\{\nabla l\cdot\nabla\phi 
+(l_{\tau}+|\nabla l|^2
-R+\frac{n}{2\tau})\phi\}\,dV_{g(\tau)}\,d\tau.
$$
By \cite{Ye1},
$$
Q(\phi)=0\tag 6.3
$$
holds for any Lipschitz function $\phi$ on $M\times [\tau_1,\tau_2]$ which 
satisfies
$$\left\{\aligned
&|\phi (q,\tau)|\le Ce^{-l(q,\tau)}\qquad\forall q\in M,\tau_1\le\tau\le\tau_2\\
&|\nabla\phi (q,\tau)|\le Ce^{-l(q,\tau)}\quad\forall q\in M,\tau_1\le\tau
\le\tau_2\endaligned\right.
$$
for some constant $C>0$. Let $h(q,\tau)=\tau^{-\frac{n}{2}}e^{-l(q,\tau)}$.
Then by (6.3), $\forall\phi\in C_0^{\infty}(M\times [\tau_1,\tau_2])$,
$$
Q(h\phi)=0\quad
\Rightarrow\quad\int_{\tau_1}^{\tau_2}\int_M[(h_{\tau}+Rh)\phi+\nabla h
\cdot\nabla\phi]\,dV_{g(\tau)}\,d\tau=0.\tag 6.4
$$
Let $0\le\theta\in C_0^{\infty}(B_0(p_0,r_0))$. For any $\rho\in (\tau_1,
\tau_2]$, let $0\le\psi\in C_0^{\infty}(B_0(p_0,r_0)\times [\tau_1,\rho])$ 
be the solution of 
$$\left\{\aligned
&\psi_{\tau}+\Delta\psi=0\quad\text{ in }B_0(p_0,r_0)\times (\tau_1,\rho)\\
&\psi(q,\tau)=0\qquad\text{ on }\1 B_0(p_0,r_0)\times 
(\tau_1,\rho)\\
&\psi(q,\rho)=\theta\quad\,\,\,\text{ in }B_0(p_0,r_0).\endaligned\right.
\tag 6.5
$$
For any $k\in\Bbb{Z}^+$, let $r_k=(2k-1)r_0/(2k)$ and $\eta_k\in 
C_0^{\infty}(B_0(p_0,r_0))$, $0\le\eta_k\le 1$, such that $\eta_k\equiv 1$ on 
$\2{B_0(p_0,r_k)}$ and $\eta_k\equiv 0$ on $M\setminus B_0(p_0,r_0)$. Then by 
(6.1), (6.4), and (6.5),
$$\align
&\int_{B_0(p_0,r_0)}(f-h)(q,\rho)\theta(q)\,dV_{g(\rho)}\\
=&\int_{B_0(p_0,r_0)}(f-h)(q,\rho)\psi(q,\rho)\,dV_{g(\rho)}
-\int_{B_0(p_0,r_0)}(f-h)(q,\tau_1)\psi(q,\tau_1)\,dV_{g(\tau_1)}\\
=&\int_{\tau_1}^{\rho}\frac{d}{d\tau}\biggl (\int_{B_0(p_0,r_0)}(f-h)\psi
\,dV_{g(\tau)}\biggr )d\tau\\
=&\int_{\tau_1}^{\rho}\int_{B_0(p_0,r_0)}[(f_{\tau}-h_{\tau})\psi +(f-h)
\psi_{\tau}+R(f-h)\psi]\,dV_{g(\tau)}d\tau\\
=&\int_{\tau_1}^{\rho}\int_{B_0(p_0,r_0)}[\psi\Delta f+(f-h)\psi_{\tau}]
\,dV_{g(\tau)}d\tau 
-\int_{\tau_1}^{\rho}\int_{B_0(p_0,r_0)}(h_{\tau}+Rh)\psi\,dV_{g(\tau)}
d\tau.\tag 6.6 
\endalign
$$
Now by (6.4),
$$\align
&\int_{\tau_1}^{\rho}\int_{B_0(p_0,r_0)}(h_{\tau}+Rh)\psi\eta_k
\,dV_{g(\tau)}d\tau\\
=&-\int_{\tau_1}^{\rho}\int_{B_0(p_0,r_0)}\nabla h\cdot\nabla (\psi\eta_k)
\,dV_{g(\tau)}d\tau\\
=&-\int_{\tau_1}^{\rho}\int_{B_0(p_0,r_0)}(\eta_k\nabla h\cdot\nabla\psi
+\psi\nabla h\cdot\nabla\eta_k)\,dV_{g(\tau)}d\tau\\
=&-\int_{\tau_1}^{\rho}\int_{B_0(p_0,r_0)}(\nabla (h\eta_k)\cdot\nabla\psi
+\psi\nabla h\cdot\nabla\eta_k-h\nabla\psi\cdot\nabla\eta_k)\,dV_{g(\tau)}
d\tau\\
=&\int_{\tau_1}^{\rho}\int_{B_0(p_0,r_0)}h\eta_k\Delta\psi\,dV_{g(\tau)}d\tau
-\int_{\tau_1}^{\rho}\int_{B_0(p_0,r_0)}\psi\nabla h\cdot\nabla\eta_k
\,dV_{g(\tau)}d\tau\\
&\qquad +\int_{\tau_1}^{\rho}\int_{B_0(p_0,r_0)}h\nabla\psi\cdot\nabla\eta_k
\,dV_{g(\tau)}d\tau.\tag 6.7
\endalign
$$
Since $\nabla h\in L^{\infty}(B_0(p_0,r_0)\times [\tau_1,\rho])$ by the proof
of Lemma 2.11, letting $k\to\infty$ in (6.7),
$$\align
&\int_{\tau_1}^{\rho}\int_{B_0(p_0,r_0)}(h_{\tau}+R)\psi
\,dV_{g(\tau)}d\tau\\
=&\int_{\tau_1}^{\rho}\int_{B_0(p_0,r_0)}h\Delta\psi\,dV_{g(\tau)}d\tau
-\int_{\tau_1}^{\rho}\int_{\1 B_0(p_0,r_0)}h\frac{\1 \psi}{\1 n}\,d\sigma
d\tau.\tag 6.8
\endalign
$$
By (6.1), (6.6) and (6.8),
$$\align
&\int_{B_0(p_0,r_0)}(f-h)(q,\rho)\theta(q)\,dV_{g(\rho)}\\
=&\int_{\tau_1}^{\rho}\int_{B_0(p_0,r_0)}[\psi\Delta f-h\Delta\psi
+(f-h)\psi_{\tau}]\,dV_{g(\tau)}d\tau
+\int_{\tau_1}^{\rho}\int_{\1 B_0(p_0,r_0)}h\frac{\1 \psi}{\1 n}\,d\sigma
d\tau\\
=&\int_{\tau_1}^{\rho}\int_{B_0(p_0,r_0)}(f-h)(\psi_{\tau}
+\Delta\psi)\,dV_{g(\tau)}d\tau\\
\le&0.\tag 6.9
\endalign
$$
We now choose a sequence of smooth functions $\theta_k\in  C_0^{\infty}
(B_0(p_0,r_0))$, $0\le\theta_k\le 1$, such that $\theta_k\to\text{sign}
(f-h)_+(q,\rho)$ as $k\to\infty$. Putting $\theta=\theta_k$ in (6.9) and
letting $k\to\infty$,
$$
\int_{B_0(p_0,r_0)}(f-h)_+(q,\rho)\,dV_{g(\rho)}\le 0\quad\forall\tau_1
<\rho\le\tau_2\quad\Rightarrow\quad f\le h\quad\text{ in }B_0(p_0,r_0)\times
[\tau_1,\tau_2].\tag 6.10
$$
Similarly we have
$$
\int_{B_0(p_0,r_0)}(h-f)_+(q,\rho)\,dV_{g(\rho)}\le 0\quad\forall\tau_1
<\rho\le\tau_2\quad\Rightarrow\quad h\le f\quad\text{ in }B_0(p_0,r_0)\times
[\tau_1,\tau_2].\tag 6.11
$$
Since $r_0>0$ is arbitrary, by (6.10) and (6.11),
$$\align
&\tau^{-\frac{n}{2}}e^{-l(q,\tau)}=h(q,\tau)=f(q,\tau)\quad\forall q\in M,
\tau_1\le\rho\le\tau_2\tag 6.12\\
\Rightarrow\quad&l(q,\tau)\in C^{\infty}(M\times[\tau_1,\tau_2]).
\endalign
$$
By (6.1) and (6.12), (6.2) follows.
\enddemo

By Theorem 6.2 and an argument similar to that of \cite{KL} and \cite{Ye1} we 
have 

\proclaim{\bf Theorem 6.3} 
Suppose $\4{V}(\tau_1)=\4{V}(\tau_2)$ for some $\tau_2>\tau_1>0$. Then
$l(q,\tau)$ satisfies
$$
2\Delta l-|\nabla l|^2+R+\frac{l-n}{\tau}=0\quad\forall q\in M,\tau_1\le\tau
\le\tau_2
$$
and
$$
R_{ij}(q,\tau)-\frac{1}{2\tau}g_{ij}(q,\tau)+\nabla_i\nabla_jl=0
$$
in $M\times (\tau_1,\tau_2)$.
\endproclaim

\proclaim{\bf Lemma 6.4}
Let $\2{\tau}>0$ and $q\in M$. Suppose $\gamma$ is the 
$\Cal{L}(q,\2{\tau})$-length minimizing $\Cal{L}$-geodesic given by Theorem 
1.11 which satisfies $\gamma (0)=p_0$ and $\gamma (\2{\tau})=q$. Then for any 
$\delta_0\in (0,1)$, there exists a constant $C=C(\delta_0)>0$ such that
$$
\frac{d_{g(\2{\tau})}(\gamma (\delta\2{\tau}),q)^2}{\2{\tau}}\le
C(1+l(q,\2{\tau}))\quad\forall \delta_0\le\delta\le 1.\tag 6.13
$$
\endproclaim
\demo{Proof}
We will use a modification of the proof of Lemma 3.2 of \cite{Ye1} to prove 
the lemma. Let $\2{\tau}>0$. Since $\Omega (\2{\tau})$ is dense in $M$ and
$l(q,\2{\tau})$ is continuous in $q$, it suffices to prove (6.13) for $q\in
\Omega (\2{\tau})$. Let $\delta_0\le\delta\le 1$. Then
$$\align
&d_{\2{\tau}}(\gamma (\delta\2{\tau}),q)\\
=&\int_0^{\2{\tau}}\frac{d}{d\rho}d_{\rho}(\gamma (\delta\rho),\gamma (\rho))
\,d\rho\\
=&\int_0^{\2{\tau}}\biggl (\frac{\1}{\1\rho}d_{\rho}(\gamma (\delta\rho),
\gamma (\rho))+\delta\nabla_Id_{\rho}(\gamma (\delta\rho),\gamma (\rho))\cdot
\gamma '(\delta\rho)+\nabla_{II}d_{\rho}(\gamma (\delta\rho),\gamma(\rho))\cdot
\gamma '(\rho)\biggr )\,d\rho\\
=&I_1+I_2+I_3\tag 6.14
\endalign
$$
where $\nabla_I$ and $\nabla_{II}$ is the gradient with respect to the first
and second argument respectively. Now by (0.6), (3.1), and Lemma 5.3,
$$\align
|\gamma'(\delta\rho)|=&|\nabla l(\gamma (\delta\rho),\delta\rho)|
\le C\biggl (\frac{l(\gamma (\delta\rho),\delta\rho)}{\delta\rho}
\biggr )^{\frac{1}{2}}
\le C(\delta_0\rho)^{-\frac{1}{2}}\biggl (\frac{L(\gamma(\delta\rho),
\delta\rho)}{2\sqrt{\delta\rho}}\biggr )^{\frac{1}{2}}\\
\le&C(\delta_0\rho)^{-\frac{1}{2}}\biggl (\frac{L(\gamma(\2{\tau}),\2{\tau})}
{2\sqrt{\delta\rho}}\biggr )^{\frac{1}{2}}
\le C(\delta_0\rho)^{-\frac{3}{4}}\2{\tau}^{\frac{1}{4}}\sqrt{l(q,\2{\tau})}.
\endalign
$$
Hence
$$
I_2\le C'\2{\tau}^{\frac{1}{2}}\sqrt{l(q,\2{\tau})}.
\tag 6.15
$$
Similarly,
$$
I_3\le C\2{\tau}^{\frac{1}{2}}\sqrt{l(q,\2{\tau})}.\tag 6.16
$$
For any $0<\rho\le\2{\tau}$, let $x(\rho)=\gamma (\delta\rho)$,
$y(\rho)=\gamma (\rho)$, and 
$$
r_0(\rho)=(l(q,\2{\tau})+1)^{-\frac{1}{2}}\rho^{\frac{5}{8}}
\2{\tau}^{-\frac{1}{8}}.
$$
Then for any $x\in B_{\rho}(x(\rho),r_0(\rho))$, by Lemma 5.3,
$$
\sqrt{l(x,\rho)}\le\sqrt{l(x(\rho),\rho)}+\frac{C}{\sqrt{\rho}}r_0(\rho)
\le\sqrt{l(x(\rho),\rho)}+C(l(q,\2{\tau})+1)^{-\frac{1}{2}}\rho^{\frac{1}{8}}
\2{\tau}^{-\frac{1}{8}}.
\tag 6.17
$$
By Lemma 5.4 there exists a constant $C>0$ such that
$$\align
&l(x(\rho),\rho)\le\delta^{-C}l(x(\rho),\delta\rho)\quad\forall 0<\rho
\le\2{\tau}\\
\Rightarrow\quad&l(x(\rho),\rho)\le\delta_0^{-C}\frac{L(q,\2{\tau})}
{2\sqrt{\delta_0\rho}}\le\delta_0^{-C-\frac{1}{2}}\rho^{-\frac{1}{2}}
\2{\tau}^{\frac{1}{2}}l(q,\2{\tau}).\tag 6.18
\endalign
$$
By (6.17) and (6.18),
$$\align
\sqrt{l(x,\rho)}\le&C\rho^{-\frac{1}{4}}\2{\tau}^{\frac{1}{4}}
\sqrt{l(q,\2{\tau})}+C(l(q,\2{\tau})+1)^{-\frac{1}{2}}\rho^{\frac{1}{8}}
\2{\tau}^{-\frac{1}{8}}\quad\forall x\in B_{\rho}(x(\rho),r_0(\rho))\\
\le&C(\rho^{-\frac{1}{4}}\2{\tau}^{\frac{1}{4}}+\rho^{\frac{1}{8}}
\2{\tau}^{-\frac{1}{8}})
\sqrt{l(q,\2{\tau})+1}\quad\forall x\in B_{\rho}(x(\rho),r_0(\rho)).
\tag 6.19
\endalign
$$
By Lemma 5.3 and (6.19),
$$
R(x,\rho)\le C\frac{l(x,\rho)}{\rho}\le C\rho^{-1}
(\rho^{-\frac{1}{4}}\2{\tau}^{\frac{1}{4}}+\rho^{\frac{1}{8}}
\2{\tau}^{-\frac{1}{8}})^2(l(q,\2{\tau})+1)
\quad\forall x\in B_{\rho}(x(\rho),r_0(\rho)).
\tag 6.20
$$
Similarly 
$$
R(x,\rho)\le C\rho^{-1}(\rho^{-\frac{1}{4}}\2{\tau}^{\frac{1}{4}}
+\rho^{\frac{1}{8}}\2{\tau}^{-\frac{1}{8}})^2(l(q,\2{\tau})+1)
\quad\forall x\in B_{\rho}(y(\rho),r_0(\rho)).\tag 6.21
$$
By Lemma 8.3(b) of \cite{P1} and (6.20), (6.21),
$$\align
\frac{\1}{\1\rho}d_{\rho}(\gamma (\delta\rho),\gamma (\rho))
\le&C[\rho^{-1}(\rho^{-\frac{1}{4}}\2{\tau}^{\frac{1}{4}}+\rho^{\frac{1}{8}}
\2{\tau}^{-\frac{1}{8}})^2(l(q,\2{\tau})+1)r_0(\rho)+r_0(\rho)^{-1}]\\
\le&C(\rho^{-\frac{7}{8}}\2{\tau}^{\frac{3}{8}}
+\rho^{-\frac{1}{8}}\2{\tau}^{-\frac{3}{8}}+\rho^{-\frac{1}{2}}
+\rho^{-\frac{5}{8}}\tau^{\frac{1}{8}})\sqrt{l(q,\2{\tau})+1}.
\endalign
$$
Hence
$$
I_1\le C\sqrt{l(q,\2{\tau})+1}
\int_0^{\2{\tau}}(\rho^{-\frac{7}{8}}\2{\tau}^{\frac{3}{8}}
+\rho^{-\frac{1}{8}}\2{\tau}^{-\frac{3}{8}}+\rho^{-\frac{1}{2}}
+\rho^{-\frac{5}{8}}\tau^{\frac{1}{8}})\,d\tau
\le C\2{\tau}^{\frac{1}{2}}\sqrt{l(q,\2{\tau})+1}.
\tag 6.22
$$
By (6.14), (6.15), (6.16), and (6.22), we get (6.13) and the lemma follows.
\enddemo

We now let $\{\2{\tau}_i\}_{i=1}^{\infty}$ be a sequence of positive numbers
such that $\2{\tau}_i\to\infty$ as $i\to\infty$. For any $i\in\Bbb{Z}^+$, 
$\tau>0$, let
$$
\4{V}_i(\tau)=\int_M\tau^{-\frac{n}{2}}e^{-\2{l}_i(q,\tau)}\,
dV_{g_i(\tau)}
$$
where $g_i(\tau)=g(\2{\tau}_i\tau)/\2{\tau}_i$ and $\2{l}_i(q,\tau)$ is the 
$l=l_{\frac{1}{2}}$ function with respect to $g_i(\tau)$. Since by 
Lemma 5.2 $\2{l}_i(q,\tau)=l(q,\2{\tau}_i\tau)$, $\4{V}_i(\tau)
=\4{V}(\2{\tau}_i\tau)$. As observed by Perelman \cite{P1}
there exists a sequence $\{q_i\}_{i=1}^{\infty}\subset M$ and a subsequence 
of $\{\2{\tau}_i\}_{i=1}^{\infty}$ which we may assume without loss of 
generality to be the sequence itself such that the sequence of pointed 
manifold $(M,g_i,q_i)$ will converge in the sense of Hamilton \cite{H6}
$0<\tau<\infty$ to some pointed manifold $(\hat{M},\hat{g},q_0)$ which also 
satisfies the backward Ricci flow as $i\to\infty$. 

That is there exists a sequence of open sets $\hat{U}_i\subset\hat{M}$ with
$q_0\in\hat{U}_i$ for all $i\in\Bbb{Z}^+$ and a sequence of diffeomorphisms
$F_i:\hat{U}_i\to\hat{V}_i$ where $q_i\in\hat{V}_i$ is open in $M$ and
$F_i(q_0)=q_i$ such that for any compact set $K\subset M$ there exists $i_0\in
\Bbb{Z}^+$ such that $K\subset\hat{U}_i$ for all $i\ge i_0$. Moreover 
if $\hat{g}_i=F^{\ast}_i(g_i)$ is the pull-back metric of $g_i$, then the 
metric $\hat{g}_i$ and all its derivatives will converge to $\hat{g}$ and
the corresponding derivatives uniformly on $K\times [a,b]$ as $i\to\infty$
for any $0<a<b<\infty$.

Moreover $\2{l}_i(q_i,1)\le n/2$ for all $i\in\Bbb{Z}^+$ and 
$\2{l}_i(q,\tau)$ converges uniformly on $\2{B_{g_i(1)}(q_i,r)}\times
[a,b]$ to some function $\hat{l}(q,\tau)$ as $i\to\infty$ for any $r>0$ and 
$0<a<b<\infty$.  By Lemma 5.3 and Lemma 5.4 
we may assume without loss of generality that $\2{l}_{i,\tau}(q,\tau)$,
$\nabla\2{l}_{i}$ converge weakly to $\2{l}_{\tau}(q,\tau)$ and
$\nabla\2{l}$ respectively as $i\to\infty$. Then $\2{l}_{\tau}$, 
$|\nabla\2{l}|\in L_{loc}^{\infty}(\hat{M}\times (0,\infty))$. 
Perelman \cite{P1} also proved that $\4{V}_i(\tau)$ decreases and converges 
to some positive constant $\4{V}_0$ which is independent of 
$\tau\in (0,\infty)$ as $i\to\infty$. Let $\hat{R}_{ij}(q,\tau)$ and 
$\hat{R}(q,\tau)$ be the Ricci curvature and scalar curvature of $\hat{M}$ 
with respect to the metric $\hat{g}(q,\tau)$. By an argument similar to the 
proof of Theorem 6.2 and Theorem 6.3 we have

\proclaim{\bf Theorem 6.5}
$\2{l}(q,\tau)\in C^{\infty}(M\times (0,\infty))$ and $\2{l}(q,\tau)$ satisfies
$$
\2{l}_{\tau}-\Delta\2{l}+|\nabla\2{l}|^2-\hat{R}+\frac{n}{2\tau}=0
$$
in $\hat{M}\times (0,\infty)$.
\endproclaim

\proclaim{\bf Theorem 6.6}
Let $\hat{q}\in\hat{M}$. Let $\2{q}_i\in M$ be such that $\hat{q}
=\lim_{i\to\infty}F_i^{-1}(\2{q}_i)$. For each $\rho\ge 1$, $i\in\Bbb{Z}^+$, 
let $\gamma_i:[0,\2{\tau}_i\rho]\to M$ be the 
$\Cal{L}(\2{q}_i,\2{\tau}_i\rho)$-length minimizing $\Cal{L}$-geodesic given 
by Theorem 1.11. Let 
$\gamma_i^{\2{\tau}_i}(w)=\gamma_i(\2{\tau}_i w)$, $0\le w\le\rho$. Then 
there exists a $\Cal{L}$-geodesic $\hat{\gamma}:(0,\rho]
\to\hat{M}$ with $\hat{\gamma}(\rho)=\hat{q}$ which is a 
$\Cal{L}_{\rho_0,\frac{1}{2}}^{\hat{\gamma}(\rho_0)}(\hat{q},\rho)$-length
minimizing $\Cal{L}$-geodesic on $[\rho_0,\rho]$ for any $\rho_0\in 
(0,\rho)$ such 
that for any $\rho_0\in (0,\rho)$ $\gamma_i^{\2{\tau}_i}(w)$ will converge 
uniformly on $\rho_0\le w\le\rho$ to a $\Cal{L}$-geodesic of $\hat{M}$ with 
$\hat{\gamma}(\rho)=\hat{q}$ as $i\to\infty$.
\endproclaim
\demo{Proof}
Let $\hat{q}\in\hat{M}$ and let $\rho\ge 1$. We choose $\2{q}_i\in M$ such 
that $\hat{q}=\lim_{i\to\infty}F_i^{-1}(\2{q}_i)$. Let $b>a>0$. Since 
$d_{g_i(w)}
(q_i,\2{q}_i)$ converges uniformly to $d_{\hat{g}(w)}(q_0,\hat{q})$ on 
$a\le w\le b$ as $i\to\infty$, there exists a constant $C_1>0$ such that
$$
d_{g_i(\rho)}(q_i,\2{q}_i)\le C_1\quad\forall i\in\Bbb{Z}^+.\tag 6.23
$$
Since $\2{l}_i$ converges to $\2{l}$ uniformly on $\2{B_{g_i(1)}(q_i,r)}
\times [a,b]$ as $i\to\infty$ for any $r>0$, $b>a>0$, there exists a 
constant $C_2>0$ such that 
$$
\2{l}_i(\2{q}_i,\rho)\le\2{l}(\hat{q},\rho)+ C_2\quad\forall i\in\Bbb{Z}^+.
\tag 6.24
$$
Let $\delta_0\in (0,1)$. By Lemma 5.1 and Lemma 6.4 $\gamma_i^{\2{\tau}_i}:
[0,1]\to M$ is a minimizing $\Cal{L}$-geodesic with $\gamma_i^{\2{\tau}_i}
(1)=\2{q}_i$ and
$$\align
&\frac{d_{g(\2{\tau}_i\rho)}(\gamma_i(\delta\2{\tau}_i\rho),\2{q}_i)^2}
{\2{\tau}_i\rho}\le C(1+l(\2{q}_i,\2{\tau}_i\rho))
\quad\forall \delta_0\le\delta\le 1,i\in\Bbb{Z}^+\\
\Rightarrow\quad&\frac{d_{g_i(\rho)}(\gamma_i^{\2{\tau}_i}
(\delta\rho),\2{q}_i)^2}
{\rho}\le C(1+\2{l}_i(\2{q}_i,\rho))\quad\forall \delta_0\le\delta\le 1,
i\in\Bbb{Z}^+.\tag 6.25
\endalign
$$
By (6.24) and (6.25) there exists a constant $C_3>0$ such that
$$
\frac{d_{g_i(\rho)}(\gamma_i^{\2{\tau}_i}(\delta\rho),\2{q}_i)^2}{\rho}
\le C_3\quad\forall \delta_0\le\delta\le 1,i\in\Bbb{Z}^+.\tag 6.26
$$
By (6.23) and (6.26),
$$
d_{g_i(\rho)}(\gamma_i^{\2{\tau}_i}(\delta\rho),q_i)\le C_1+\sqrt{C_3\rho}
\quad\forall \delta_0\le\delta\le 1,i\in\Bbb{Z}^+.
$$
Hence by the Hamilton compactness theorem \cite{H6} and Lemma 5.1
$\gamma_i^{\2{\tau}_i}$ will converge uniformly on $\rho_0\le\tau\le\rho$ 
to a $\Cal{L}$-geodesic of $\hat{M}$ with $\hat{\gamma}(\rho)=\hat{q}$ as 
$i\to\infty$ for any $\rho_0\in (0,\rho)$. Since $\delta_0$ is arbtiary, 
$\lim_{i\to\infty}
\gamma_i^{\2{\tau}_i}(\delta\rho)$ exists for any $\delta\in (0,1)$. For any
$0<\tau\le\rho$, let
$$
\hat{\gamma}(\tau)=\lim_{i\to\infty}\gamma_i^{\2{\tau}_i}(\tau).
$$
Then $\hat{\gamma}:(0,\rho)\to\hat{M}$ is a $\Cal{L}$-geodesic of $\hat{M}$
with $\hat{\gamma}(\rho)=\hat{q}$. Since each $\gamma_i^{\2{\tau}_i}$ is a 
$\Cal{L}(\2{q}_i,\rho)$-length minimizing $\Cal{L}$-geodesic,
$\left.\hat{\gamma}\right|_{[\rho_0,\rho]}$ is a $\Cal{L}_{\rho_0,
\frac{1}{2}}^{\hat{\gamma}(\rho_0)}
(\hat{q},\rho)$-length minimizing $\Cal{L}$-geodesic on $[\rho_0,\rho]$ for 
any $\rho_0\in (0,\rho)$ and the theorem follows.
\enddemo

\enddocument